\documentclass[11pt]{amsart}
\usepackage{amsmath}
\usepackage{amsxtra}
\usepackage{amscd}
\usepackage{amsthm}
\usepackage{amsfonts}
\usepackage{amssymb}
\usepackage{eucal}

\newtheorem{thm}{Theorem}
\newtheorem{conj}{Conjecture}
\newtheorem{cor}{Corollary}
\newtheorem{lem}{Lemma}
\newtheorem{prop}{Proposition}

\theoremstyle{remark}

\theoremstyle{definition}

\numberwithin{equation}{section}

\newcommand{\nc}{\newcommand}

\def\bb{{\mathfrak b}}

\def\ggg{{\mathfrak g}} 
\def\hh{{\mathfrak h}}

\def\mm{{\mathfrak m}} 
\def\nn{{\mathfrak n}}

\def\zz{{\mathfrak z}}

\def\D{{\mc D}} 
 
\def\F{{\mc F}}

\def\K{{\mc K}}

\def\O{{\mc O}} 

\def\T{{\mc T}}

\def\({\left(}
\def\){\right)}
\def\[{\left[}
\def\]{\right]}
\def\<{\left\langle}
\def\>{\right\rangle}

\def\tensor{\otimes}

\def\liesl{{\mathfrak{sl}}}

\def\Gr{\mathop{\rm Gr}\nolimits}
\def\Aut{\mathop{\rm Aut}\nolimits}

\def\Gal{\mathop{\rm Gal}\nolimits}
\def\mod{\on{-mod}}

\def\ge{\geqslant}

\def\qed{\hfill$\Box$\smallskip}
\def\proof{{\bf Proof.\ \ }}

\numberwithin{equation}{section}

\nc{\Ref}[1]{{$($\ref{#1}$)$}}

\nc{\thmref}[1]{Theorem~\ref{#1}}
\nc{\secref}[1]{Section~\ref{#1}}
\nc{\lemref}[1]{Lemma~\ref{#1}}
\nc{\propref}[1]{Proposition~\ref{#1}}
\nc{\corref}[1]{Corollary~\ref{#1}}
\nc{\remref}[1]{Remark~\ref{#1}}
\nc{\conjref}[1]{Conjecture~\ref{#1}}
\nc{\on}{\operatorname}
\nc{\ch}{\mbox{ch}}
\nc{\Z}{{\mathbb Z}}
\nc{\C}{{\mathbb C}}
\nc{\pone}{{\mathbb P}^1}
\nc{\pa}{\partial}
\nc{\arr}{\rightarrow}
\nc{\larr}{\longrightarrow}
\nc{\al}{\alpha}
\nc{\ri}{\rangle}
\nc{\lef}{\langle}
\nc{\la}{\lambda}
\nc{\ep}{\epsilon}
\nc{\su}{\widehat{{\mathfrak s}{\mathfrak l}}_2}
\nc{\sw}{{\mathfrak s}{\mathfrak l}}
\nc{\g}{{\mathfrak g}}
\nc{\h}{{\mathfrak h}}
\nc{\n}{{\mathfrak n}}
\nc{\nhat}{\widehat{\n}}
\nc{\De}{\Delta}
\nc{\gt}{\widetilde{\g}}
\nc{\Ga}{\Gamma}
\nc{\one}{{\mathbf 1}}
\nc{\z}{{\mathfrak Z}}
\nc{\La}{\Lambda}
\nc{\wt}{\widetilde}
\nc{\wh}{\widehat}
\nc{\cri}{_{\kappa_c}}
\nc{\hvee}{h^\vee}
\nc{\sun}{\widehat{\sw}_N}
\nc{\si}{\sigma}
\nc{\el}{\ell}
\nc{\bi}{\bibitem}
\nc{\om}{\omega}
\nc{\ol}{\overline}
\nc{\ds}{\displaystyle}
\nc{\dzz}{\frac{dz}{z}}
\nc{\Res}{\on{Res}}
\nc{\mc}{\mathcal}
\nc{\Cal}{\mathcal}
\nc{\ot}{\otimes}
\nc{\ga}{\gamma}
\nc{\scs}{\scriptstyle}
\nc{\us}{\underset}
\nc{\opl}{\oplus}
\nc{\beq}{\begin{equation}}
\nc{\Fq}{{\mathbb F}_q}
\nc{\Mq}{{\mathcal M}}
\nc{\Rep}{\on{Rep}}
\nc{\sssec}{\subsubsection}
\nc{\ssec}{\subsection}
\nc{\lan}{\langle}
\nc{\ran}{\rangle}

\nc{\Vect}{\on{Vect}}
\nc{\ghat}{\widehat{\g}}
\nc{\Tloc}{\T^\g_{\on{loc}}}
\nc{\vac}{|0\ran}
\nc{\Wick}{{\mb :}}
\nc{\mb}{\mathbf}
\nc{\delz}{\partial_z}
\nc{\cali}{\mathcal}
\nc{\li}{\mathfrak l}
\nc{\lt}{\widetilde{\li}}
\nc{\astar}{a^*}
\nc{\cA}{{\mc A}}
\nc{\ka}{\kappa}

\nc{\OO}{{\mc O}}
\nc{\AutO}{\on{Aut} \O}
\nc{\AutOO}{\on{Aut}_+ \O}
\nc{\DerO}{\on{Der} \O}
\nc{\DerpO}{\on{Der}_+ \O}
\nc{\Au}{{\mc A}ut}
\nc{\mf}{\mathfrak}
\nc{\hhat}{\wh{\h}}

\nc{\gr}{\on{gr}}
\nc{\Spe}{\on{Spec}}
\nc{\rv}{\crho}
\nc{\can}{\on{can}}
\nc{\CC}{{\mc C}}
\nc{\Op}{\on{Op}_G(D)}
\nc{\MOp}{\on{MOp}_G(D)}
\nc{\Db}{D}
\nc{\ww}{w}

\nc{\ConD}{\on{Conn}(\Omega^{\crho})_{\Db}}
\nc{\ConDL}{\on{Conn}(\Omega^{\rho})_{\Db}}
\nc{\ConDtL}{\on{Conn}(\Omega^{\rho})_{\Db^\times}}
\nc{\OpD}{\on{Op}_G(\Db)}
\nc{\crho}{\check{\rho}}
\nc{\chal}{\check{\al}}
\nc{\cchi}{\check{\chi}}

\nc{\bone}{{\mb 1}}

\nc{\Q}{{\mathbb Q}}

\nc{\Qp}{{\mathbb Q}_p}

\nc{\AD}{{\mathbb A}}

\nc{\bs}{\backslash}

\nc{\Loc}{\on{Loc}}

\nc{\Ql}{{\mathbb Q}_\ell}

\nc{\ab}{{\mathfrak a}}

\nc{\ppart}{(\!(t)\!)}

\nc{\Vir}{\on{Vir}}

\nc{\Ll}{{\mc L}}

\nc{\Proj}{{\mc P}roj}

\nc{\zpart}{(\!(z)\!)}
\nc{\wpart}{(\!(w)\!)}
\nc{\zwpart}{(\!(z-w)\!)}
\nc{\partzw}{(\!(z,w)\!)}

\nc{\cmu}{\check\mu}

\nc{\abs}{\on{abs}}

\nc{\Pp}{{\mathbb C}{\mathbb P}^1}

\nc{\vect}{|-2\rangle}

\nc{\cla}{\check\lambda}

\nc{\com}{\check\omega}

\nc{\CF}{{\mc F}}
\nc{\CP}{{\mc P}}
\nc{\CH}{{\mc H}}
\nc{\CM}{{\mc M}}
\nc{\CV}{{\mc V}}
\nc{\CL}{{\mc L}}

\nc{\BV}{{\mathbb V}}

\nc{\sF}{{\mathsf{F}}}
\nc{\sG}{{\mathsf{G}}}

\nc{\Hom}{\on{Hom}}
\nc{\Ext}{\on{Ext}}

\nc{\Hecke}{{\on{Hecke}}}
\nc{\Isom}{\on{Isom}}

\nc{\nMOp}{\on{MOp}^{\nilp}_{^L G}}

\nc{\QCoh}{\on{QCoh}}

\nc{\nilp}{{\on{nilp}}}

\nc{\BM}{{\mathbb M}}

\nc{\reg}{{\on{reg}}}

\nc{\nOp}{\on{Op}^{\nilp}}
\nc{\nOpcla}{\on{Op}^{\nilp,\cla}}
\nc{\nOpla}{\on{Op}^{\nilp,\la}}

\nc{\nMOpw}{\on{MOp}^{\nilp,w_0}_{^L G}}

\nc{\Bun}{\on{Bun}}
\nc{\out}{{\on{out}}}

\nc{\gmod}{\ghat_{\ka_c}\mod}
\nc{\gmodx}{\ghat_{\ka_c,x}\mod}
\nc{\gmody}{\ghat_{\ka_c,y}\mod}

\nc{\opmod}{\zz\mod}

\nc{\ord}{\on{ord}}
\nc{\bsl}{\backslash}
\nc{\Conn}{\on{Conn}(\Omega^{\crho})}
\nc{\Ind}{\on{Ind}}
\nc{\Con}{\on{Conn}(\Omega^{-\rho})}
\nc{\Co}{\on{Conn}(\Omega^{\rho})}
\nc{\crit}{{\ka_c}}

\begin{document}

\title{Ramifications of the geometric Langlands Program}

\author{Edward Frenkel}\thanks{Supported by the DARPA grant
HR0011-04-1-0031 and by the NSF grant DMS-0303529.}

\address{Department of Mathematics, University of California,
Berkeley, CA 94720, USA}

\date{November 2006. Based on the lectures given by the author at the
CIME Summer School "Representation Theory and Complex Analysis",
Venice, June 2004. To appear in the Proceedings of the School
published in Lecture Notes in Mathematics, Springer Verlag.}

\maketitle

\tableofcontents

\section*{Introduction}

The Langlands Program, conceived as a bridge between Number Theory and
Automorphic Representations \cite{Langlands}, has recently expanded
into such areas as Geometry and Quantum Field Theory and exposed a
myriad of unexpected connections and dualities between seemingly
unrelated disciplines. There is something deeply mysterious in the
ways the Langlands dualities manifest themselves and this is what
makes their study so captivating.

In this review we will focus on the geometric Langlands correspondence
for complex algebraic curves, which is a particular brand of the
general theory. Its origins and the connections with the classical
Langlands correspondence are discussed in detail elsewhere (see, in
particular, the reviews \cite{F:bull,F:rev}), and we will not try to
repeat this here. The general framework is the following: let $X$ be a
smooth projective curve over $\C$ and $G$ be a simple Lie group over
$\C$. Denote by $^L G$ the Langlands dual group of $G$ (we recall this
notion in \secref{dual grp}). Suppose that we are given a principal
$^L G$-bundle $\F$ on $X$ equipped with a flat connection. This is
equivalent to $\F$ being a holomorphic principal $^L G$-bundle
equipped with a holomorphic connection $\nabla$ (which is
automatically flat as the complex dimension of $X$ is equal to
one). The pair $(\F,\nabla)$ may also be thought of as a $^L G$-local
system on $X$, or as a homomorphism $\pi_1(X) \to {}^L G$
(corresponding to a base point in $X$ and a trivialization of the
fiber of $\F$ at this point).

The global Langlands correspondence is supposed to assign to $E =
(\F,\nabla)$ an object $\on{Aut}_E$, called {\bf Hecke eigensheaf
with eigenvalue} $E$, on the moduli stack $\Bun_G$ of holomorphic
$G$-bundles on $X$:
$$
\boxed{\begin{matrix} \text{holomorphic $^L G$-bundles} \\
\text{with connection on } X \end{matrix}} \quad \longrightarrow \quad
\boxed{\text{Hecke eigensheaves on } \Bun_G}
$$
$$
E \mapsto \on{Aut}_E
$$
(see, e.g., \cite{F:rev}, Sect. 6.1, for the definition of Hecke
eigensheaves). It is expected that there is a unique irreducible Hecke
eigensheaf $\on{Aut}_E$ (up to isomorphism) if $E$ is sufficiently
generic.

The Hecke eigensheaves $\Aut_E$ have been constructed, and the
Langlands correspondence proved, in \cite{FGV1,Ga:vanish} for $G=GL_n$
and an arbitrary irreducible $GL_n$-local system, and in \cite{BD} for
an arbitrary simple Lie group $G$ and those $^L G$-local systems which
admit the structure of a $^L G$-oper (which is recalled below).

Recently, A. Kapustin and E. Witten \cite{KW} have related the
geometric Langlands correspondence to the $S$-duality of
supersymmetric four-dimensional Yang-Mills theories, bringing into the
realm of the Langlands correspondence new ideas and insights from
quantum physics.

\medskip

So far, we have considered the {\bf unramified} $^L G$-local
systems. In other words, the corresponding flat connection has no
poles. But what should happen if we allow the connection to be
singular at finitely many points of $X$?

This {\bf ramified geometric Langlands correspondence} is the subject
of this paper. Here are the most important adjustments that one needs
to make in order to formulate this correspondence:
\begin{itemize}

\item The moduli stack $\Bun_G$ of $G$-bundles has to be replaced by
the moduli stack of $G$-bundles together with the level structures at
the ramification points. We call them the {\em enhanced} moduli
stacks. Recall that a level structure of order $N$ is a trivialization
of the bundle on the $N$th infinitesimal neighborhood of the
point. The order of the level structure should be at least the order
of the pole of the connection at this point.

\item At the points at which the connection has regular singularity
(pole of order $1$) one can take instead of the level structure, a
parabolic structure, i.e., a reduction of the fiber of the bundle to a
Borel subgroup of $G$.

\item The Langlands correspondence will assign to a flat $^L G$-bundle
  $E=(\F,\nabla)$ with ramification at the points $y_1,\ldots,y_n$ a
  {\bf category} ${\mc Aut}_E$ of Hecke eigensheaves on the
  corresponding enhanced moduli stack with eigenvalue $E|_{X \bs \{
  y_1,\ldots,y_n \}}$, which is a subcategory of the category of
  (twisted) ${\mc D}$-modules on this moduli stack.
\end{itemize}

If $E$ is unramified, then we may consider the category ${\mc Aut}_E$
on the moduli stack $\Bun_G$ itself. We then expect that for generic
$E$ this category is equivalent to the category of vector spaces: its
unique (up to isomorphism) irreducible object is $\Aut_E$ discussed
above, and all other objects are direct sums of copies of
$\Aut_E$. Because this category is expected to have such a simple
structure, it makes sense to say that the unramified geometric
Langlands correspondence assigns to an unramified $^L G$-local system
on $X$ a single Hecke eigensheaf, rather than a category. This is not
possible for general ramified local systems.

The questions that we are facing now are

\begin{itemize}
\item[(1)] How to construct the categories of Hecke eigensheaves for
ramified local systems?

\item[(2)] How to describe them in terms of the Langlands dual group
$^L G$?
\end{itemize}

In this article I will review an approach to these questions which has
been developed by D. Gaitsgory and myself in \cite{FG:local}.

The idea goes back to the construction of A. Beilinson and V. Drinfeld
\cite{BD} of the unramified geometric Langlands correspondence, which
may be interpreted in terms of a localization functor. Functors of
this type were introduced by A. Beilinson and J. Bernstein \cite{BB}
in representation theory of simple Lie algebras. In our situation this
functor sends representations of the affine Kac-Moody algebra $\ghat$
to twisted ${\mc D}$-modules on $\Bun_G$, or its enhanced versions. As
explained in \cite{F:rev} (see also \cite{vertex}), these ${\mc
D}$-modules may be viewed as sheaves of conformal blocks (or
coinvariants) naturally arising in the framework of Conformal Field
Theory.

The affine Kac-Moody algebra $\ghat$ is the universal one-dimensional
central extension of the loop algebra $\g\ppart$. The representation
categories of $\ghat$ have a parameter $\ka$, called the level, which
determines the scalar by which a generator of the one-dimensional
center of $\ghat$ acts on representations. We consider a particular
value $\ka_c$ of this parameter, called the {\bf critical level}. The
completed enveloping algebra of an affine Kac-Moody algebra acquires
an unusually large center at the critical level and this makes the
structure of the corresponding category $\gmod$ very rich and
interesting. B. Feigin and I have shown \cite{FF:gd,F:wak} that this
center is canonically isomorphic to the algebra of functions on the
space of $^L G$-{\bf opers} on $D^\times$. Opers are bundles on
$D^\times$ with flat connection and an additional datum (as defined by
Drinfeld-Sokolov \cite{DS} and Beilinson-Drinfeld \cite{BD}; we recall
the definition below). Remarkably, their structure group turns out to
be not $G$, but the Langlands dual group $^L G$, in agreement with the
general Langlands philosophy.

This result means that the category $\gmod$ of (smooth)
$\ghat$-modules of critical level ``lives'' over the space
$\on{Op}_{^L G}(D^\times)$ of $^L G$-opers on the punctured disc
$D^\times$. For each $\chi \in \on{Op}_{^L G}(D^\times)$ we have a
``fiber'' category $\gmod_\chi$ whose objects are $\ghat$-modules on
which the center acts via the central character corresponding to
$\chi$. Applying the localization functors to these categories, and
their $K$-equivariant subcategories $\gmod_\chi^K$ for various
subgroups $K \subset G[[t]]$, we obtain categories of Hecke
eigensheaves on the moduli spaces of $G$-bundles on $X$ with level (or
parabolic) structures.

Thus, the localization functor gives us a powerful tool for converting
{\bf local} categories of representations of $\ghat$ into {\bf global}
categories of Hecke eigensheaves. This is a new phenomenon which does
not have any obvious analogues in the classical Langlands
correspondence.

The simplest special case of this construction gives us the
Beilinson-Drinfeld Hecke eigensheaves $\Aut_E$ on $\Bun_G$
corresponding to unramified $^L G$-local systems admitting the oper
structure. Motivated by this, we wish to apply the localization
functors to more general categories $\gmod_\chi^K$ of $\ghat$-modules
of critical level, corresponding to opers on $X$ with singularities,
or ramifications.

These categories $\gmod_\chi$ are assigned to $^L G$-opers $\chi$ on
the punctured disc $D^\times$. It is important to realize that the
formal loop group $G\ppart$ naturally acts on each of these categories
via its adjoint action on $\ghat_{\ka_c}$ (because the center is
invariant under the adjoint action of $G\ppart$). Thus, we assign to
each oper $\chi$ a categorical representation of $G\ppart$ on
$\gmod_\chi$.

This is analogous to the classical {\bf local Langlands
correspondence}. Let $F$ be a local non-archimedian field, such as the
field $\Fq\ppart$ or the field of $p$-adic numbers. Let $W'_F$ be the
Weil-Deligne group of $F$, which is a version of the Galois group of
$F$ (we recall the definition in \secref{langlands param}). The local
Langlands correspondence relates the equivalence classes of
irreducible (smooth) representations of the group $G(F)$ (or
``$L$-packets'' of such representations) and the equivalence classes
of (admissible) homomorphisms $W'_F \to {}^L G$. In the geometric
setting we replace these homomorphisms by flat $^L G$-bundles on
$D^\times$ (or by $^L G$-opers), the group $G(F)$ by the loop group
$G\ppart$ and representations of $G(F)$ by categorical representations
of $G\ppart$.

This analogy is very suggestive, as it turns out that the structure of
the categories $\gmod_\chi$ (and their $K$-equivariant subcategories
$\gmod_\chi^K$) is similar to the structure of irreducible
representations of $G(F)$ (and their subspaces of $K$-invariants). We
will see examples of this parallelism in Sects. 7 and 8 below. This
means that what we are really doing is developing a {\bf local
Langlands correspondence for loop groups}.

To summarize, our strategy \cite{FG:local} for constructing the global
geometric Langlands correspondence has two parts:

\begin{itemize}
\item[(1)] the local part: describing the structure of the categories
of $\ghat$-modules of critical level, and

\item[(2)] the global part: applying the localization
functor to these categories to obtain the categories of Hecke
eigensheaves on enhanced moduli spaces of $G$-bundles.
\end{itemize}

We expect that these localization functors are equivalences of
categories (at least, in the generic situation), and therefore we can
infer a lot of information about the global categories by studying the
local categories $\gmod_\chi$ of $\ghat$-modules. Thus, the local
categories $\gmod_\chi$ take the center stage.

In this paper I review the results and conjectures of
\cite{FG:exact}--\cite{FG:weyl} with the emphasis on unramified and
tamely ramified local systems. (I also discuss the case of irregular
singularities at the end.) In particular, our study of the categories
of $\ghat_{\ka_c}$-modules leads us to the following conjecture. (For
related results, see \cite{AB,ABG,Bez1,Bez}.)

Suppose that $E = (\F,\nabla)$, where $\F$ is a $^L G$-bundle and
$\nabla$ is a connection on $\F$ with regular singularity at a single
point $y \in X$ and unipotent monodromy (this is easy to generalize to
multiple points). Let $M = \exp(2\pi i u)$, where $u \in {}^L G$ be a
representative of the conjugacy class of the monodromy of $\nabla$
around $y$. Denote by $\on{Sp}_u$ the {\bf Springer fiber} of $u$, the
variety of Borel subalgebras of $^L \g$ containing $u$.  The category
${\mc Aut}_E$ of Hecke eigensheaves with eigenvalue $E$ may then be
realized as a subcategory of the category of ${\mc D}$-modules on the
moduli stack of $G$-bundles on $X$ with parabolic structure at the
point $y$. We have the following conjectural description of the
derived category of ${\mc Aut}_E$:
$$
D^b({\mc Aut}_E) \simeq D^b(\QCoh(\on{Sp}_u^{\on{DG}})),
$$
where $\QCoh(\on{Sp}_u^{\on{DG}})$ is the category of quasicoherent
sheaves on a suitable ``DG enhancement'' of
$\on{Sp}_u^{\on{DG}}$. This is a category of differential graded (DG)
modules over a sheaf of DG algebras whose zeroth cohomology is the
structure sheaf of ${\mc B}^{\on{DG}}_M$ (we discuss this in detail in
\secref{global}).

Thus, we expect that the geometric Langlands correspondence attaches
to a $^L G$-local system on a Riemann surface with regular singularity
at a puncture, a category which is closely related to the variety of
Borel subgroups containing the monodromy around the puncture. We hope
that further study of the categories of $\ghat$-modules will help us
to find a similar description of the Langlands correspondence for
connections with irregular singularities.

\medskip

The paper is organized as follows. In Sect. 1 we review the
Beilinson-Drinfeld construction in the unramified case, in the
framework of localization functors from representation categories of
affine Kac-Moody algebras to ${\mc D}$-modules on on $\Bun_G$. This
will serve as a prototype for our construction of more general
categories of Hecke eigensheaves, and it motivates us to study
categories of $\ghat$-modules of critical level. We wish to interpret
these categories in the framework of the local geometric Langlands
correspondence for loop groups. In order to do that, we first recall
in Sect. 2 the setup of the classical Langlands correspondence. Then
in Sect. 3 we explain the passage to the geometric context. In Sect. 4
we describe the structure of the center at the critical level and the
isomorphism with functions on opers. In Sect. 5 we discuss the
connection between the local Langlands parameters ($^L G$-local
systems on the punctured disc) and opers. We introduce the categorical
representations of loop groups corresponding to opers and the
corresponding categories of Harish-Chandra modules in Sect. 6. We
discuss these categories in detail in the unramified case in Sect. 7,
paying particular attention to the analogies between the classical and
the geometric settings. In Sect. 8 we do the same in the tamely
ramified case. We then apply localization functor to these categories
in Sect. 9 to obtain various results and conjectures on the global
Langlands correspondence, both for regular and irregular singularities.

\medskip

Much of the material of this paper is borrowed from my new
book \cite{newbook}, where I refer the reader for more details, in
particular, for background on representation theory of affine
Kac-Moody algebras of critical level.

\medskip

Finally, I note that in a forthcoming paper \cite{GW} the geometric
Langlands correspondence with tame ramification is studied from the
point of view of dimensional reduction of four-dimensional
supersymmetric Yang-Mills theory.

\medskip

{\bf Acknowledgments.} I thank D. Gaitsgory for his collaboration on
our joint works which are reviewed in this article. I am also grateful
to R. Bezrukavnikov, V. Ginzburg, D. Kazhdan and E. Witten for useful
discussions.

I thank the organizers of the CIME Summer School in Venice,
especially, A. D'Agnolo, for the invitation to give lectures on this
subject at this enjoyable conference.

\section{The unramified global Langlands correspondence}
\label{unram first}

Our goal in this section is to construct Hecke eigensheaves $\Aut_E$
corresponding to unramified $^L G$-local systems $E = (\F,\nabla)$ on
$X$. By definition, $\Aut_E$ is a ${\mc D}$-module on $\Bun_G$. We
would like to construct $\Aut_E$ by applying a localization functor to
representations of affine Kac-Moody algebra $\ghat$.

Throughout this paper, unless specified otherwise, we let $\g$ be a
simple Lie algebra and $G$ the corresponding connected and
simply-connected algebraic group.

The key observation used in constructing the localization functor is
that for a simple Lie group $G$ the moduli stack $\Bun_G$ of
$G$-bundles on $X$ has a realization as a double quotient. Namely, let
$x$ be a point of $X$. Denote by $\K_x$ the completion of the field of
rational functions on $X$ at $x$, and by $\OO_x$ its ring of
integers. If we choose a coordinate $t$ at $x$, then we may identify
$\K_x \simeq \C\ppart, \OO_x \simeq \C[[t]]$. But in general there is
no preferred coordinate, and so it is better not to use these
identifications. Now let $G(\K_x) \simeq G\ppart$ be the formal loop
group corresponding to the punctured disc $D_x^\times$ around $x$. It
has two subgroups: one is $G(\OO_x) \simeq G[[t]]$ and the other is
$G_{\out}$, the group of algebraic maps $X \bs x \to G$. Then,
according to \cite{BL,DSimp}, the algebraic stack $\Bun_G$ is
isomorphic to the double quotient
\begin{equation}    \label{global uniform}
\Bun_G \simeq G_{\out}\bs G(\K_x)/G(\OO_x).
\end{equation}
Intuitively, any $G$-bundle may be trivialized on the formal disc
$D_x$ and on $X \bs x$. The transition function is then an element of
$G(\K_x)$, which characterizes the bundle uniquely up to the right
action of $G(\OO_x)$ and the left action of $G_{\out}$ corresponding
to changes of trivializations on $D_x$ and $X \bs x$, respectively.

The localization functor that we need is a special case of the
following general construction. Let $\g$ be a Lie algebra and $K$ a
Lie group $(\g,K)$ whose Lie algebra is contained in $\g$. The pair
$(\g,K)$ is called a Harish-Chandra pair. We will assume that $K$ is
connected. A $\g$-module $M$ is called $K$-equivariant if the action
of the Lie subalgebra $\on{Lie} K \subset \g$ on $M$ may be
exponentiated to an action of the Lie group $K$. Let $\g\mod^K$ be the
category of $K$-equivariant $\g$-modules.

Now suppose that $H$ is another subgroup of $G$. Let ${\mc D}_{H\bs
G/K}\mod$ be the category of ${\mc D}$-modules on $H\bs G/K$. Then
there is a localization functor \cite{BB,BD} (see also
\cite{F:rev,vertex})
$$
\Delta: \g\mod^K \to {\mc D}_{H\bs G/K}\mod.
$$
Now let $\ghat$ be a one-dimensional central extension of $\g$ which
becomes trivial when restricted to the Lie subalgebras $\on{Lie} K$
and $\on{Lie} H$. Suppose
that this central extension can be exponentiated to
a central extension $\wh{G}$ of the corresponding Lie group $G$. Then
we obtain a $\C^\times$-bundle $H \bs \wh{G}/K$ over $H \bs G/K$. Let
$\Ll$ be the corresponding line bundle and ${\mc D}_{\Ll}$ the
sheaf of differential operators acting on $\Ll$. Then we have a
functor
$$
\Delta_{\Ll}: \ghat\mod^K \to {\mc D}_{\Ll}\mod.
$$

In our case we take the formal loop group $G(\K_x)$, and the subgroups
$K = G(\OO_x)$ and $H = G_{\out}$ of $G(\K_x)$. We also consider the
so-called critical central extension of $G(\K_x)$. Let us first
discuss the corresponding central extension of the Lie algebra $\g
\otimes \K_x$. Choose a coordinate $t$ at $x$ and identify $\K_x
\simeq \C\ppart$. Then $\g \otimes \K_x$ is identified with
$\g\ppart$. Let $\ka$ be an invariant bilinear form on $\g$. The {\bf
affine Kac-Moody algebra} $\ghat_\ka$ is defined as the central
extension
$$
0 \to \C {\mb 1} \to \ghat_\ka \to \g\ppart \to 0.
$$
As a vector space, it is equal to the direct sum $\g\ppart \oplus \C {\mb
1}$, and the commutation relations read
\begin{equation}    \label{KM rel1}
[A \otimes f(t),B \otimes g(t)] = [A,B] \otimes f(t) g(t) -
(\kappa(A,B) \on{Res} f dg) {\mb 1},
\end{equation}
and ${\mb 1}$ is a central element, which commutes with everything
else. For a simple Lie algebra $\g$ all invariant inner products are
proportional to each other. Therefore the Lie algebras $\ghat_\ka$ are
isomorphic to each other for non-zero inner products $\ka$.

Note that the restriction of the second term in \eqref{KM rel1} to the
Lie subalgebra $\g \otimes t^N \C[[t]]$, where $N \in \Z_+$, is equal
to zero, and so it remains a Lie subalgebra of $\ghat_\ka$. A
$\ghat_\ka$-module is called {\bf smooth} if every vector in it is
annihilated by this Lie subalgebra for sufficiently large $N$. We
define the category $\ghat_\ka\mod$ whose objects are smooth
$\ghat_\ka$-modules on which the central element ${\mb 1}$ acts as the
identity. The morphisms are homomorphisms of representations of
$\ghat_\ka$. Throughout this paper, unless specified otherwise, by a
``$\ghat_\ka$-module'' will always mean a module on which the central
element ${\mb 1}$ acts as the identity.\footnote{Note that we could
have ${\mb 1}$ act instead as $\la$ times the identity for $\la \in
\C^\times$; but the corresponding category would just be equivalent to
the category $\ghat_{\la \ka}\mod$.} We will refer to $\ka$ as the
{\bf level}.

Now observe that formula \eqref{KM rel1} is independent of the choice
of coordinate $t$ at $x \in X$ and therefore defines a central
extension of $\g \otimes \K_x$, which we denote by $\ghat_{\ka,x}$. One
can show that this central extension may be exponentiated to a central
extension of the group $G(\K_x)$ if $\ka$ satisfies a certain
integrality condition, namely, $\ka = k \ka_0$, where $k \in \Z$ and
$\ka_0$ is the inner product normalized by the condition that the
square of the length of the maximal root is equal to $2$. A particular
example of the inner product which satisfies this condition is the
{\bf critical level} $\ka_c$ defined by the formula
\begin{equation}    \label{kacrit}
\ka_c(A,B) = - \frac{1}{2} \on{Tr}_{\g} \on{ad}A \on{ad}B.
\end{equation}
Thus, $\ka_c$ is equal to minus one half of the Killing form on
$\g$.\footnote{It is also equal to $-h^\vee \ka_0$, where $h^\vee$ is the
dual Coxeter number of $\g$.} When $\ka=\ka_c$ representation theory
of $\ghat_\ka$ changes dramatically because the completed enveloping
algebra of $\ghat_\ka$ acquires a large center (see below).

Let $\wh{G}_x$ be the corresponding critical central extension of
$G(\K_x)$. It is known (see \cite{BD}) that in this case the
corresponding line bundle $\Ll$ is the square root $K^{1/2}$ of the
canonical line bundle on $\Bun_G$.\footnote{Recall that by our
assumption $G$ is simply-connected. In this case there is a unique
square root.} Now we are ready to apply the localization functor in
the situation where our group is $G(\K_x)$, with the two subgroups
$K=G(\OO_x)$ and $H=G_{\out}$, so that the double quotient $H \bs G/K$
is $\Bun_G$.\footnote{Since $\Bun_G$ is an algebraic stack, one needs
to be careful in applying the localization functor. The appropriate
formalism has been developed in \cite{BD}.} We choose $\Ll =
K^{1/2}$. Then we have a localization functor
$$
\Delta_{\ka_c,x}: \gmodx^{G(\OO_x)} \to {\mc D}_{\ka_c}\mod.
$$
As explained in Part III of \cite{F:rev} (see also \cite{vertex}), the
${\mc D}$-modules obtained by applying this functor may be viewed as
sheaves of conformal blocks (or coinvariants) naturally arising in the
framework of Conformal Field Theory.

We will apply this functor to a particular $\ghat_{\ka_c,x}$-module.
To construct this module, let us first define the vacuum module over
$\ghat_{\ka_c,x}$ as the induced module
$$
\BV_{0,x} = {\rm Ind}_{\g \tensor\O_x \oplus \C {\mb
1}}^{\ghat_{\ka,x}} \C,
$$
where $\g \tensor\O_x$ acts by $0$ on $\C$ and ${\mb 1}$ acts as the
identity. According to the results of \cite{FF:gd,F:wak}, we have
$$
\on{End}_{\ghat_{\ka_c}} \BV_{0,x} \simeq \on{Fun} \on{Op}_{^L
  G}(D_x),
$$
where $\on{Op}_{^L G}(D_x)$ is the space of $^L G$-opers on the formal
disc $D_x = \on{Spec} \OO_x$ around $x$. We will give the definition of
this space and discuss this isomorphsm in the next section.

Now, given $\chi_x \in \on{Op}_{^L G}(D_x)$, we obtain a maximal ideal
$I(\chi_x)$ in $\on{End}_{\ghat_{\ka_c}} \BV_{0,x}$. Let
$\BV_0(\chi_x)$ be the $\ghat_{\ka_c,x}$-module which is the quotient
of $\BV_{0,x}$ by the image of $I(\chi_x)$ (it is non-zero, as
explained in \secref{cat unram}). The module $\BV_{0,x}$ is clearly
$G(\OO_x)$-equivariant, and hence so is $\BV_0(\chi_x)$. Therefore
$\BV_0(\chi_x)$ is an object of the category $\gmodx^{G(\OO_x)}$.

We now apply the localization functor $\Delta_{\ka_c,x}$ to
$\BV_0(\chi_x)$. The following theorem is due to Beilinson and
Drinfeld \cite{BD}.

\begin{thm}    \label{bd}
{\em (1)} The $\D_{\ka_c}$-module $\Delta_{\ka_c,x}(\BV_0(\chi_x))$ is
non-zero if and only if there exists a global $^L \g$-oper on $X$,
$\chi \in \on{Op}_{^L G}(X)$ such that $\chi_x \in \on{Op}_{^L
G}(D_x)$ is the restriction of $\chi$ to $D_x$.

{\em (2)} If this holds, then $\Delta_{\ka_c,x}(\BV_0(\chi_x))$ depends
only on $\chi$ and is independent of the choice of $x$ in the sense
that for any other point $y \in X$, if $\chi_y = \chi|_{D_y}$, then
$\Delta_{\ka_c,x}(\BV_0(\chi_x)) \simeq
\Delta_{\ka_c,y}(\BV_0(\chi_y))$.

{\em (3)} For any $\chi = (\F,\nabla,\F_{^L B}) \in \on{Op}_{^L
G}(X)$ the $\D_{\ka_c}$-module $\Delta_{\ka_c,x}(\BV_0(\chi_x))$ is
a non-zero Hecke eigensheaf with the eigenvalue $E_\chi = (F,\nabla)$.
\end{thm}

Thus, for any $\chi \in \on{Op}_{^L G}(X)$, the $\D_{\ka_c}$-module
$\Delta_{\ka_c,x}(\BV_0(\chi_x))$ is the sought-after Hecke eigensheaf
$\Aut_{E_\chi}$ corresponding to the $^L G$-local system under the
global geometric Langlands correspondence.\footnote{More
precisely, $\Aut_{E_\chi}$ is the ${\mc D}$-module
$\Delta_{\ka_c,x}(\BV_0(\chi_x)) \otimes K^{-1/2}$, but here and below
we will ignore the twist by $K^{1/2}$.} For an outline of the proof of
this theorem from \cite{BD}, see \cite{F:rev}, Sect. 9.4.

A drawback of this construction is that not all $^L G$-local systems
on $X$ admit the structure of an oper. In fact, under our assumption
that $G$ is simply-connected (and so $^L G$ is of adjoint type), the
local systems, or flat bundles $(\F,\nabla)$, on a smooth projective
curve $X$ that admit an oper structure correspond to a unique $^L
G$-bundle on $X$ described as follows (see \cite{BD}). Let
$\Omega^{1/2}_X$ be a square root of the canonical line bundle
$\Omega_X$. There is a unique (up to an isomorphism) non-trivial
extension
$$
0 \to \Omega_X^{1/2} \to \F_0 \to \Omega_X^{-1/2} \to 0.
$$
Let $\F_{PGL_2}$ be the $PGL_2$-bundle corresponding to the rank two
vector bundle $\F_0$. Note that it does not depend on the choice of
$\Omega_X^{1/2}$. This is the oper bundle for $PGL_2$. We define the
oper bundle $\F_{^L G}$ for a general simple Lie group $^L G$ of adjoint
type as the push-forward of $\F_{PGL_2}$ with respect to a principal
embedding $PGL_2 \hookrightarrow G$ (see \secref{canon represent}).

For each flat connection $\nabla$ on the oper bundle $\F_{^L G}$ there
exists a unique $^L B$-reduction $\F_{^L B}$ satisfying the oper
condition. Therefore $\Op$ is a subset of $\on{Loc}_{^L G}(X)$, which
is the fiber of the forgetful map $\on{Loc}_{^L G}(X) \to \Bun_{^L G}$
over $\F_{^L G}$.

\thmref{bd} gives us a construction of Hecke eigensheaves for $^L
G$-local system that belong to the locus of opers. For a general $^L
G$-local system outside this locus, the above construction may be
generalized as discussed at the end of \secref{outside the locus}
below.

Thus, \thmref{bd}, and its generalization to other unramified $^L
G$-local systems, give us an effective tool for constructing Hecke
eigensheaves on $\Bun_G$. It is natural to ask whether it can be
generalized to the ramified case if we consider more general
representations of $\ghat_{\ka_c,x}$. The goal of this paper is to
explain how to do that.

We will see below that the completed universal enveloping algebra of
$\ghat_{\ka_c,x}$ contains a large center. It is isomorphic to the
algebra $\on{Fun} \on{Op}_{^L G}(D_x^\times)$ of functions on the
space $\on{Op}_{^L G}(D_x^\times)$ of $^L G$-opers on the punctured
disc $D_x^\times$. For $\chi_x \in \on{Op}_{^L G}(D_x^\times)$, let
$\gmodx_{\chi_x}$ be the full subcategory of $\gmodx$ whose objects
are $\ghat_{\ka_c,x}$-modules on which the center acts according to
the character corresponding to $\chi_x$.

The construction of Hecke eigensheaves now breaks into two steps:

\begin{itemize}
\item[(1)] we study the Harish-Chandra categories $\gmodx_{\chi_x}^K$
for various subgroups $K \subset G(\OO_x)$;

\item[(2)] we apply the localization functors to these categories.
\end{itemize}

The simplest case of this construction is precisely the
Beilinson-Drinfeld construction explained above. In this case we take
$\chi_x$ to be a point in the subspace $\on{Op}_{^L G}(D_x) \subset
\on{Op}_{^L G}(D_x^\times)$. Then the category
$\gmodx_{\chi_x}^{G(\OO_x)}$ is equivalent to the category of vector
spaces: its unique up to an isomorphism irreducible object is the
above $\BV_0(\chi_x)$, and all other objects are direct sums of copies
of $\BV_0(\chi_x)$ (see \cite{FG:exact} and \thmref{single nabla}
below). Therefore the localization functor $\Delta_{\ka_c,x}$ is
determined by $\Delta_{\ka_c,x}(\BV_0(\chi_x))$, which is described in
\thmref{bd}. It turns out to be the desired Hecke eigensheaf
$\Aut_{E_{\chi}}$. Moreover, we expect that the functor
$\Delta_{\ka_c,x}$ sets up an equivalence between
$\gmodx_{\chi_x}^{G(\OO_x)}$ and the category of Hecke eigensheaves on
$\Bun_G$ with eigenvalue $E_\chi$.

For general opers $\chi_x$, with ramification, the (local) categories
$\gmodx_{\chi_x}^K$ are more complicated, as we will see below, and so
are the corresponding (global) categories of Hecke eigensheaves. In
order to understand the structure of the global categories, we need to
study first of local categories of $\ghat_{\ka_c,x}$-modules. Using
the localization functor, we can then understand the structure of the
global categories. We will consider examples of the local categories
in the following sections.

It is natural to view our study of the local categories
$\gmodx_{\chi_x}$ and $\gmodx_{\chi_x}^K$ as a geometric analogue of
the local Langlands correspondence. We will explain this point of view
in the next section.

\section{Classical local Langlands correspondence}

The local Langlands correspondence relates smooth representations of
reductive algebraic groups over local fields and representations of
the Galois group of this field. In this section we define these
objects and explain the main features of this correspondence. As the
material of this section serves motivational purposes, we will only
mention those aspects of this story that are most relevant for us. For
a more detailed treatment, we refer the reader to the informative
surveys \cite{Vogan,Kudla} and references therein.

The local Langlands correspondence may be formulated for any local
non-archimedian field. There are two possibilities: either $F$ is the
field $\Qp$ of $p$-adic numbers or a finite extension of $\Qp$, or $F$
is the field $\Fq\ppart$ of formal Laurent power series with
coefficients in $\Fq$, the finite field with $q$ elements (where $q$
is a power of a prime number). For the sake of definiteness, in what
follows we will restrict ourselves to the second case.

\subsection{Langlands parameters}    \label{langlands param}

Consider the group $GL_n(F)$, where $F=\Fq\ppart$. A representation of
$GL_n(F)$ on a complex vector space $V$ is a homomorphism $\pi:
GL_n(F) \to \on{End} V$ such that $\pi(gh) = \pi(g) \pi(h)$ and
$\pi(1) = \on{Id}$. Define a topology on $GL_n(F)$ by stipulating that
the base of open neighborhoods of $1 \in GL_n(F)$ is formed by the
congruence subgroups
$$
K_N = \{ g \in GL_n(\Fq[[t]])\, | \, g \equiv 1 \; \on{mod} \; t^N \},
\qquad N \in \Z_+.
$$
For each $v \in V$ we obtain a map $\pi(\cdot) v: GL_n(F) \to V, g
\mapsto \pi(g)v$. A representation $(V,\pi)$ is called {\bf smooth}
\index{smooth representation} if the map $\pi(\cdot) v$ is continuous
for each $v$, where we give $V$ the discrete topology. In other words,
$V$ is smooth if for any vector $v \in V$ there exists $N \in \Z_+$
such that
$$
\pi(g) v = v, \qquad \forall g \in K_N.
$$

We are interested in describing the equivalence classes of irreducible
smooth representations of $GL_n(F)$. Surprisingly, those turn out to
be related to objects of a different kind: $n$-dimensional
representations of the Galois group of $F$.

Recall that the algebraic closure of $F$ is a field obtained by
adjoining to $F$ the roots of all polynomials with coefficients in
$F$. However, in the case when $F = \Fq\ppart$ some of the extensions
of $F$ may be non-separable. We wish to avoid the non-separable
extensions, because they do not contribute to the Galois group. Let
$\ol{F}$ be the maximal separable extension inside a given algebraic
closure of $F$. It is uniquely defined up to an isomorphism.

Let $\Gal(\ol{F}/F)$ be the {\bf absolute Galois group} \index{Galois
group} of $F$. Its elements are the automorphisms $\sigma$ of the
field $\ol{F}$ such that $\sigma(y) = y$ for all $y \in F$.

Now set $F=\Fq\ppart$. Observe that we have a natural map
$\Gal(\ol{F}/F) \to \Gal(\ol{\mathbb F}_q/\Fq)$ obtained by applying
an automorphism of $F$ to $\ol{\mathbb F}_q \subset F$. The group
$\Gal(\ol{\mathbb F}_q/\Fq)$ is isomorphic to the profinite completion
$\wh\Z$ of $\Z$ (see, e.g., \cite{F:rev}, Sect. 1.3). Its subgroup $\Z
\subset \wh\Z$ is generated by the {\bf geometric Frobenius element}
which is inverse to the automorphism $x \mapsto x^q$ of $\ol{\mathbb
F}_q$. Let $W_F$ be the preimage of the subgroup $\Z \subset
\Gal(\ol{\mathbb F}_q/\Fq)$. This is the {\bf Weil group} of
$F$. \index{Weil group} Denote by $\nu$ be the corresponding
homomorphism $W_F \to \Z$.

Now let $W'_F = W_F \ltimes \C$ be the semi-direct product of $W_F$
and the one-dimensional complex additive group $\C$, where $W_F$ acts
on $\C$ by the formula
\begin{equation}    \label{rel in weil}
\sigma x \sigma^{-1} = q^{\nu(\sigma)} x, \qquad \sigma \in W_F, x
\in \C.
\end{equation}
This is the {\bf Weil-Deligne group} of $F$. \index{Weil-Deligne
  group}

An $n$-dimensional complex representation of $W'_F$ is by definition a
homomorphism $\rho': W'_F \to GL_n(\C)$ which may be described as a
pair $(\rho,N)$, where $\rho$ is an $n$-dimensional representation of
$W_F$, $N \in GL_n(\C)$, and we have $\rho(\sigma) N \rho(\sigma)^{-1}
= q^{\nu(\sigma)} \rho(N)$ for all $\sigma \in W_F$. The group $W_F$
is topological, with respect to the Krull topology (in which the open
neighborhoods of the identity are the normal subgroups of finite
index). The representation $(\rho,N)$ is called {\bf admissible} if
$\rho$ is continuous (equivalently, factors through a finite quotient
of $W_F$) and semisimple, and $N$ is a unipotent element of
$GL_n(\C)$.

The group $W'_F$ was introduced by P. Deligne \cite{Del}. The idea is
that by adjoining the unipotent element $N$ to $W_F$ one obtains a
group whose complex admissible representations are the same as
continuous $\ell$-adic representations of $W_F$ (where $\ell \neq p$
is a prime).

\subsection{The local Langlands correspondence for $GL_n$}
\label{llc for gln}

Now we are ready to state the local Langlands correspondence for the
group $GL_n$ over a local non-archimedian field $F$. It is a
bijection between two different sorts of data. One is the set of the
equivalence classes of irreducible smooth representations of
$GL_n(F)$. The other is the set of equivalence classes of
$n$-dimensional admissible representations of $W'_F$. We represent it
schematically as follows:

$$
\boxed{\begin{matrix} \text{$n$-dimensional admissible} \\
\text{representations of } W'_F \end{matrix}} \quad
\Longleftrightarrow  \quad \boxed{\begin{matrix} \text{irreducible
      smooth} \\ \text{representations of } GL_n(F)  \end{matrix}}
$$

This correspondence is supposed to satisfy an overdetermined system of
constraints which we will not recall here (see, e.g., \cite{Kudla}).

The local Langlands correspondence for $GL_n$ is a theorem. In the
case when $F=\Fq\ppart$ it has been proved in \cite{LRS}, and when
$F=\Qp$ or its finite extension in \cite{HT} and also in
\cite{Henniart}.

\subsection{Generalization to other reductive groups}    \label{dual
  grp}

Let us replace the group $GL_n$ by an arbitrary connected reductive
group $G$ over a local non-archimedian field $F$. The group $G(F)$ is
also a topological group, and there is a notion of smooth
representation of $G(F)$ on a complex vector space. It is natural to
ask whether one can relate irreducible smooth representations of
$G(F)$ to representations of the Weil-Deligne group $W'_F$. This
question is addressed in the general local Langlands conjectures. It
would take us too far afield to try to give here a precise formulation
of these conjectures. So we will only indicate some of the objects
involved referring the reader to the articles \cite{Vogan,Kudla} where
these conjectures are described in great detail.

Recall that in the case when $G=GL_n$ the irreducible smooth
representations are parametrized by admissible homomorphisms $W'_F \to
GL_n(\C)$. In the case of a general reductive group $G$, the
representations are conjecturally parametrized by admissible
homomorphisms from $W'_F$ to the so-called {\bf Langlands dual group}
\index{Langlands dual group} $^L G$, which is defined over $\C$.

In order to explain the notion of the Langlands dual group, consider
first the group $G$ over the closure $\ol{F}$ of the field $F$. All
maximal tori $T$ of this group are conjugate to each other and are
necessarily split, i.e., we have an isomorphism $T(\ol{F}) \simeq
(\ol{F}^\times)$. For example, in the case of $GL_n$, all maximal tori
are conjugate to the subgroup of diagonal matrices. We associate to
$T(\ol{F})$ two lattices: the weight lattice $X^*(T)$ of homomorphisms
$T(\ol{F}) \to \ol{F}^{\times}$ and the coweight lattice $X_*(T)$ of
homomorphisms $\ol{F}^{\times} \to T(\ol{F})$. They contain the sets
of roots $\Delta \subset X^*(T)$ and coroots $\Delta^\vee \subset
X_*(T)$, respectively. The quadruple
$(X^*(T),X_*(T),\Delta,\Delta^\vee)$ is called the root datum for $G$
over $\ol{F}$. The root datum determines $G$ up to an isomorphism
defined over $\ol{F}$. The choice of a Borel subgroup $B(\ol{F})$
containing $T(\ol{F})$ is equivalent to a choice of a basis in
$\Delta$, namely, the set of simple roots $\Delta_s$, and the
corresponding basis $\Delta^\vee_s$ in $\Delta^\vee$.

Now, given $\gamma \in \on{Gal}(\ol{F}/F)$, there is $g \in G(\ol{F})$
such that $g(\gamma(T(\ol{F}))g^{-1} = T(\ol{F})$ and
$g(\gamma(B(\ol{F}))g^{-1} = B(\ol{F})$. Then $g$ gives rise to an
automorphism of the based root data
$(X^*(T),X_*(T),\Delta_s,\Delta^\vee_s)$. Thus, we obtain an action of
$\on{Gal}(\ol{F}/F)$ on the based root data.

Let us now exchange the lattices of weights and coweights and the sets
of simple roots and coroots. Then we obtain the based root data
$$
(X_*(T),X^*(T),\Delta^\vee_s,\Delta_s)
$$
of a reductive algebraic group over $\C$ which is denoted by $^L
G^\circ$. For instance, the group $GL_n$ is self-dual, the dual of
$SO_{2n+1}$ is $Sp_{2n}$, the dual of $Sp_{2n}$ is $SO_{2n+1}$, and
$SO_{2n}$ is self-dual.

The action of $\on{Gal}(\ol{F}/F)$ on the based root data gives rise
to its action on $^L G^\circ$. The semi-direct product $^L G =
\on{Gal}(\ol{F}/F) \ltimes {}^L G^\circ$ is called the {\bf Langlands
dual group} of $G$.

According to the local Langlands conjecture, the equivalence classes
of irreducible smooth representations of $G(F)$ are, roughly speaking,
parameterized by the equivalence classes of admissible homomorphisms
$W'_F \to {}^L G$. In fact, the conjecture is more subtle: one needs
to consider simultaneously representations of all inner forms of $G$,
and a homomorphism $W'_F \to {}^L G$ corresponds in general not to a
single irreducible representation of $G(F)$, but to a finite set of
representations called an $L$-{\bf packet}. To distinguish between
them, one needs additional data (see \cite{Vogan} and \secref{tame
reps} below for more details). But in the first approximation one can
say that the essence of the local Langlands correspondence is that

\medskip

\noindent {\em Irreducible smooth representations of $G(F)$ are
parameterized in terms of admissible homomorphisms $W'_F \to {}^L G$.}

\section{Geometric local Langlands correspondence over $\C$}
\label{llc}

We now wish to find a generalization of the local Langlands
conjectures in which we replace the field $F=\Fq\ppart$ by the field
$\C\ppart$. We would like to see how the ideas and patterns of the
Langlands correspondence play out in this new context, with the hope
of better understanding the deep underlying structures behind this
correspondence.

So let $G$ be a connected simply-connected algebraic group over $\C$,
and $G(F)$ the loop group $G\ppart = G(\C\ppart)$. Thus, we wish to
study smooth representations of the loop group $G\ppart$ and try to
relate them to some ``Langlands parameters'', which we expect, by
analogy with the case of local non-archimedian fields described above,
to be related to the Galois group of $\C\ppart$ and the Langlands dual
group $^L G$.

\subsection{Geometric Langlands parameters}    \label{fund grp}

Unfortunately, the Galois group of $\C\ppart$ is too small: it is
isomorphic to the pro-finite completion $\wh\Z$ of $\Z$. This is not
surprising from the point of view of the analogy between the Galois
groups and the fundamental groups (see, e.g., \cite{F:rev},
Sect. 3.1). The topological fundamental group of the punctured disc is
$\Z$, and the algebraic fundamental group is its pro-finite
completion.

However, we may introduce additional Langlands parameters by using a
more geometric perspective on homomorphisms from the fundamental group
to $^L G$. Those may be viewed as $^L G$-local systems. In general,
$^L G$-local systems on a compact variety $Z$ are the same as flat $^L
G$-bundles $(\F,\nabla)$ on $Z$. If the variety is not compact (as in
the case of $D^\times$), then we should impose the additional
condition that the connection has {\bf regular singularities} (pole of
order at most $1$) at infinity. In our case we obtain $^L G$-bundles
on $D^\times$ with a connection that has regular singularity at the
origin. Then the monodromy of the connection gives rise to a
homomorphism from $\pi_1(D^\times)$ to $^L G$. Now we generalize this
by allowing connections with {\bf arbitrary}, that is regular and {\bf
irregular}, singularities at the origin. Thus, we want to use as the
general Langlands parameters, the equivalence classes of pairs
$(\F,\nabla)$, where $\F$ is a $^L G$-bundle on $D^\times$ and
$\nabla$ is an arbitrary connection on $\F$.

Any bundle $\F$ on $D^\times$ may be trivialized. Then $\nabla$ may be
represented by the first-order differential operator
\begin{equation}    \label{gauge eq classes}
\nabla = \pa_t + A(t), \qquad A(t) \in {}^L \g\ppart.
\end{equation}
where $^L \g$ is the Lie algebra of the Langlands dual group $^L
G$. Changing the trivialization of $\F$ amounts to a gauge
transformation
$$
\nabla \mapsto \nabla' = \pa_t + g A g^{-1} - (\pa_t g)
g^{-1}
$$
with $g \in {}^L G\ppart$. Therefore the set of equivalence
classes of $^L G$-bundles with a connection on $D^\times$ is in
bijection with the set of gauge equivalence classes of operators
\eqref{gauge eq classes}. We denote this set by $\Loc_{^L
G}(D^\times)$. Thus, we have
\begin{equation}    \label{Loc}
\Loc_{^L G}(D^\times) = \{ \pa_t + A(t), \; A(t) \in {}^L \g\ppart
\}/{}^L G\ppart.
\end{equation}

We declare that the local Langlands parameters in the complex
setting should be the points of $\Loc_{^L G}(D^\times)$: the
equivalence classes of flat $^L G$-bundles on $D^\times$ or, more
concretely, the gauge equivalence classes \eqref{Loc} of first-order
differential operators.

Having settled the issue of the Langlands parameters, we have to
decide what it is that we will be parameterizing. Recall that in the
classical setting the homomorphism $W'_F \to {}^L G$ parameterized
irreducible smooth representations of the group $G(F)$, $F =
\Fq\ppart$. We start by translating this notion to the representation
theory of loop groups.

\subsection{Representations of the loop group}

The loop group $G\ppart$ contains the congruence subgroups
\begin{equation}    \label{congr sub}
K_N = \{ g \in G[[t]] \, | \, g \equiv 1 \; \on{mod} \; t^N \}, \qquad
N \in \Z_+.
\end{equation}
It is natural to call a representation of $G\ppart$ on a complex
vector space $V$ {\bf smooth} \index{smooth representation} if for any
vector $v \in V$ there exists $N \in \Z_+$ such that $K_N \cdot v =
v$. This condition may be interpreted as the continuity condition, if
we define a topology on $G\ppart$ by taking as the base of open
neighborhoods of the identity the subgroups $K_N, N \in \Z_+$, as
before.

But our group $G$ is now a complex Lie group (not a finite group), and
so $G\ppart$ is an infinite-dimensional Lie group. More precisely, we
view $G\ppart$ as an ind-group, i.e., as a group object in the category
of ind-schemes. At first glance, it is natural to consider the
algebraic representations of $G\ppart$. We observe that $G\ppart$ is
generated by the ``parahoric'' algebraic groups $P_i$ corresponding to
the affine simple roots. For these subgroups the notion of algebraic
representation makes perfect sense. A representation of $G\ppart$ is
then said to be algebraic if its restriction to each of the $P_i$'s is
algebraic.

However, this naive approach leads us to the following discouraging
fact: an irreducible smooth representation of $G\ppart$, which is
algebraic, is necessarily trivial (see \cite{BD}, 3.7.11(ii)).
Thus, we find that the class of algebraic representations of loop
groups turns out to be too restrictive. We could relax this condition
and consider differentiable representations, i.e., the representations
of $G\ppart$ considered as a Lie group. But it is easy to see that the
result would be the same. Replacing $G\ppart$ by its central extension
$\wh{G}$ would not help us much either: irreducible integrable
representations of $\wh{G}$ are parameterized by dominant integral
weights, and there are no extensions between them \cite{Kac}. These
representations are again too sparse to be parameterized by the
geometric data considered above. Therefore we should look for other
types of representations.

Going back to the original setup of the local Langlands
correspondence, we recall that there we considered representations of
$G(\Fq\ppart)$ on $\C$-vector spaces, so we could not possibly use the
algebraic structure of $G(\Fq\ppart)$ as an ind-group over
$\Fq$. Therefore we cannot expect the class of algebraic (or
differentiable) representations of the complex loop group $G\ppart$ to
be meaningful from the point of view of the Langlands
correspondence. We should view the loop group $G\ppart$ as an abstract
topological group, with the topology defined by means of the congruence
subgroups, in other words, consider its smooth representations as an
{\em abstract} group.

So we need to search for some geometric objects that encapsulate
representations of our groups and make sense both over a finite field
and over the complex field.

\subsection{From functions to sheaves}    \label{from functions to
  sheaves}

We start by revisiting smooth representations of the group $G(F)$,
where $F = \Fq\ppart$. We realize such representations more concretely
by considering their matrix coefficients. Let $(V,\pi)$ be an
irreducible smooth representation of $G(F)$. We define the {\bf
contragredient} representation $V^\vee$ as the linear span of all
smooth vectors in the dual representation $V^*$. This span is stable
under the action of $G(F)$ and so it admits a smooth representation
$(V^\vee,\pi^\vee)$ of $G(F)$. Now let $\phi$ be a $K_N$-invariant
vector in $V^\vee$. Then we define a linear map
$$
V \to C(G(F)/K_N), \qquad v \mapsto f_v,
$$
where
$$
f_v(g) = \langle \pi^\vee(g) \phi,v \rangle.
$$
Here $C(G(F)/K_N)$ denotes the vector space of $\C$-valued locally
constant functions on $G(F)/K_N$. The group $G(F)$ naturally acts on
this space by the formula $(g \cdot f)(h) = f(g^{-1} h)$, and the
above map is a morphism of representations, which is non-zero, and
hence injective, if $(V,\pi)$ is irreducible.

Thus, we realize our representation in the space of functions on the
quotient $G(F)/K_N$. More generally, we may realize representations in
spaces of functions on the quotient $G\ppart/K$ with values in a
finite-dimensional vector space, by considering a
finite-dimen\-sio\-nal subrepresentation of $K$ inside $V$ rather than
the trivial one.

An important observation here is that $G(F)/K$, where $F = \Fq\ppart$)
and $K$ is a compact subgroup of $G(F)$, is not only a set, but it is
a set of points of an algebraic variety (more precisely, an
ind-scheme) defined over the field $\Fq$. For example, for $K_0 =
G(\Fq[[t]])$, which is the maximal compact subgroup, the quotient
$G(F)/K_0$ is the set of $\Fq$-points of the ind-scheme called
the {\bf affine Grassmannian}. \index{affine Grassmannian}

Next, we recall an important idea going back to Grothendieck that
functions on the set of $\Fq$-points on an algebraic variety $X$
defined over $\Fq$ can often be viewed as the ``shadows'' of the
so-called $\ell$-adic sheaves on $X$. We will not give the definition
of these sheaves, referring the reader to \cite{Milne,Weil}.  The
Grothendieck {\bf fonctions-faisceaux} dictionary (see, e.g.,
\cite{Laumon}) is formulated as follows. Let ${\mc F}$ be an
$\ell$-adic sheaf and $x$ be an ${\mathbb F}_{q_1}$-point of $X$,
where $q_1=q^m$. Then one has the Frobenius conjugacy class
$\on{Fr}_x$ acting on the stalk ${\mc F}_x$ of ${\mc F}$ at $x$. Hence
we can define a function $\text{\tt f}_{q_1}({\mc F})$ on the set of
${\mathbb F}_{q_1}$-points of $V$, whose value at $x$ is
$\on{Tr}(\on{Fr}_x,{\mc F}_x)$. This function takes values in the
algebraic closure $\ol{\mathbb Q}_\ell$ of ${\mathbb Q}_\ell$. But
there is not much of a difference between $\ol{\mathbb Q}_\ell$-valued
functions and $\C$-valued functions: since they have the same
cardinality, $\ol{\mathbb Q}_\ell$ and $\C$ may be identified as
abstract fields. Besides, in most interesting cases, the values
actually belong to $\ol\Q$, which is inside both $\ol{\mathbb Q}_\ell$
and $\C$.

More generally, if $\K$ is a complex of $\ell$-adic sheaves, one
defines a function $\text{\tt f}_{q_1}({\mc K})$ on $V({\mathbb
F}_{q_1})$ by taking the alternating sums of the traces of $\on{Fr}_x$
on the stalk cohomologies of $\K$ at $x$. The map $\K \to \text{\tt
f}_{q_1}({\mc K})$ intertwines the natural operations on sheaves with
natural operations on functions (see \cite{Laumon},
Sect. 1.2).

Let $K_0({\mc S}h_X)$ be the complexified Grothendieck group of the
category of $\ell$-adic sheaves on $X$. Then the above construction
gives us a map
$$
K_0({\mc S}h_X) \to \prod_{m\geq 1} X({\mathbb F}_{q^m}),
$$
and it is known that this map is injective (see
\cite{Laumon}).

Therefore we may hope that the functions on the quotients $G(F)/K_N$
which realize our representations come by this constructions from
$\ell$-adic sheaves, or more generally, from complexes of $\ell$-adic
sheaves, on $X$.

Now, the notion of constructible sheaf (unlike the notion of a
function) has a transparent and meaningful analogue for a complex
algebraic variety $X$, namely, those sheaves of $\C$-vector spaces
whose restrictions to the strata of a stratification of the variety
$X$ are locally constant. The affine Grassmannian and more general
ind-schemes underlying the quotients $G(F)/K_N$ may be defined both
over $\Fq$ and $\C$. Thus, it is natural to consider the categories of
such sheaves (or, more precisely, their derived categories) on these
ind-schemes over $\C$ as the replacements for the vector spaces of
functions on their points realizing smooth representations of the
group $G(F)$.

We therefore naturally come to the idea, advanced in \cite{FG:local},
that the representations of the loop group $G\ppart$ that we need to
consider are not realized on vector spaces, but on {\bf categories},
such as the derived category of coherent sheaves on the affine
Grassmannian. Of course, such a category has a Grothendieck group, and
the group $G\ppart$ will act on the Grothendieck group as well, giving
us a representation of $G\ppart$ on a vector space. But we obtain much
more structure by looking at the categorical representation. The
objects of the category, as well as the action, will have a geometric
meaning, and thus we will be using the geometry as much as possible.

Let us summarize: to each local Langlands parameter $\chi \in \Loc_{^L
G}(D^\times)$ we wish to attach a category ${\mc C}_\chi$ equipped
with an action of the loop group $G\ppart$. But what kind of categories
should these ${\mc C}_\chi$ be and what properties do we expect them
to satisfy?

To get closer to answering these questions, we wish to discuss two more
steps that we can make in the above discussion to get to the types of
categories with an action of the loop group that we will consider in
this paper.

\subsection{A toy model}    \label{finite groups}

At this point it is instructive to detour slightly and consider a toy
model of our construction.  Let $G$ be a split reductive group over
$\Z$, and $B$ its Borel subgroup. A natural representation of $G(\Fq)$
is realized in the space of complex (or $\ol{\mathbb Q}_\ell$-) valued
functions on the quotient $G(\Fq)/B(\Fq)$. It is natural to ask what
is the ``correct'' analogue of this representation if we replace the
field $\Fq$ by the complex field and the group $G(\Fq)$ by
$G(\C)$. This may be viewed as a simplified version of our quandary,
since instead of considering $G(\Fq\ppart)$ we now look at $G(\Fq)$.

The quotient $G(\Fq)/B(\Fq)$ is the set of $\Fq$-points of the
algebraic variety defined over $\Z$ called the flag variety of $G$ and
defined by $\on{Fl}$. Our discussion in the previous section suggests
that we first need to replace the notion of a function on
$\on{Fl}(\Fq)$ by the notion of an $\ell$-adic sheaf on the variety
$\on{Fl}_{\Fq} = \on{Fl} \underset{\Z}\otimes \Fq$.

Next, we replace the notion of an $\ell$-adic sheaf on $\on{Fl}$
considered as an algebraic variety over $\Fq$, by the notion of a
constructible sheaf on $\on{Fl}_{\C} = \on{Fl} \underset{\Z}\otimes
\C$ which is an algebraic variety over $\C$. The complex algebraic
group $G_{\C}$ naturally acts on $\on{Fl}_{\C}$ and hence on this
category. Now we make two more reformulations of this category.

First of all, for a smooth complex algebraic variety $X$ we have a
{\bf Riemann-Hilbert correspondence} \index{Riemann-Hilbert
correspondence} which is an equivalence between the derived category
of constructible sheaves on $X$ and the derived category of ${\mc
D}$-modules on $X$ that are holonomic and have regular singularities.

Here we consider the sheaf of algebraic differential operators on $X$
and sheaves of modules over it, which we simply refer to as
$\D$-modules. The simplest example of a $\D$-module is the sheaf of
sections of a vector bundle on $V$ equipped with a flat
connection. The flat connection enables us to multiply any section by
a function and we can use the flat connection to act on sections by
vector fields. The two actions generate an action of the sheaf of
differential operators on the sections of our bundle. The sheaf of
horizontal sections of this bundle is then a locally constant sheaf of
X. We have seen above that there is a bijection between the set of
isomorphism classes of rank $n$ bundles on $X$ with connection having
regular singularities and the set of isomorphism classes of locally
constant sheaves on $X$ of rank $n$, or equivalently,
$n$-dimensional representations of $\pi_1(X)$. This bijection may be
elevated to an equivalence of the corresponding categories, and the
general Riemann-Hilbert correspondence is a generalization of this
equivalence of categories that encompasses more general $\D$-modules.

The Riemann-Hilbert correspondence allows us to associate to any
holonomic $\D$-module on $X$ a complex of constructible sheaves on
$X$, and this gives us a functor between the corresponding derived
categories which turns out to be an equivalence if we restrict
ourselves to the holonomic $\D$-modules with regular singularities
(see \cite{Dmodules,GM} for more details).

Thus, over $\C$ we may pass from constructible sheaves to
$\D$-modules. In our case, we consider the category of (regular
holonomic) $\D$-modules on the flag variety $\on{Fl}_\C$. This
category carries a natural action of $G_{\C}$.

Finally, let us observe that the Lie algebra $\g$ of $G_{\C}$ acts on
the flag variety infinitesimally by vector fields. Therefore, given a
$\D$-module $\F$ on $\on{Fl}_\C$, the space of its global sections
$\Gamma(\on{Fl}_\C,\F)$ has the structure of $\g$-module. We obtain a
functor $\Gamma$ from the category of $\D$-modules on $\on{Fl}_\C$ to
the category of $\g$-modules. A. Beilinson and J. Bernstein have
proved that this functor is an equivalence between the category of all
$\D$-modules on $\on{Fl}_\C$ (not necessarily regular holonomic) and
the category ${\mc C}_0$ of $\g$-modules on which the center of the
universal enveloping algebra $U(\g)$ acts through the augmentation
character.

Thus, we can now answer our question as to what is a meaningful
geometric analogue of the representation of the finite group $G(\Fq)$
on the space of functions on the quotient $G(\Fq)/B(\Fq)$. The answer
is the following: it is an {\bf abelian category} equipped with an
action of the algebraic group $G_{\C}$. This category has two
incarnations: one is the category of $\D$-modules on the flag variety
$\on{Fl}_\C$, and the other is the category ${\mc C}_0$ of modules
over the Lie algebra $\g$ with the trivial central character. Both
categories are equipped with natural actions of the group $G_{\C}$.

Let us pause for a moment and spell out what exactly we mean when we
say that the group $G_{\C}$ acts on the category ${\mc C}_0$. For
simplicity, we will describe the action of the corresponding group
$G(\C)$ of $\C$-points of $G_{\C}$.\footnote{More generally, for any
$\C$-algebra $R$, we have an action of $G(R)$ on the corresponding
base-changed category over $R$. Thus, we are naturally led to the
notion of an algebraic group (or, more generally, a group scheme)
acting on an abelian category, which is spelled out in
\cite{FG:local}, Sect. 20.} This means the following: each element $g
\in G$ gives rise to a functor $F_g$ on ${\mc C}_0$ such that $F_1$ is
the identity functor, and the functor ${\mc F}_{g^{-1}}$ is
quasi-inverse to $F_g$. Moreover, for any pair $g,h \in G$ we have a
fixed isomorphism of functors $i_{g,h}: F_{gh} \to F_g \circ F_h$ so
that for any triple $g,h,k \in G$ we have the equality $i_{h,k}
i_{g,hk} = i_{g,h} i_{gh,k}$ of isomorphisms $F_{ghk} \to F_g \circ
F_h \circ F_k$.

The functors $F_g$ are defined as follows. Given a representation
$(V,\pi)$ of $\g$ and an element $g \in G(\C)$, we define a new
representation $F_g((V,\pi)) = (V,\pi_g)$, where by definition
$\pi_g(x) = \pi(\on{Ad}_g(x))$. Suppose that $(V,\pi)$ is
irreducible. Then it is easy to see that $(V,\pi_g) \simeq (V,\pi)$ if
and only if $(V,\pi)$ is integrable, i.e., is obtained from an
algebraic representation of $G$.\footnote{In general, we could obtain
a representation of a central extension of $G$, but if $G$ is
reductive, it does not have non-trivial central extensions.} This is
equivalent to this representation being finite-dimensional. But a
general representation $(V,\pi)$ is infinite-dimensional, and so it
will not be isomorphic to $(V,\pi_g)$, at least for some $g \in G$.

Now we consider morphisms in ${\mc C}_0$, which are just
$\g$-homomorphisms. Given a $\g$-homomorphism between representations
$(V,\pi)$ and $(V',\pi')$, i.e., a linear map $T: V \to V'$ such that
$T \pi(x) = \pi'(x) T$ for all $x \in \g$, we set $F_g(T) = T$. The
isomorphisms $i_{g,h}$ are all identical in this case.

\subsection{Back to loop groups}    \label{back to
  groups}

In our quest for a complex analogue of the local Langlands
correspondence we need to decide what will replace the notion of a
smooth representation of the group $G(F)$, where $F = \Fq\ppart$. As
the previous discussion demonstrates, we should consider
representations of the complex loop group $G\ppart$ on various
categories of $\D$-modules on the ind-schemes $G\ppart/K$, where $K$
is a ``compact'' subgroup of $G\ppart$, such as $G[[t]]$ or the
Iwahori subgroup (the preimage of a Borel subgroup $B \subset G$ under
the homomorphism $G[[t]] \to G$), or the categories of representations
of the Lie algebra $\g\ppart$. Both scenarios are viable, and they
lead to interesting results and conjectures which we will discuss in
detail in \secref{global}, following \cite{FG:local}. In this
paper we will concentrate on the second scenario and consider
categories of modules over the loop algebra $\g\ppart$.

The group $G\ppart$ acts on the category of representations of
$\g\ppart$ in the way that we described in the previous section. An
analogue of a smooth representation of $G(F)$ is a category of smooth
representations of $\g\ppart$.  Let us observe however that we could
choose instead the category of smooth representations of the central
extension of $\g\ppart$, namely, $\ghat_\ka$.

The group $G\ppart$ acts on the Lie algebra $\ghat_\ka$ for any $\ka$,
because the adjoint action of the central extension of $G\ppart$
factors through the action of $G\ppart$.  We use the action of
$G\ppart$ on $\ghat_\ka$ to construct an action of $G\ppart$ on the
category $\ghat_\ka\mod$, in the same way as in \secref{finite
groups}.

Now recall the space $\on{Loc}_{^L G}(D^\times)$ of the Langlands
parameters that we defined in \secref{fund grp}. Elements of
$\on{Loc}_{^L G}(D^\times)$ have a concrete description as gauge
equivalence classes of first order operators $\pa_t + A(t), A(t) \in
{}^L \g\ppart$, modulo the action of $^L G\ppart$ (see formula
\eqref{Loc}).

We can now formulate the local Langlands correspondence over $\C$ as
the following problem:

\medskip

\noindent{\em To each local Langlands parameter $\chi
\in \on{Loc}_{^L G}(D^\times)$ associate a subcategory
$\ghat_\ka\mod_\chi$ of $\ghat_\ka\mod$ which is stable under the
action of the loop group $G\ppart$.}

\medskip

We wish to think of the category $\ghat_\ka\mod$ as ``fibering'' over
the space of local Langlands parameters $\on{Loc}_{^L G}(D^\times)$,
with the categories $\ghat_\ka\mod_\chi$ being the ``fibers'' and the
group $G\ppart$ acting along these fibers. From this point of view the
categories $\ghat_\ka\mod_\chi$ should give us a ``spectral
decomposition'' of the category $\ghat_\ka\mod$ over $\on{Loc}_{^L
G}(D^\times)$.

In the next sections we will present a concrete proposal made in
\cite{FG:local} describing these categories in the special case when
$\ka=\ka_c$, the critical level.

\section{Center and opers}    \label{center and opers}

In Section 1 we have introduced the category $\ghat_\ka\mod$ whose
objects are smooth $\ghat_\ka$-modules on which the central element
${\mb 1}$ acts as the identity. As explained at the end of the
previous section, we wish to show that this category ``fibers'' over
the space of the Langlands parameters, which are gauge equivalence
classes of $^L G$-connections on the punctured disc $D^\times$ (or
perhaps, something similar). Moreover, the loop group $G\ppart$ should
act on this category ``along the fibers''.

Any abelian category may be thought of as ``fibering'' over the
spectrum of its center. Hence the first idea that comes to mind is to
describe the center of the category $\ghat_\ka\mod$ in the hope that
its spectrum is related to the Langlands parameters. As we will see,
this is indeed the case for a particular value of $\ka$.

\subsection{Center of an abelian category}

Let us first recall what is the center of an abelian category. Let
${\mc C}$ be an abelian category over $\C$. The center $Z({\mc C})$ is
by definition the set of endomorphisms of the identity functor on
${\mc C}$. Let us recall such such an endomorphism is a system of
endomorphisms $e_M \in \on{Hom}_{\mc C}(M,M)$, for each object $M$ of
${\mc C}$, which is compatible with the morphisms in ${\mc C}$: for
any morphism $f: M \to N$ in ${\mc C}$ we have $f \circ e_M = e_N
\circ f$. It is clear that $Z({\mc C})$ has a natural structure of a
commutative algebra over $\C$.

Let $S = \on{Spec} Z({\mc C})$. This is an affine algebraic variety
such that $Z({\mc C})$ is the algebra of functions on $S$. Each point
$s \in S$ defines an algebra homomorphism (equivalently, a character)
$\rho_s: Z({\mc C}) \to \C$ (evaluation of a function at the point
$s$). We define the full subcategory ${\mc C}_s$ of ${\mc C}$ whose
objects are the objects of ${\mc C}$ on which $Z({\mc C})$ acts
according to the character $\rho_s$. It is instructive to think
of the category ${\mc C}$ as ``fibering'' over $S$, with the fibers
being the categories ${\mc C}_s$.

Now suppose that ${\mc C} = A\mod$ is the category of left modules
over an associative $\C$-algebra $A$. Then $A$ itself, considered as a
left $A$-module, is an object of ${\mc C}$, and so we obtain a
homomorphism $$Z({\mc C}) \to Z(\on{End}_A A) = Z(A^{\on{opp}}) =
Z(A),$$ where $Z(A)$ is the center of $A$. On the other hand, each
element of $Z(A)$ defines an endomorphism of each object of $A\mod$,
and so we obtain a homomorphism $Z(A) \to Z({\mc C})$. It is easy to
see that these maps set mutually inverse isomorphisms between $Z({\mc
C})$ and $Z(A)$.

If $\g$ is a Lie algebra, then the category $\g\mod$ of $\g$-modules
coincides with the category $U(\g)\mod$ of $U(\g)$-modules, where
$U(\g)$ is the universal enveloping algebra of $\g$. Therefore the
center of the category $\g\mod$ is equal to the center of $U(\g)$,
which by abuse of notation we denote by $Z(\g)$.

Now consider the category $\ghat_\ka\mod$. Let
us recall from \secref{unram first} that objects of $\ghat_\ka\mod$ are
$\ghat_\ka$-modules $M$ on which the central element ${\mb 1}$ acts as
the identity and which are {\bf smooth}, that is for any vector $v \in
M$ we have
\begin{equation}    \label{smooth1}
(\g \otimes t^N \C[[t]]) \cdot v = 0
\end{equation}
for sufficiently large $N$.

Thus, we see that there are two properties that its objects
satisfy. Therefore it does not coincide with the category of all
modules over the universal enveloping algebra $U(\ghat_\ka)$ (which is
the category of all $\ghat_\ka$-modules). We need to modify this
algebra.

First of all, since ${\mb 1}$ acts as the identity, the action of
$U(\ghat_\ka)$ factors through the quotient
$$
U_\ka(\ghat) \overset{\on{def}}= U_\ka(\ghat)/({\mb 1} - 1).
$$
Second, the smoothness condition \eqref{smooth1} implies that the
action of $U_\ka(\ghat)$ extends to an action of its completion defined
as follows.

Define a linear topology on $U_\ka(\ghat)$ by using as the
basis of neighborhoods for $0$ the following left ideals:
$$I_N = U_\ka(\ghat) (\ggg\tensor t^N\C[[t]]), \qquad
N\ge 0.$$ Let $\wt{U}_\ka(\ghat)$ be the completion of $U_\ka(\ghat)$
with respect to this topology. Note that, equivalently, we can write
$$
\wt{U}_\ka(\ghat) = \lim_{\longleftarrow} U_\ka(\ghat)/I_N.
$$
Even though the $I_N$'s are only left ideals (and not two-sided
ideals), one checks that the associative product structure on
$U_\ka(\ghat)$ extends by continuity to an associative product
structure on $\wt{U}_\ka(\ghat)$ (this follows from the fact that the
Lie bracket on $U_\ka(\ghat)$ is continuous in the above
topology). Thus, $\wt{U}_\ka(\ghat)$ is a complete topological
algebra. It follows from the definition that the category
$\ghat_\ka\mod$ coincides with the category of discrete modules over
$\wt{U}_\ka(\ghat)$ on which the action of $\wt{U}_\ka(\ghat)$ is
pointwise continuous (this is precisely equivalent to the
condition \eqref{smooth1}).

It is now easy to see that the center of our category $\ghat_\ka\mod$
is equal to the center of the algebra $\wt{U}_\ka(\ghat)$, which we
will denote by $Z_\ka(\ghat)$. The argument is similar to the one we
used above: though $\wt{U}_\ka(\ghat)$ itself is not an object of
$\ghat_\ka\mod$, we have a collection of objects
$\wt{U}_\ka(\ghat)/I_N$. Using this collection, we obtain an
isomorphism between the center of category $\ghat_\ka\mod$ and the
inverse limit of the algebras $Z(\on{End}_{\ghat_\ka}
\wt{U}_\ka(\ghat)/I_N)$, which, by definition, coincides with
$Z_\ka(\ghat)$.

Now we can formulate our first question:

\medskip

\begin{center}

{\em describe the center $Z_\ka(\ghat)$ for all levels $\ka$.}

\end{center}

\medskip

In order to answer this question we need to introduce the concept of
$G$-opers.

\subsection{Opers}    \label{opers1}

Let $G$ be a simple algebraic group of adjoint type, $B$ its Borel
subgroup and $N = [B,B]$ its unipotent radical, with the corresponding
Lie algebras $\n \subset \bb\subset \g$.

Thus, $\g$ is a simple Lie algebra, and as such it has the Cartan
decomposition
$$\ggg = \nn_- \oplus \hh \oplus \nn_+.$$ We will choose generators
$e_1, \dots, e_\ell$ (resp., $f_1, \dots, f_\ell$) of $\n_+$ (resp.,
$\n_-$). We have $\n_{\al_i} = \C e_i, \n_{-\al_i} = \C f_i$. We take
$\bb = \hh \oplus \nn_+$ as the Lie algebra of $B$. Then $\n$ is the
Lie algebra of $N$. In what follows we will use the notation $\n$ for
$\n_+$.

Let $[\n,\n]^\perp \subset \g$ be the orthogonal complement of
$[\n,\n]$ with respect to a non-degenerate invariant bilinear form
$\ka_0$. We have
$$
[\n,\n]^\perp/\bb \simeq \bigoplus_{i=1}^\ell \n_{-\al_i}.
$$
Clearly, the group $B$ acts on $\n^\perp/\bb$. Our first observation
is that there is an open $B$-orbit ${\bf O}\subset \n^\perp/\bb
\subset \g/\bb$, consisting of vectors whose projection on each
subspace $\n_{-\al_i}$ is non-zero. This orbit may also be described
as the $B$-orbit of the sum of the projections of the generators $f_i,
i=1,\ldots,\ell$, of any possible subalgebra $\n_-$, onto
$\g/\bb$. The action of $B$ on ${\bf O}$ factors through an action of
$H = B/N$. The latter is simply transitive and makes ${\bf O}$ into an
$H$-torsor.

Let $X$ be a smooth curve and $x$ a point of $X$. As before, we denote
by $\O_x$ the completed local ring and by $\K_x$ its field of
fractions. The ring $\O_x$ is isomorphic, but not canonically, to
$\C[[t]]$. Then $D_x = \on{Spec} \OO_x$ is the disc without a
coordinate and $D_x^\times = \on{Spec} \K_x$ is the corresponding
punctured disc.

Suppose now that we are given a principal $G$-bundle $\F$ on a smooth
curve $X$, or $D_x$, or $D_x^\times$, together with a connection
$\nabla$ (automatically flat) and a reduction $\F_B$ to the Borel
subgroup $B$ of $G$. Then we define the relative position of $\nabla$
and $\F_B$ (i.e., the failure of $\nabla$ to preserve $\F_B$) as
follows. Locally, choose any flat connection $\nabla'$ on $\F$
preserving $\F_B$, and take the difference $\nabla - \nabla'$, which
is a section of $\g_{\F_B} \otimes \Omega_X$. We project it onto
$(\g/\bb)_{\F_B} \otimes \Omega_X$. It is clear that the resulting
local section of $(\g/\bb)_{\F_B} \otimes \Omega_X$ are independent of
the choice $\nabla'$. These sections patch together to define a global
$(\g/\bb)_{\F_B}$-valued one-form on $X$, denoted by $\nabla/\F_B$.

Let $X$ be a smooth curve, or $D_x$, or $D_x^\times$. Suppose we are
given a principal $G$-bundle $\F$ on $X$, a connection $\nabla$ on
$\F$ and a $B$-reduction $\F_B$. We will say that $\F_B$ is {\bf
transversal} to $\nabla$ if the one-form $\nabla/\F_B$ takes values in
${\bf O}_{\F_B}
\subset(\g/\bb)_{\F_B}$. Note that ${\bf O}$ is $\C^\times$-invariant,
so that ${\bf O} \otimes \Omega_X$ is a well-defined subset of
$(\g/\bb)_{\F_B} \otimes \Omega_X$.

Now, a $G$-{\bf oper} on $X$ is by definition a triple
$(\F,\nabla,\F_B)$, where $\F$ is a principal $G$-bundle $\F$ on $X$,
$\nabla$ is a connection on $\F$ and $\F_B$ is a $B$-reduction of
$\F$, such that $\F_B$ is transversal to $\nabla$.

\medskip

This definition is due to A. Beilinson and V. Drinfeld \cite{BD} (in
the case when $X$ is the punctured disc opers were introduced earlier
by V. Drinfeld and V. Sokolov in \cite{DS}).

\medskip

Equivalently, the transversality condition may be reformulated as
saying that if we choose a local trivialization of $\F_B$ and a local
coordinate $t$ then the connection will be of the form
\begin{equation}    \label{form of nabla1}
\nabla = \pa_t + \sum_{i=1}^\ell \psi_i(t) f_i + {\mb v}(t),
\end{equation}
where each $\psi_i(t)$ is a nowhere vanishing function, and ${\mb
v}(t)$ is a $\bb$-valued function.

If we change the trivialization of $\F_B$, then this operator will get
transformed by the corresponding $B$-valued gauge transformation. This
observation allows us to describe opers on the disc $D_x = \on{Spec}
\O_x$ and the punctured disc $D_x^\times = \on{Spec} \K_x$ in a more
explicit way. The same reasoning will work on any sufficiently small
analytic subset $U$ of any curve, equipped with a local coordinate
$t$, or on a Zariski open subset equipped with an \'etale
coordinate. For the sake of definiteness, we will consider now the
case of the base $D_x^\times$.

Let us choose a coordinate $t$ on $D_x$, i.e., an isomorphism $\O_x
\simeq \C[[t]]$. Then we identify $D_x$ with $D = \on{Spec} \C[[t]]$
and $D_x^\times$ with $D^\times = \on{Spec} \C\ppart$. The space
$\on{Op}_G(D^\times)$ of $G$-opers on $D^\times$ is then the quotient
of the space of all operators of the form \eqref{form of nabla1},
where $\psi_i(t) \in \C\ppart, \psi_i(0) \neq 0, i=1,\ldots,\ell$, and
${\mb v}(t) \in \bb\ppart$, by the action of the group $B\ppart$ of
gauge transformations:
$$
g \cdot (\pa_t + A(t)) = \pa_t + g A(t) g^{-1} - g^{-1}\pa_t g.
$$

Let us choose a splitting $\imath: H \to B$ of the homomorphism $B \to
H$. Then $B$ becomes the semi-direct product $B = H \ltimes N$. The
$B$-orbit ${\bf O}$ is an $H$-torsor, and so we can use $H$-valued
gauge transformations to make all functions $\psi_i(t)$ equal to
$1$. In other words, there is a unique element of $H\ppart$, namely,
the element $\prod_{i=1}^\ell \check\omega_i(\psi_i(t))$, where
$\check\omega_i: \C^\times \to H$ is the $i$th fundamental coweight of
$G$, such that the corresponding gauge transformation brings our
connection operator to the form
\begin{equation}    \label{another form of nabla1}
\nabla = \pa_t + \sum_{i=1}^\ell f_i + {\mb v}(t), \qquad {\mb v}(t)
\in \bb\ppart.
\end{equation}
What remains is the group of $N$-valued gauge transformations. Thus,
we obtain that $\on{Op}_G(D^\times)$ is equal to the quotient of the
space $\wt{\on{Op}}_G(D^\times)$ of operators of the form
\eqref{another form of nabla1} by the action of the group $N\ppart$ by
gauge transformations:
$$
\on{Op}_G(D^\times) = \wt{\on{Op}}_G(D^\times)/N\ppart.
$$

\begin{lem}[\cite{DS}]    \label{free1}
The action of $N\ppart$ on $\wt{\on{Op}}_G(D^\times)$ is free.
\end{lem}

\subsection{Canonical representatives}    \label{canon represent}

Now we construct canonical representatives in the $N\ppart$-gauge
classes of connections of the form \eqref{another form of nabla1},
following \cite{BD}. Observe that the operator $\on{ad} \crho$ defines
a gradation on $\g$, called the {\bf principal gradation}, with
respect to which we have a direct sum decomposition $\g = \bigoplus_i
\g_i$. In particular, we have $\bb = \bigoplus_{i\geq 0} \bb_i$, where
$\bb_0 = \h$.

Let now
$$
p_{-1} = \sum_{i=1}^\ell f_i.
$$
The operator $\on{ad} p_{-1}$ acts from $\bb_{i+1}$ to $\bb_{i}$
injectively for all $i\geq 0$. Hence we can find for each $i\geq 0$ a
subspace $V_i \subset \bb_i$, such that $\bb_i = [p_{-1},\bb_{i+1}]
\oplus V_i$. It is well-known that $V_i \neq 0$ if and only if $i$ is
an {\bf exponent} of $\g$, and in that case $\dim V_i$ is equal to the
multiplicity of the exponent $i$. In particular, $V_0=0$.

Let $V = \bigoplus_{i\in E} V_i \subset \n$, where $E = \{
d_1,\ldots,d_\ell \}$ is the set of exponents of $\g$ counted with
multiplicity. They are equal to the orders of the generators of the
center of $U(\g)$ minus one.  We note that the multiplicity of each
exponent is equal to one in all cases except the case $\g=D_{2n},
d_n=2n$, when it is equal to two.

There is a special choice of the transversal subspace $V =
\bigoplus_{i \in E} V_i$. Namely, there exists a unique element $p_1$
in $\n$, such that $\{ p_{-1},2\rv,p_1 \}$ is an $\sw_2$-triple. This
means that they have the same relations as the generators $\{ e,h,f
\}$ of $\liesl_2$. We have $p_1 = \sum_{i=1}^\ell m_i e_i$, where
$e_i$'s are generators of $\n_+$ and $m_i$ are certain coefficients
uniquely determined by the condition that $\{ p_{-1},2\rv,p_1 \}$ is
an $\sw_2$-triple.

Let $V^{\can} = \bigoplus_{i \in E} V^{\can}_{i}$ be the space of $\on{ad}
p_1$-invariants in $\n$. Then $p_1$ spans $V^{\can}_{1}$. Let $p_j$
be a linear generator of $V^{\can}_{d_j}$. If the multiplicity of $d_j$
is greater than one, then we choose linearly independent vectors in
$V^{\can}_{d_j}$.

Each $N\ppart$-equivalence class contains a unique operator of the
form $\nabla = \pa_t + p_{-1} + {\mathbf v}(t)$, where ${\mathbf v}(t)
\in V^{\can}[[t]]$, so that we can write
$$
{\mathbf v}(t) = \sum_{j=1}^\ell v_j(t) \cdot p_j, \qquad v_j(t) \in
\C[[t]].
$$
It is easy to find (see, e.g., \cite{F:rev}, Sect. 8.3) that under
changes of coordinate $t$, $v_1$ transforms as a projective
connection, and $v_j, j>1$, transforms as a $(d_j+1)$-differential on
$D_x$. Thus, we obtain an isomorphism
\begin{equation}    \label{repr1}
\on{Op}_G(D^\times) \simeq  {\mc P}roj(D^\times) \times
\bigoplus_{j=2}^\el \Omega_{\K}^{\otimes(d_j+1)},
\end{equation}
where $\Omega_{\K}^{\otimes n}$ is the space of $n$-differentials on
$D^\times$ and ${\mc P}roj(D^\times)$ is the $\Omega_{\K}^{\otimes
2}$-torsor of projective connections on $D^\times$.

We have an analogous isomorphism with $D^\times$ replaced by formal
disc $D$ or any smooth algebraic curve $X$.

\subsection{Description of the center}

Now we are ready to describe the center of the completed universal
enveloping algebra $\wt{U}_{\ka_c}(\ghat)$. The following assertion is
proved in \cite{newbook}, using results of \cite{Kac:laplace}:

\begin{prop}    \label{center away}
The center of $\wt{U}_{\ka}(\ghat)$ consists of the scalars for $\ka
\neq \ka_c$.
\end{prop}

Let us denote the center of $\wt{U}_{\ka_c}(\ghat)$ by $Z(\ghat)$. The
following theorem was proved in \cite{FF:gd,F:wak} (it was conjectured
by V. Drinfeld).

\begin{thm}    \label{Z ghat}
The center $Z(\ghat)$ is isomorphic to the algebra $\on{Fun}
\on{Op}_{^L G}(D^\times)$ in a way compatible with the action of the
group of coordinate changes.
\end{thm}

This implies the following result. Let $x$ be a point of a smooth
curve $X$. Then we have the affine algebra $\ghat_{\ka_c,x}$ as
defined in \secref{unram first} and the corresponding completed
universal enveloping algebra of critical level. We denote its center
by $Z(\ghat_x)$.

\begin{cor}
The center $Z(\ghat_x)$ is isomorphic to the algebra $\on{Fun}
\on{Op}_{^L G}(D^\times_x)$ of functions on the space of $^L G$-opers
on $D_x^\times$.
\end{cor}

\section{Opers vs. local systems}    \label{opers and loc syst}

We now go back to the question posed at the end of \secref{llc}:
let
\begin{equation}    \label{Loc as quotient}
\Loc_{^L G}(D^\times) = \left\{ \pa_t + A(t), \; A(t) \in {}^L \g\ppart
\right\}/\; {}^L G\ppart
\end{equation}
be the set of gauge equivalence classes of $^L G$-connections on the
punctured disc $D^\times = \on{Spec} \C\ppart$. We had argued in
\secref{llc} that $\Loc_{^L G}(D^\times)$ should be taken as the
space of Langlands parameters for the loop group $G\ppart$. Recall
that the loop group $G\ppart$ acts on the category $\ghat_{\ka}\mod$
of (smooth) $\ghat$-modules of level $\ka$ (see \secref{unram first}
for the definition of this category). We asked the following question:

\medskip

\noindent{\em Associate to each local Langlands parameter $\sigma \in
\on{Loc}_{^L G}(D^\times)$ a subcategory $\ghat_\ka\mod_\sigma$ of
$\ghat_\ka\mod$ which is stable under the action of the loop group
$G\ppart$.}

\medskip

Even more ambitiously, we wish to represent the category
$\ghat_\ka\mod$ as ``fibering'' over the space of local Langlands
parameters $\on{Loc}_{^L G}(D^\times)$, with the categories
$\ghat_\ka\mod_\sigma$ being the ``fibers'' and the group $G\ppart$
acting along these fibers. If we could do that, then we would think of
this fibration as a ``spectral decomposition'' of the category
$\ghat_\ka\mod$ over $\on{Loc}_{^L G}(D^\times)$.

At the beginning of \secref{center and opers} we proposed a
possible scenario for solving this problem. Namely, we observed that
any abelian category may be thought of as ``fibering'' over the
spectrum of its center. Hence our idea was to describe the center of
the category $\ghat_\ka\mod$ (for each value of $\ka$) and see if its
spectrum is related to the space $\Loc_{^L G}(D^\times)$ of Langlands
parameters.

We have identified the center of the category $\ghat_\ka\mod$ with the
center $Z_\ka(\ghat)$ of the associative algebra $\wt{U}_\ka(\ghat)$,
the completed enveloping algebra of $\ghat$ of level $\ka$, defined in
\secref{center and opers}. Next, we described the algebra
$Z_\ka(\ghat)$. According to \propref{center away}, if $\ka \neq
\ka_c$, the critical level, then $Z_\ka(\ghat) = \C$. Therefore our
approach cannot work for $\ka \neq \ka_c$. However, we found that the
center $Z_{\ka_c}(\ghat)$ at the critical level is highly non-trivial
and indeed related to $^L G$-connections on the punctured disc.

Now, following the works \cite{FG:exact}--\cite{FG:weyl} of
D. Gaitsgory and myself, I will use these results to formulate more
precise conjectures on the local Langlands correspondence for loop
groups and to provide some evidence for these conjectures. I will then
discuss the implications of these conjectures for the global geometric
Langlands correspondence.\footnote{Note that A. Beilinson has another
proposal \cite{Be} for local geometric Langlands correspondence, using
representations of affine Kac-Moody algebras of levels {\em less than
critical}. It would be interesting to understand the connection
between his proposal and ours.}

According to \thmref{Z ghat}, $Z_{\ka_c}(\ghat)$ is isomorphic to the
algebra $\on{Fun} \on{Op}_{^L G}(D^\times)$ of functions on the space
of $^L G$-opers on the punctured disc $D^\times$. This isomorphism is
compatible with various symmetries and structures on both algebras,
such as the action of the group of coordinate changes. There is a
one-to-one correspondence between points $\chi \in \on{Op}_{^L
G}(D^\times)$ and homomorphisms (equivalently, characters) $$\on{Fun}
\on{Op}_{^L G}(D^\times) \to \C,$$ corresponding to evaluating a
function at $\chi$. Hence points of $\on{Op}_{^L G}(D^\times)$
parametrize {\bf central characters} $Z_{\ka_c}(\ghat) \to \C$.

Given a $^L G$-oper $\chi \in \on{Op}_{^L G}(D^\times)$, define the
category $$\gmod_\chi$$ as a full subcategory of $\gmod$ whose
objects are $\ghat$-modules of critical level (hence
$\wt{U}_{\ka_c}(\ghat)$-modules) on which the center $Z_{\ka_c}(\ghat)
\subset \wt{U}_{\ka_c}(\ghat)$ acts according to the central character
corresponding to $\chi$. More generally, for any closed algebraic
subvariety $Y \subset \on{Op}_{^L G}(D^\times)$ (not necessarily a
point), we have an ideal
$$
I_Y \subset \on{Fun} \on{Op}_{^L G}(D^\times) \simeq Z_{\ka_c}(\ghat)
$$
of those functions that vanish on $Y$. We then have a full subcategory
$\gmod_Y$ of $\gmod$ whose objects are $\ghat$-modules of critical
level on which $I_Y$ acts by $0$. This category is an example of a
``base change'' of the category $\gmod$ with respect to the morphism
$Y \to \on{Op}_{^L G}(D^\times)$. It is easy to generalize this
definition to an arbitrary affine scheme $Y$ equipped with a morphism
$Y \to \on{Op}_{^L G}(D^\times)$.\footnote{The corresponding base
changed categories $\gmod_Y$ may then be ``glued'' together, which
allows us to define the base changed category  $\gmod_Y$ for any
scheme $Y$ mapping to $\on{Op}_{^L G}(D^\times)$.}

Since the algebra $\on{Op}_{^L G}(D^\times)$ acts on the category
$\gmod$, one can say that the category $\ghat_\ka\mod$ ``fibers'' over
the space $\on{Op}_{^L G}(D^\times)$, in such a way that the
fiber-category corresponding to $\chi \in \on{Op}_{^L G}(D^\times)$ is
the category $\gmod_\chi$.\footnote{The precise notion of an abelian
category fibering over a scheme is spelled out in \cite{Ga}.}

Recall that the group $G\ppart$ acts on $\wt{U}_{\ka_c}(\ghat)$ and on
the category $\gmod$. One can show (see \cite{BD}, Remark 3.7.11(iii))
that the action of $G\ppart$ on $Z_{\ka_c}(\ghat) \subset
\wt{U}_{\ka_c}(\ghat)$ is trivial. Therefore the subcategories
$\gmod_\chi$ (and, more generally, $\gmod_Y$) are stable under the
action of $G\ppart$. Thus, the group $G\ppart$ acts ``along the
fibers'' of the ``fibration'' $\gmod \to \on{Op}_{^L G}(D^\times)$
(see \cite{FG:local}, Sect. 20, for more details).

The fibration $\gmod \to \on{Op}_{^L G}(D^\times)$ almost gives us
the desired local Langlands correspondence for loop groups. But there
is one important difference: we asked that the category $\gmod$ fiber
over the space $\Loc_{^L G}(D^\times)$ of local systems on
$D^\times$. We have shown, however, that $\gmod$ fibers over the space
$\on{Op}_{^L G}(D^\times)$ of $^L G$-opers.

What is the difference between the two spaces? While a $^L G$-local
system is a pair $({\mc F},\nabla)$, where ${\mc F}$ is an $^L
G$-bundle and $\nabla$ is a connection on ${\mc F}$, an $^L G$-oper is
a triple $({\mc F},\nabla,{\mc F}_{^L B})$, where ${\mc F}$ and
$\nabla$ are as before, and ${\mc F}_{^L B}$ is an additional piece of
structure, namely, a reduction of ${\mc F}$ to a (fixed) Borel
subgroup $^L B \subset {} ^L G$ satisfying the transversality
condition explained in \secref{opers1}. Thus, for any curve $X$ we
clearly have a forgetful map
$$
\on{Op}_{^L G}(X) \to \Loc_{^L G}(X).
$$
The fiber of this map over $({\mc F},\nabla) \in \Loc_{^L G}(X)$
consists of all $^L B$-reductions of ${\mc F}$ satisfying the
transversality condition with respect to $\nabla$.

For a general $X$ it may well be that this map is not surjective,
i.e., that the fiber of this map over a particular local system $({\mc
F},\nabla)$ is empty. For example, if $X$ is a projective curve and
$^L G$ is a group of adjoint type, then there is a unique $^L
G$-bundle ${\mc F}_{^L G}$ such that the fiber over $({\mc F}_{^L
G},\nabla)$ is non-empty, as we saw in \secref{unram first}.

The situation is quite different when $X=D^\times$. In this case any
$^L G$-bundle ${\mc F}$ may be trivialized. A connection $\nabla$
therefore may be represented as a first order operator $\pa_t + A(t),
A(t) \in {}^L \g\ppart$. However, the trivialization of ${\mc F}$ is
not unique; two trivializations differ by an element of $^L
G\ppart$. Therefore the set of equivalence classes of pairs $({\mc
F},\nabla)$ is identified with the quotient \eqref{Loc as quotient}.

Suppose now that $({\mc F},\nabla)$ carries an oper reduction ${\mc
F}_{^L B}$. Then we consider only those trivializations of ${\mc F}$
which come from trivializations of ${\mc F}_{^L B}$. There are fewer
of those, since two trivializations now differ by an element of $^L
B\ppart$ rather than $^L G\ppart$. Due to the oper transversality
condition, the connection $\nabla$ must have a special form with
respect to any of those trivializations, namely,
$$
\nabla = \pa_t + \sum_{i=1}^\ell \psi_i(t) f_i + {\mb v}(t),
$$
where each $\psi_i(t) \neq 0$ and ${\mb
v}(t) \in {}^L \bb\ppart$ (see \secref{opers1}). Thus,
we obtain a concrete realization of the space of opers as a space of
gauge equivalence classes
\begin{equation}    \label{Op as quotient}
\on{Op}_{^L G}(D^\times) = \left. \left\{ \pa_t + \sum_{i=1}^\ell
\psi_i(t) f_i + {\mb v}(t), \; \psi_i \neq 0, {\mb v}(t) \in {}^L
\bb\ppart \right\} \right/ \; {}^L B\ppart.
\end{equation}
Now the map
$$
\al: \on{Op}_{^L G}(D^\times) \to \Loc_{^L G}(D^\times)
$$
simply takes a $^L B\ppart$-equivalence class of operators of the form
\eqref{Op as quotient} to its $^L G\ppart$-equivalence class.

Unlike the case of projective curves $X$ discussed above, we expect
that the map $\al$ is {\bf surjective} for any simple Lie group $^L
G$. In the case of $G=SL_n$ this follows from the results of
P. Deligne \cite{D}, and we conjecture it to be true in general.

\begin{conj}    \label{surjective map}
The map $\al$ is surjective for any simple Lie group $^L G$.
\end{conj}

Now we find ourselves in the following situation: we {\em expect}
that there exists a category ${\mc C}$ fibering over the space
$\Loc_{^L G}(D^\times)$ of ``true'' local Langlands parameters,
equipped with a fiberwise action of the loop group $G\ppart$. The
fiber categories ${\mc C}_\sigma$ corresponding to various $\sigma \in
\Loc_{^L G}(D^\times)$ should satisfy various, not yet specified,
properties. This should be the ultimate form of the local Langlands
correspondence. On the other hand, we have {\em constructed} a
category $\gmod$ which fibers over a close cousin of the space
$\Loc_{^L G}(D^\times)$, namely, the space $\on{Op}_{^L G}(D^\times)$
of $^L G$-opers, and is equipped with a fiberwise action of the loop
group $G\ppart$.

What should be the relationship between the two?

The idea of \cite{FG:local} is that the second fibration is
a ``base change'' of the first one, that is we have a Cartesian
diagram
\begin{equation}    \label{base change}
\begin{CD}
\gmod @>>> {\mc C} \\
@VVV @VVV \\
\on{Op}_{^L G}(D^\times) @>{\al}>> \Loc_{^L G}(D^\times)
\end{CD}
\end{equation}
that commutes with the action of $G\ppart$ along the fibers of the two
vertical maps. In other words,
$$
\gmod \simeq {\mc C} \underset{\Loc_{^L G}(D^\times)}\times
\on{Op}_{^L G}(D^\times).
$$
At present, we do not have a definition of ${\mc C}$, and therefore we
cannot make this isomorphism precise. But we will use it as our
guiding principle. We will now discuss various corollaries of this
conjecture and various pieces of evidence that make us believe that it
is true.

In particular, let us fix a Langlands parameter $\sigma \in \Loc_{^L
G}(D^\times)$ that is in the image of the map $\al$ (according to
\conjref{surjective map}, all Langlands parameters are). Let $\chi$
be a $^L G$-oper in the preimage of $\sigma$, $\al^{-1}(\sigma)$. Then,
according to the above conjecture, the category $\gmod_\chi$ is
equivalent to the ``would be'' Langlands category ${\mc C}_\sigma$
attached to $\sigma$. Hence we may take $\gmod_\chi$ as the {\bf
definition} of ${\mc C}_\sigma$.

The caveat is, of course, that we need to ensure that this definition
is independent of the choice of $\chi$ in $\al^{-1}(\sigma)$. This
means that for any two $^L G$-opers, $\chi$ and $\chi'$, in the
preimage of $\sigma$, the corresponding categories, $\gmod_\chi$ and
$\gmod_{\chi'}$, should be equivalent to each other, and this
equivalence should commute with the action of the loop group
$G\ppart$. Moreover, we should expect that these equivalences are
compatible with each other as we move along the fiber
$\al^{-1}(\sigma)$. We will not try to make this condition more
precise here (however, we will explain below in \conjref{equiv reg
opers} what this means for regular opers).

Even putting the questions of compatibility aside, we arrive at the
following rather non-trivial conjecture (see \cite{FG:local}).

\begin{conj}    \label{equiv loc syst}
Suppose that $\chi,\chi' \in \on{Op}_{^L G}(D^\times)$ are such
that $\al(\chi) = \al(\chi')$, i.e., that the flat $^L G$-bundles
on $D^\times$ underlying the $^L G$-opers $\chi$ and $\chi'$ are
isomorphic to each other. Then there is an equivalence between the
categories $\gmod_\chi$ and $\gmod_{\chi'}$ which commutes with
the actions of the group $G\ppart$ on the two categories.
\end{conj}

Thus, motivated by our quest for the local Langlands correspondence,
we have found an unexpected symmetry in the structure of the category
$\gmod$ of $\ghat$-modules of critical level.

\section{Harish--Chandra categories}

As explained in \secref{llc}, the local Langlands correspondence
for the loop group $G\ppart$ should be viewed as a categorification of
the local Langlands correspondence for the group $G(F)$, where $F$ is
a local non-archimedian field. This means that the categories ${\mc
C}_\sigma$, equipped with an action of $G\ppart$, that we are trying
to attach to the Langlands parameters $\sigma \in \Loc_{^L
G}(D^\times)$ should be viewed as categorifications of the smooth
representations of $G(F)$ on complex vector spaces attached to the
corresponding local Langlands parameters discussed in \secref{dual
grp}. Here we use the term ``categorification'' to indicate that we
expect the Grothendieck groups of the categories $\CC_\sigma$ to
``look like'' irreducible smooth representations of $G(F)$. We begin
by taking a closer look at the structure of these representations.

\subsection{Spaces of $K$-invariant vectors}    \label{Kinv}

It is known that an irreducible smooth representation $(R,\pi)$ of
$G(F)$ is automatically {\bf admissible}, in the sense that for any
open compact subgroup $K$, such as the $N$th congruence subgroup $K_N$
defined in \secref{langlands param}, the space $R^{\pi(K)}$ of
$K$-invariant vectors in $R$ is finite-dimensional. Thus, while most
of the irreducible smooth representations $(R,\pi)$ of $G(F)$ are
infinite-dimensional, they are filtered by the finite-dimensional
subspaces $R^{\pi(K)}$ of $K$-invariant vectors, where $K$ are smaller
and smaller open compact subgroups. The space $R^{\pi(K)}$ does not
carry an action of $G(F)$, but it carries an action of the {\bf Hecke
algebra} $H(G(F),K)$.

By definition, $H(G(F),K)$ is the space of compactly supported
$K$ bi-invariant functions on $G(F)$. It is given an algebra structure
with respect to the {\bf convolution product}
\begin{equation}    \label{conv}
(f_1 \star f_2)(g) = \int_{G(F)} f_1(gh^{-1}) f_2(h) \; dh,
\end{equation}
where $dh$ is the Haar measure on $G(F)$ normalized in such a way
that the volume of the subgroup $K_0 = G(\OO)$ is equal to $1$ (here
$\OO$ is the ring of integers of $F$; e.g., for $F = \Fq\ppart$ we
have $\OO = \Fq[[t]]$). The algebra $H(G(F),K)$ acts on the space
$R^{\pi(K)}$ by the formula
\begin{equation}    \label{conv action}
f \star v = \int_{G(F)} f_1(gh^{-1}) (\pi(h) \cdot v) \; dh, \qquad v
\in R^{\pi(K)}.
\end{equation}

Studying the spaces of $K$-invariant vectors and their ${\mc
H}(G(F),K)$-module structure gives us an effective tool for analyzing
representations of the group $G(F)$, where $F = \Fq\ppart$.

Can we find a similar structure in the categorical local Langlands
correspondence for loop groups?

\subsection{Equivariant modules}    \label{eq mod}

In the categorical setting a representation $(R,\pi)$ of the group
$G(F)$ is replaced by a category equipped with an action of $G\ppart$,
such as $\gmod_\chi$. The open compact subgroups of $G(F)$ have
obvious analogues for the loop group $G\ppart$ (although they are, of
course, not compact with respect to the usual topology on
$G\ppart$). For instance, we have the ``maximal compact subgroup''
$K_0 = G[[t]]$, or, more generally, the $N$th congruence subgroup
$K_N$, whose elements are congruent to $1$ modulo
$t^N\C[[t]]$. Another important example is the analogue of the {\bf
Iwahori subgroup}. This is the subgroup of $G[[t]]$, which we denote
by $I$, whose elements $g(t)$ have the property that their value at
$0$, that is $g(0)$, belong to a fixed Borel subgroup $B \subset G$.

Now, for a subgroup $K \subset G\ppart$ of this type, an analogue of a
$K$-invariant vector in the categorical setting is an object of our
category, i.e., a smooth $\ghat_{\ka_c}$-module $(M,\rho)$, where
$\rho: \ghat_{\ka_c} \to \on{End} M$, which is stable under the action
of $K$. Recall from \secref{back to groups} that for any $g \in G\ppart$
we have a new $\ghat_{\ka_c}$-module $(M,\rho_g)$, where $\rho_g(x) =
\rho(\on{Ad}_g(x))$. We say that $(M,\rho)$ is stable under $K$, or
that $(M,\rho)$ is {\bf weakly} $K$-{\bf equivariant}, if there is a
compatible system of isomorphisms between $(M,\rho)$ and $(M,\rho_k)$
for all $k \in K$. More precisely, this means that for each $k \in K$
there exists a linear map $T^M_k: M \to M$ such that $$T^M_k \rho(x)
(T^M_k)^{-1} = \rho(\on{Ad}_k(x))$$ for all $x \in \ghat_{\ka_c}$, and
we have
$$
T^M_1 = \on{Id}_M, \qquad T^M_{k_1} T^M_{k_2} = T^M_{k_1 k_2}.
$$

Thus, $M$ becomes a representation of the group $K$.\footnote{In
general, it is reasonable to modify the last condition to allow for a
non-trivial two-cocycle and hence a non-trivial central extension of
$K$; however, in the case of interest $K$ does not have any
non-trivial central extensions.} Consider the corresponding
representation of the Lie algebra ${\mathfrak k} = \on{Lie} K$ on
$M$. Let us assume that the embedding ${\mathfrak k} \hookrightarrow
\g\ppart$ lifts to ${\mathfrak k} \hookrightarrow \ghat_{\ka_c}$
(i.e., that the central extension cocycle is trivial on ${\mathfrak
k}$). This is true, for instance, for any subgroup contained in $K_0 =
G[[t]]$, or its conjugate. Then we also have a representation of
${\mathfrak k}$ on $M$ obtained by restriction of $\rho$. In general,
the two representations do not have to coincide. If they do coincide,
then the module $M$ is called {\bf strongly} $K$-{\bf equivariant}, or
simply $K$-{\bf equivariant}.

The pair $(\ghat_{\ka_c},K)$ is an example of {\bf Harish-Chandra
pair}, that is a pair $(\g,H)$ consisting of a Lie algebra $\g$ and a
Lie group $H$ whose Lie algebra is contained in $\g$. The
$K$-equivariant $\ghat_{\ka_c}$-modules are therefore called
$(\ghat_{\ka_c},K)$ {\bf Harish-Chandra modules}. These are (smooth)
$\ghat_{\ka_c}$-modules on which the action of the Lie algebra
$\on{Lie} K \subset \ghat_{\ka_c}$ may be exponentiated to an action
of $K$ (we will assume that $K$ is connected). We denote by $\gmod^K$
and $\gmod^K_\chi$ the full subcategories of $\gmod$ and $\gmod_\chi$,
respectively, whose objects are $(\ghat_{\ka_c},K)$ Harish-Chandra
modules.

We will stipulate that the analogues of $K$-invariant vectors in the
category $\gmod_\chi$ are $(\ghat_{\ka_c},K)$ Harish-Chandra
modules. Thus, while the categories $\gmod_\chi$ should be viewed as
analogues of smooth irreducible representations $(R,\pi)$ of the group
$G(F)$, the categories $\gmod^K_\chi$ are analogues of the spaces of
$K$-invariant vectors $R^{\pi(K)}$.

Next, we discuss the categorical analogue of the Hecke algebra
$H(G(F),K)$.

\subsection{Categorical Hecke algebras}    \label{cat hecke}

We recall that $H(G(F),K)$ is the algebra of compactly supported
$K$ bi-invariant functions on $G(F)$. We realize it as the algebra of
left $K$-invariant compactly supported functions on $G(F)/K$. In
\secref{finite groups} we have already discussed the question of
categorification of the algebra of functions on a homogeneous space
like $G(F)/K$. Our conclusion was that the categorical analogue of
this algebra, when $G(F)$ is replaced by the complex loop group
$G\ppart$, is the category of ${\mc D}$-modules on $G\ppart/K$. More
precisely, this quotient has the structure of an ind-scheme which is a
direct limit of finite-dimensional algebraic varieties with respect to
closed embeddings. The appropriate notion of (right) ${\mc D}$-modules
on such ind-schemes is formulated in \cite{BD} (see also
\cite{FG:exact,FG:local}). As the categorical analogue of the algebra
of left $K$-invariant functions on $G(F)/K$, we take the category
${\mc H}(G\ppart,K)$ of $K$-equivariant ${\mc D}$-modules on the
ind-scheme $G\ppart/K$ (with respect to the left action of $K$ on
$G\ppart/K$). We call it the {\bf categorical Hecke algebra}
associated to $K$.

It is easy to define the convolution of two objects of the category
${\mc H}(G\ppart,K)$ by imitating formula \eqref{conv}. Namely, we
interpret this formula as a composition of the operations of pulling
back and integrating functions. Then we apply the same operations to
${\mc D}$-modules, thinking of the integral as push-forward. However,
here one encounters two problems. The first problem is that for a
general group $K$ the morphisms involved will not be proper, and so we
have to choose between the $*$- and $!$-push-forward. This problem
does not arise, however, if $K$ is such that $I \subset K \subset
G[[t]]$, which will be our main case of interest. The second, and more
serious, issue is that in general the push-forward is not an exact
functor, and so the convolution of two ${\mc D}$-modules will not be a
${\mc D}$-module, but a complex, more precisely, an object of the
corresponding $K$-equivariant (bounded) derived category
$D^b(G\ppart/K)^K$ of ${\mc D}$-modules on $G\ppart/K$. We will not
spell out the exact definition of this category here, referring the
interested reader to \cite{BD} and
\cite{FG:local}. The exception is the case of the subgroup $K_0 =
G[[t]]$, when the convolution functor is exact and so we may restrict
ourselves to the abelian category of $K_0$-equivariant ${\mc
D}$-modules on $G\ppart/K_0$.

Now the category $D^b(G\ppart/K)^K$ has a monoidal structure, and as
such it acts on the derived category of $(\ghat_{\ka_c},K)$
Harish-Chandra modules (again, we refer the reader to
\cite{BD,FG:local} for the precise definition). In the special case
when $K=K_0$, we may restrict ourselves to the corresponding abelian
categories. This action should be viewed as the categorical analogue
of the action of $H(G(F),K)$ on the space $R^{\pi(K)}$ of
$K$-invariant vectors discussed above.

Our ultimate goal is understanding the ``local Langlands categories''
${\mc C}_\sigma$ associated to the ``local Langlands parameters $\sigma
\in \Loc_{^L G}(D^\times)$. We now have a candidate for the category
${\mc C}_\sigma$, namely, the category $\gmod_\chi$, where $\sigma =
\al(\chi)$. Therefore $\gmod_\chi$ should be viewed as a
categorification of a smooth representation $(R,\pi)$ of $G(F)$. The
corresponding category $\gmod^K_\chi$ of $(\ghat_{\ka_c},K)$
Harish-Chandra modules should therefore be viewed as a
categorification of $R^{\pi(K)}$. This category (or, more precisely,
its derived category) is acted upon by the categorical Hecke algebra
${\mc H}(G\ppart,K)$. We summarize this analogy in the following
table.

\smallskip

\begin{center}
\begin{tabular}{l|l}
\hline
& \\
{\bf Classical Theory} & {\bf Geometric Theory} \\
& \\
\hline
& \\
Representation of $G(F)$ & Representation of $G\ppart$ \\
on a vector space $R$ & on a category $\gmod_\chi$ \\
& \\
\hline
& \\
A vector in $R$ & An object of $\gmod_\chi$ \\
& \\
\hline
& \\
The subspace $R^{\pi(K)}$ of & The subcategory $\gmod^K_\chi$ of \\
$K$-invariant vectors of $R$ & $(\ghat_{\ka_c},K)$ Harish-Chandra
modules\\
& \\
\hline
& \\
Hecke algebra $H(G(F),K)$ & Categorical Hecke algebra ${\mc
  H}(G\ppart,K)$ \\
acts on $R^{\pi(K)}$ & acts on $\gmod^K_\chi$ \\
& \\
\hline
\end{tabular}
\end{center}

\smallskip

Now we may test our proposal for the local Langlands correspondence by
studying the categories $\gmod^K_\chi$ of Harish-Chandra modules and
comparing their structure to the structure of the spaces $R^{\pi(K)}$
of $K$-invariant vectors of smooth representations of $G(F)$ in the
known cases. Another possibility is to test \conjref{equiv loc syst}
when applied to the categories of Harish-Chandra modules.

In the next section we consider the case of the ``maximal compact
subgroup'' $K_0 = G[[t]]$ and find perfect agreement with the
classical results about unramified representations of $G(F)$. We then
take up the more complicated case of the Iwahori subgroup $I$. There
we also find the conjectures and results of \cite{FG:local} to be
consistent with the known results about representations of $G(F)$ with
Iwahori fixed vectors.

\section{Local Langlands correspondence: unramified case}
\label{unram case}

We first take up the case of the ``maximal compact subgroup'' $K_0 =
G[[t]]$ of $G\ppart$ and consider the categories $\gmod_\chi$ which
contain non-trivial $K_0$-equivariant objects.

\subsection{Unramified representations of $G(F)$}    \label{unr rep}

These categories are analogues of smooth representations of the group
$G(F)$, where $F$ is a local non-archimedian field (such as
$\Fq\ppart$) that contain non-zero $K_0$-invariant vectors. Such
representations are called {\bf unramified}. The classification of the
irreducible unramified representations of $G(F)$ is the simplest case
of the local Langlands correspondence discussed in Sections \ref{llc
for gln} and \ref{dual grp}. Namely, we have a bijection between the
sets of equivalence classes of the following objects:
\begin{equation}    \label{unramified langlands}
\boxed{\begin{matrix} \text{unramified admissible} \\
\text{homomorphisms } W'_F \to {}^L G \end{matrix}} \quad
\Longleftrightarrow  \quad \boxed{\begin{matrix} \text{irreducible
unramified} \\ \text{representations of } G(F)  \end{matrix}}
\end{equation}
where $W'_F$ is the Weil-Deligne group introduced in \secref{langlands
param}.

By definition, unramified homomorphisms $W'_F \longrightarrow {}^L G$
are those which factor through the quotient
$$
W'_F \to W_F \to \Z
$$
(see \secref{langlands param} for the definitions of these groups and
homomorphisms). It is admissible if its image in $^L G$ consists of
semi-simple elements. Therefore the set on the left hand side of
\eqref{unramified langlands} is just the set of conjugacy classes of
semi-simple elements of $^L G$. Thus, the above bijection may be
reinterpreted as follows:
\begin{equation}   \label{unramified langlands1}
\boxed{\begin{matrix} \text{semi-simple conjugacy} \\
\text{classes in } {}^L G \end{matrix}} \quad
\Longleftrightarrow  \quad \boxed{\begin{matrix} \text{irreducible
unramified} \\ \text{representations of } G(F)  \end{matrix}}
\end{equation}

To construct this bijection, we look at the Hecke algebra
$H(G(F),K_0)$. According to the Satake isomorphism \cite{Satake}, in
the interpretation of Langlands \cite{Langlands}, this algebra is
commutative and isomorphic to the representation ring of the Langlands
dual group $^L G$:
\begin{equation}    \label{satake}
H(G(F),K_0) \simeq \on{Rep} {}^L G.
\end{equation}
We recall that $\on{Rep} {}^L G$ consists of finite linear
combinations $\sum_i a_i [V_i]$, where the $V_i$ are
finite-dimensional representations of $^L G$ (without loss of
generality we may assume that they are irreducible) and $a_i \in \C$,
with respect to the multiplication
$$
[V] \cdot [W] = [V \otimes W].
$$
Because $\on{Rep} {}^L G$ is commutative, its irreducible modules are
all one-dimensional. They correspond to characters $\on{Rep} {}^L G
\to \C$. We have a bijection
\begin{equation}    \label{ss classes}
\boxed{\begin{matrix} \text{semi-simple conjugacy} \\
\text{classes in } {}^L G \end{matrix}} \quad
\Longleftrightarrow  \quad \boxed{\begin{matrix} \text{characters} \\
\text{of } \on{Rep} {}^L G \end{matrix}}
\end{equation}
where the character $\phi_\ga$ corresponding to the conjugacy class
$\ga$ is given by the formula\footnote{It is customary to multiply the
right hand side of this formula, for irreducible representation $V$,
by a scalar depending on $q$ and the highest weight of $V$, but this
is not essential for our discussion.}
$$
\phi_\ga: [V] \mapsto \on{Tr}(\ga,V).
$$

Now, if $(R,\pi)$ is a representation of $G(F)$, then the space
$R^{\pi(K_0)}$ of $K_0$-invariant vectors in $V$ is a module over
$H(G(F),K_0)$. It is easy to show that this sets up a one-to-one
correspondence between equivalence classes of irreducible unramified
representations of $G(F)$ and irreducible
$H(G(F),K_0)$-modules. Combining this with the bijection \eqref{ss
classes} and the isomorphism \eqref{satake}, we obtain the
sought-after bijections \eqref{unramified langlands} and
\eqref{unramified langlands1}.

In particular, we find that, because the Hecke algebra $H(G(F),K_0)$
is commutative, the space $R^{\pi(K_0)}$ of $K_0$-invariants of an
irreducible representation, which is an irreducible
$H(G(F),K_0)$-module, is either zero or one-dimensional. If it is
one-dimensional, then $H(G(F),K_0)$ acts on it by the character
$\phi_\ga$ for some $\ga$:
\begin{equation}    \label{Hecke eigenfunction}
H_V \star v = \on{Tr}(\ga,V) v, \qquad v \in R^{\pi(K_0)}, [V] \in
\on{Rep} {}^L G,
\end{equation}
where $H_V$ is the element of $H(G(F),K_0)$ corresponding to $[V]$
under the isomorphism \eqref{satake} (see formula \eqref{conv action}
for the definition of the convolution action).

We now discuss the categorical analogues of these statements.

\subsection{Unramified categories $\ghat_{\ka_c}$-modules}
\label{weyl}

In the categorical setting, the role of an irreducible representation
$(R,\pi)$ of $G(F)$ is played by the category $\gmod_\chi$ for some
$\chi \in \on{Op}_{^L G}(D^\times)$. The analogue of an unramified
representation is a category $\gmod_\chi$ which contains non-zero
$(\ghat_{\ka_c},G[[t]])$ Harish-Chandra modules. This leads us to the
following question: for what $\chi \in \on{Op}_{^L G}(D^\times)$ does
the category $\gmod_\chi$ contain non-zero $(\ghat_{\ka_c},G[[t]])$
Harish-Chandra modules?

We saw in the previous section that $(R,\pi)$ is unramified if and
only if it corresponds to an unramified Langlands parameter, which is
a homomorphism $W'_F \to {}^L G$ that factors through $W'_F \to
\Z$. Recall that in the geometric setting the Langlands parameters are
$^L G$-local systems on $D^\times$. The analogues of unramified
homomorphisms $W'_F \to {}^L G$ are those local systems on $D^\times$
which extend to the disc $D$, in other words, have no singularity at
the origin $0 \in D$. Note that there is a unique, up to isomorphism
local system on $D$. Indeed, suppose that we are given a regular
connection on a $^L G$-bundle ${\mc F}$ on $D$. Let us trivialize the
fiber ${\mc F}_0$ of ${\mc F}$ at $0 \in D$. Then, because $D$ is
contractible, the connection identifies ${\mc F}$ with the trivial
bundle on $D$. Under this identification the connection itself becomes
trivial, i.e., represented by the operator $\nabla = \pa_t$.

Therefore all regular $^L G$-local systems (i.e., those which extend
to $D$) correspond to a single point of the set $\Loc_{^L
G}(D^\times)$, namely, the equivalence class of the trivial local
system $\sigma_0$.\footnote{Note however that the trivial $^L G$-local
system on $D$ has a non-trivial group of automorphisms, namely, the
group $^L G$ itself (it may be realized as the group of automorphisms
of the fiber at $0 \in D$). Therefore if we think of $\Loc_{^L
G}(D^\times)$ as a stack rather than as a set, then the trivial local
system corresponds to a substack $\on{pt}/{}^L G$.}  {}From the point
of view of the realization of $\Loc_{^L G}(D^\times)$ as the quotient
\eqref{Loc} this simply means that there is a unique $^L G\ppart$
gauge equivalence class containing all regular connections of the form
$\pa_t + A(t)$, where $A(t) \in {}^L \g[[t]]$.

The gauge equivalence class of regular connections is the unique local
Langlands parameter that we may view as unramified in the geometric
setting. Therefore, by analogy with the unramified Langlands
correspondence for $G(F)$, we expect that the category
$\gmod_\chi$ contains non-zero $(\ghat_{\ka_c},G[[t]])$
Harish-Chandra modules if and only if the $^L G$-oper $\chi \in
\on{Op}_{^L G}(D^\times)$ is $^L G\ppart$ gauge equivalent to the
trivial connection, or, in other words, $\chi$ belongs to the fiber
$\al^{-1}(\sigma_0)$ over $\sigma_0$.

What does this fiber look like? Let $P^+$ be the set of dominant
integral weights of $G$ (equivalently, dominant integral coweights of
$^L G$). In \cite{FG:local} we defined, for each $\la \in P^+$, the
space $\on{Op}_{^L G}^{\la}$ of $^L B[[t]]$-equivalence classes of
operators of the form
\begin{equation}    \label{nilp oper mu3}
\nabla = \pa_t + \sum_{i=1}^\ell t^{\langle \check\al_i,\la \rangle}
\psi_i(t) f_i + {\mb v}(t),
\end{equation}
where $\psi_i(t) \in \C[[t]], \psi_i(0) \neq 0$, ${\mb v}(t) \in
{}^L \bb[[t]]$.

\begin{lem}    \label{fiber over trivial}
Suppose that the local system underlying an oper $\chi \in \on{Op}_{^L
G}(D^\times)$ is trivial. Then $\chi$ belongs to the disjoint union of
the subsets $\on{Op}_{^L G}^{\la} \subset \on{Op}_{^L G}(D^\times),
\la \in P^+$.
\end{lem}

\proof
It is clear from the definition that any oper in $\on{Op}_{^L
G}^{\la}$ is regular on the disc $D$ and is therefore $^L G\ppart$
gauge equivalent to the trivial connection.

Now suppose that we have an oper $\chi = (\CF,\nabla,\CF_{^L B})$ such
that the underlying $^L G$-local system is trivial. Then $\nabla$ is
$^L G\ppart$ gauge equivalent to a regular connection, that is one of
the form $\pa_t + A(t)$, where $A(t) \in {}^L \g[[t]]$. We have the
decomposition $^L G\ppart = {}^L G[[t]] {}^L B\ppart$. The gauge
action of $^L G[[t]]$ clearly preserves the space of regular
connections. Therefore if an oper connection $\nabla$ is $^L G\ppart$
gauge equivalent to a regular connection, then its $^L B\ppart$ gauge
class already must contain a regular connection. The oper condition
then implies that this gauge class contains a connection operator of
the form \eqref{nilp oper mu3} for some dominant integral weight $\la$
of $^L G$. Therefore $\chi \in
\on{Op}_{^L G}^{\la}$.
\qed

Thus, we see that the set of opers corresponding to the (unique)
unramified Langlands parameter is the disjoint union $\bigsqcup_{\la
\in P^+} \on{Op}_{^L G}^{\la}$. We call such opers ``unramified''.
The following result then confirms our expectation that the category
$\gmod_\chi$ is ``unramified'', that is contains non-zero
$G[[t]]$-equivariant objects, if and only if $\chi$ is unramified (see
\cite{FG:fusion} for a proof).

\begin{lem}    \label{inv vect}
The category $\gmod_\chi$ contains a non-zero $(\ghat_{\ka_c},G[[t]])$
Harish-Chandra module if and only if
\begin{equation}    \label{unram subscheme}
\chi \in \bigsqcup_{\la \in P^+} \on{Op}_{^L G}^{\la}.
\end{equation}
\end{lem}

The next question is to describe the category $\gmod^{G[[t]]}_\chi$ of
$(\ghat_{\ka_c},G[[t]])$ modules for $\chi \in \on{Op}_{^L G}^{\la}$.

\subsection{Categories of $G[[t]]$-equivariant modules}    \label{cat
unram}

Let us recall from \secref{unr rep} that the space of $K_0$-invariant
vectors in an unramified irreducible representation of $G(F)$ is
always one-dimensional. We have proposed that the category
$\gmod^{G[[t]]}_\chi$ should be viewed as a categorical analogue of
this space. Therefore we expect it to be the simplest possible abelian
category: the category of $\C$-vector spaces. Here we assume that
$\chi$ belongs to the union of the spaces $\on{Op}_{^L G}^{\la}$,
where $\la \in P^+$, for otherwise the category $\gmod^{G[[t]]}_\chi$
would be trivial (zero object is the only object).

In this subsection we will prove, following \cite{FG:exact} (see also
\cite{BD}), that our expectation is in fact correct provided that $\la
= 0$, in which case $\on{Op}_{^L G}^{0} = \on{Op}_{^L G}(D)$, and so
$$
\chi \in \on{Op}_{^L G}(D) \subset \on{Op}_{^L G}(D^\times).
$$
We will also conjecture that this is true for $\chi \in \on{Op}_{^L
G}^{\la}$ for all $\la \in P^+$.

Recall the vacuum module
$$
\BV_0 = \on{Ind}_{\g[[t]] \oplus \C {\mb 1}}^{\ghat_{\ka_c}} \C.
$$
According to \cite{FF:gd,F:wak}, we have
\begin{equation}    \label{End V and reg opers}
\on{End}_{\ghat_{\ka_c}} \BV_0 \simeq \on{Fun} \on{Op}_{^L G}(D).
\end{equation}
Let $\chi \in \on{Op}_{^L G}(D) \subset \on{Op}_{^L
G}(D^\times)$. Then $\chi$ defines a character of the algebra
$\on{End}_{\ghat_{\ka_c}} \BV_0$. Let $\BV_0(\chi)$ be the quotient of
$\BV_0$ by the kernel of this character. Then we have the following
result.

\begin{thm}    \label{single nabla}
Let $\chi \in \on{Op}_{^L G}(D) \subset \on{Op}_{^L
G}(D^\times)$. Then the category $\gmod^{G[[t]]}_\chi$ is equivalent
to the category of vector spaces: its unique, up to isomorphism,
irreducible object is $\BV_0(\chi)$ and any other object is isomorphic
to the direct sum of copies of $\BV_0(\chi)$.
\end{thm}

This theorem provides the first piece of evidence for \conjref{equiv
loc syst}: we see that the categories $\gmod^{G[[t]]}_\chi$ are
equivalent to each other for all $\chi \in \on{Op}_{^L G}(D)$.

It is more convenient to consider, instead of an individual regular
$^L G$-oper $\chi$, the entire family $\on{Op}_{^L G}^{0} =
\on{Op}_{^L G}(D)$ of regular opers on the disc
$D$. Let $\gmod_{\on{reg}}$ be the full subcategory of the category
$\gmod$ whose objects are $\ghat_{\ka_c}$-modules on which the action
of the center $Z(\ghat)$ factors through the homomorphism
$$
Z(\ghat) \simeq \on{Fun} \on{Op}_{^L G}(D^\times) \to \on{Fun}
\on{Op}_{^L G}(D).
$$
Note that the category $\gmod_{\on{reg}}$ is an example of a category
$\gmod_V$ introduced in \secref{opers and loc syst}, in the case when
$V = \on{Op}_{^L G}(D)$.

Let $\gmod^{G[[t]]}_{\on{reg}}$ be the corresponding
$G[[t]]$-equivariant category. It is instructive to think of
$\gmod_{\on{reg}}$ and $\gmod^{G[[t]]}_{\on{reg}}$ as categories
fibered over $\on{Op}_{^L G}(D)$, with the fibers over $\chi \in
\on{Op}_{^L G}(D)$ being $\gmod_\chi$ and $\gmod^{G[[t]]}_\chi$,
respectively.

We will now describe the category $\gmod^{G[[t]]}_{\on{reg}}$. This
description will in particular imply \thmref{single nabla}.

In order to simplify our formulas, in what follows we will use the
following notation for $\on{Fun} \on{Op}_{^L G}(D)$:
$$
\zz = \zz(\ghat) = \on{Fun} \on{Op}_{^L G}(D).
$$

Let $\opmod$ be the category of modules over the commutative algebra
$\zz$. Equivalently, this is the category of quasicoherent sheaves on
the space $\on{Op}_{^L G}(D)$.

By definition, any object of $\gmod^{G[[t]]}_{\reg}$ is a
$\zz$-module. Introduce the functors
\begin{align*}
\sF: &\gmod^{G[[t]]}_{\reg} \to \opmod, \qquad M \mapsto
\Hom_{\ghat_{\ka,c}}(\BV_0,M), \\ \sG: &\opmod \to
\gmod^{G[[t]]}_{\reg}, \qquad {\mc F}
\mapsto \BV_0 \underset{\zz}\otimes {\mc F}.
\end{align*}

The following theorem has been proved in \cite{FG:exact}, Theorem 6.3
(important results in this direction were obtained earlier in
\cite{BD}).

\begin{thm} \label{equivalence of categories}
The functors $\sF$ and $\sG$ are mutually inverse equivalences of
categories
\begin{equation}    \label{equiv reg}
\gmod^{G[[t]]}_{\reg} \simeq \opmod.
\end{equation}
\end{thm}

This immediately implies \thmref{single nabla}. Indeed, for
each $\chi \in \on{Op}_{^L G}(D)$ the category $\gmod^{G[[t]]}_\chi$
is the full subcategory of $\gmod^{G[[t]]}_{\reg}$ which are
annihilated, as $\zz$-modules, by the maximal ideal $I_\chi$ of
$\chi$. By \thmref{equivalence of categories}, this category is
equivalent to the category of $\zz$-modules annihilated by
$I_\chi$. But this is the category of $\zz$-modules supported
(scheme-theoretically) at the point $\chi$, which is equivalent to the
category of vector spaces.

\subsection{The action of the spherical Hecke algebra}
\label{action of sph hecke}

In \secref{unr rep} we discussed irreducible unramified
representations of the group $G(F)$, where $F$ is a local
non-archimedian field. We have seen that such representations are
parameterized by conjugacy classes of the Langlands dual group $^L
G$. Given such a conjugacy class $\ga$, we have an irreducible
unramified representation $(R_\ga,\pi_\ga)$, which contains a
one-dimensional subspace $(R_\ga)^{\pi_\ga(K_0)}$ of $K_0$-invariant
vectors. The spherical Hecke algebra $H(G(F),K_0)$, which is
isomorphic to $\on{Rep} {}^L G$ via the Satake isomorphism, acts on
this space by a character $\phi_\ga$, see formula \eqref{Hecke
eigenfunction}.

In the geometric setting, we have argued that for any $\chi \in
\on{Op}_{^L G}(D)$ the category $\gmod_\chi$, equipped with an
action of the loop group $G\ppart$, should be viewed as a
categorification of $(R_\ga,\pi_\ga)$. Furthermore, its subcategory
$\gmod^{G[[t]]}_\chi$ of $(\ghat_{\ka_c},G[[t]])$ Harish-Chandra
modules should be viewed as a categorification of the one-dimensional
space $(R_\ga)^{\pi_\ga(K_0)}$. According to \thmref{single nabla},
the latter category is equivalent to the category of vector spaces,
which is consistent with our expectations.

We now discuss the categorical analogue of the action of the spherical
Hecke algebra.

As explained in \secref{cat hecke}, the categorical analogue of the
spherical Hecke algebra is the category of $G[[t]]$-equivariant ${\mc
D}$-modules on the {\bf affine Grassmannian} $\on{Gr} =
G\ppart/G[[t]]$. We refer the reader to \cite{BD,FG:local} for the
precise definition of $\on{Gr}$ and this category. There is an
important property that is satisfied in the unramified case: the
convolution functors with these ${\mc D}$-modules are exact, which
means that we do not need to consider the derived category; the
abelian category of such ${\mc D}$-modules will do. Let us denote this
abelian category by ${\mc H}(G\ppart,G[[t]])$.

According to the results of \cite{MV}, this category carries a natural
structure of tensor category, which is equivalent to the tensor
category ${\mc Rep} \, {}^L G$ of representations of $^L G$. This should
be viewed as a categorical analogue of the Satake isomorphism. Thus,
for each object $V$ of ${\mc Rep} \, {}^L G$ we have an object of ${\mc
H}(G\ppart,G[[t]])$ which we denote by $\CH_V$. What should be the
analogue of the Hecke eigenvector property \eqref{Hecke
eigenfunction}?

As we explained in \secref{cat hecke}, the category ${\mc
H}(G\ppart,G[[t]])$ naturally acts on the category
$\gmod^{G[[t]]}_\chi$, and this action should be viewed as a
categorical analogue of the action of $H(G(F),K_0)$ on
$(R_\ga)^{\pi_\ga(K_0)}$.

Now, by \thmref{single nabla}, any object of $\gmod^{G[[t]]}_\chi$
is a direct sum of copies of $\BV_0(\chi)$. Therefore it is sufficient
to describe the action of ${\mc H}(G\ppart,G[[t]])$ on
$\BV_0(\chi)$. This action is described by the following statement,
which follows from \cite{BD}: there exists a family of isomorphisms
\begin{equation}    \label{cat action}
\al_V: \CH_V \star \BV_0(\chi) \overset\sim\longrightarrow \underline{V}
\otimes \BV_0(\chi), \qquad V \in {\mc Rep} \, {}^L G,
\end{equation}
where $\underline{V}$ is the vector space underlying the
representation $V$. Moreover, these isomorphisms are compatible with
the tensor product structure on $\CH_V$ (given by the convolution) and
on $\underline{V}$ (given by tensor product of vector spaces).

In view of \thmref{single nabla}, this is not surprising. Indeed, it
follows from the definition that $\CH_V \star \BV_0(\chi)$ is again an
object of the category $\gmod^{G[[t]]}_\chi$. Therefore it must be
isomorphic to $U_V \otimes_\C \BV_0(\chi)$, where $U_V$ is a vector
space. But then we obtain a functor ${\mc H}(G\ppart,G[[t]]) \to
{\mc Vect}, {\mc H}_V \mapsto U_V$. It follows from the construction
that this is a tensor functor. Therefore the standard Tannakian
formalism implies that $U_V$ is isomorphic to $\underline{V}$.

The isomorphisms \eqref{cat action} should be viewed as the
categorical analogues of the Hecke eigenvector conditions \eqref{Hecke
eigenfunction}. The difference is that while in \eqref{Hecke
eigenfunction} the action of elements of the Hecke algebra on a
$K_0$-invariant vector in $R_\ga$ amounts to multiplication by a
scalar, the action of an object of the Hecke category ${\mc
H}(G\ppart,G[[t]])$ on the $G[[t]]$-equivariant object $\BV_0(\chi)$ of
$\gmod_\chi$ amounts to multiplication by a {\em vector space},
namely, the vector space underlying the corresponding representation
of $^L G$. It is natural to call a module satisfying this property a
{\bf Hecke eigenmodule}. Thus, we obtain that $\BV_0(\chi)$ is a Hecke
eigenmodule. This is in agreement with our expectation that the
category $\gmod^{G[[t]]}_\chi$ is a categorical version of the space 
of $K_0$-invariant vectors in $R_\ga$.

One ingredient that is missing in the geometric case is the conjugacy
class $\ga$ of $^L G$. We recall that in the classical Langlands
correspondence this was the image of the Frobenius element of the
Galois group $\on{Gal}(\overline{\mathbb F}_q/\Fq)$, which does not
have an analogue in the geometric setting where our ground field is
$\C$, which is algebraically closed. So while unramified local systems
in the classical case are parameterized by the conjugacy classes
$\ga$, there is only one, up to an isomorphism, unramified local
system in the geometric case. However, this local system has a large
group of automorphisms, namely, $^L G$ itself. One can argue that what
replaces $\ga$ in the geometric setting is the action of this group
$^L G$ by automorphisms of the category $\gmod_\chi$, which we will
discuss in the next two sections.

\subsection{Categories of representations and ${\mc D}$-modules}
\label{reps and Dmod}

When we discussed the procedure of categorification of representations
in \secref{back to groups}, we saw that there are two possible
scenarios for constructing categories equipped with an action of the
loop group $G\ppart$. In the first one we consider categories of ${\mc
D}$-modules on the ind-schemes $G\ppart/K$, where $K$
is a ``compact'' subgroup of $G\ppart$, such as $G[[t]]$ or the
Iwahori subgroup. In the second one we consider categories of
representations $\gmod_\chi$. So far we have focused exclusively on
the second scenario, but it is instructive to also discuss categories
of the first type.

In the toy model considered in \secref{finite groups} we discussed the
category of $\g$-modules with fixed central character and the category
of ${\mc D}$-modules on the flag variety $G/B$. We have argued that
both could be viewed as categorifications of the representation of the
group $G(\Fq)$ on the space of functions on $(G/B)(\Fq)$. These
categories are equivalent, according to the Beilinson-Bernstein
theory, with the functor of global sections connecting the two. Could
something like this be true in the case of affine Kac-Moody algebras
as well?

The affine Grassmannian $\Gr = G\ppart/G[[t]]$ may be viewed as the
simplest possible analogue of the flag variety $G/B$ for the loop
group $G\ppart$. Consider the category of ${\mc D}$-modules on
$G\ppart/G[[t]]$ (see \cite{BD,FG:local} for the precise
definition). We have a functor of global sections from this category
to the category of $\g\ppart$-modules. In order to obtain
$\ghat_{\ka_c}$-modules, we need to take instead the category ${\mc
D}_{\ka_c}\mod$ of ${\mc D}$-modules twisted by a line bundle ${\mc
L}_{\ka_c}$. This is the unique line bundle ${\mc L}_{\ka_c}$ on $\Gr$
which carries an action of $\ghat_{\ka_c}$ (such that the central
element ${\mb 1}$ is mapped to the identity) lifting the natural
action of $\g\ppart$ on $\Gr$. Then for any object $\CM$ of ${\mc
D}_{\ka_c}\mod$, the space of global sections $\Gamma(\Gr,\CM)$ is a
$\ghat_{\ka_c}$-module. Moreover, it is known (see \cite{BD,FG:exact})
that $\Gamma(\Gr,\CM)$ is in fact an object of
$\gmod_{\reg}$. Therefore we have a functor of global sections
$$
\Gamma: {\mc D}_{\ka_c}\mod \to \gmod_{\reg}.
$$
We note that the categories ${\mc D}\mod$ and ${\mc D}_{\ka_c}\mod$
are equivalent under the functor $\CM \mapsto \CM \otimes {\mc
L}_{\ka_c}$. But the corresponding global sections functors are
very different.

However, unlike in the Beilinson-Bernstein scenario, the functor
$\Gamma$ cannot possibly be an equivalence of categories. There are
two reasons for this. First of all, the category $\gmod_{\reg}$ has a
large center, namely, the algebra $\zz = \on{Fun}
\on{Op}_{^L G}(D)$, while the center of the category ${\mc
D}_{\ka_c}\mod$ is trivial.\footnote{Recall that we are under the
assumption that $G$ is a connected simply-connected algebraic group,
and in this case $\Gr$ has one connected component. In general, the
center of the category ${\mc D}_{\ka_c}\mod$ has a basis enumerated by
the connected components of $\Gr$ and is isomorphic to the group
algebra of the finite group $\pi_1(G)$.} The second, and more serious,
reason is that the category ${\mc D}_{\ka_c}\mod$ carries an
additional symmetry, namely, an action of the tensor category ${\mc
Rep} {}^L G$ of representations of the Langlands dual group $^L G$,
and this action trivializes under the functor $\Gamma$ as we explain
presently.

Over $\on{Op}_{^L G}(D)$ there exists a canonical principal $^L
G$-bundle, which we will denote by $\CP$. By definition, the fiber of
$\CP$ at $\chi = (\F,\nabla,\F_{^L B}) \in \on{Op}_{^L G}(D)$ is
$\F_0$, the fiber at $0 \in D$ of the $^L G$-bundle $\F$ underlying
$\chi$. For an object $V \in {\mc Rep} \, {}^L G$ let us denote by $\CV$
the associated vector bundle over $\on{Op}_{^L G}(D)$, i.e.,
$$
\CV = \CP \underset{^L G}\times V.
$$

Next, consider the category ${\mc D}_{\ka_c}\mod^{G[[t]]}$ of
$G[[t]]$-equivariant ${\mc D}_{\ka_c}$-modules on $\Gr$. It is
equivalent to the category
$$
{\mc D}\mod^{G[[t]]} = {\mc H}(G\ppart,G[[t]])
$$
considered above. This is a tensor category, with respect to the
convolution functor, which is equivalent to the category ${\mc Rep} \,
{}^L G$. We will use the same notation ${\mc H}_V$ for the object of
${\mc D}_{\ka_c}\mod^{G[[t]]}$ corresponding to $V \in {\mc Rep} \, {}^L
G$.  The category ${\mc D}_{\ka_c}\mod^{G[[t]]}$ acts on ${\mc
D}_{\ka_c}\mod$ by convolution functors
$$
\CM \mapsto {\mc H}_V \star \CM
$$
which are exact. This amounts to a tensor action of the category ${\mc
Rep} {} ^L G$ on ${\mc D}_{\ka_c}\mod$.

Now, A. Beilinson and V. Drinfeld have proved in \cite{BD} that there
are functorial isomorphisms
$$
\Gamma(\Gr,{\mc H}_V \star \CM) \simeq \Gamma(\Gr,\CM)
\underset{\zz}\otimes \CV, \qquad V \in 
{\mc Rep} \, {}^L G,
$$
compatible with the tensor structure. Thus, we see that there are
non-isomorphic objects of ${\mc D}_{\ka_c}\mod$, which the functor
$\Gamma$ sends to isomorphic objects of
$\gmod_{\reg}$. Therefore the category ${\mc D}_{\ka_c}\mod$ and
the functor $\Gamma$ need to be modified in order to have a chance to
obtain a category equivalent to $\gmod_{\reg}$.

In \cite{FG:local} it was shown how to modify the category ${\mc
D}_{\ka_c}\mod$, by simultaneously "adding" to it $\zz$ as a center,
and "dividing" it by the above ${\mc Rep} \, {} ^L G$-action. As the
result, we obtain a candidate for a category that can be equivalent to
$\gmod_{\reg}$.  This is the category of {\bf Hecke eigenmodules} on
$\Gr$, denoted by ${\mc D}_{\ka_c}^{\Hecke}\mod_{\reg}$.

By definition, an object of ${\mc D}_{\ka_c}^{\Hecke}\mod_{\reg}$ is
an object of ${\mc D}_{\ka_c}\mod$, equipped with an action of the
algebra $\zz$ by endomorphisms and a system of
isomorphisms
$$
\al_V: {\mc H}_V \star \CM \overset{\sim}\longrightarrow
\CV \underset{\zz}\otimes \CM, \qquad V \in 
{\mc Rep} \, {}^L G,
$$
compatible with the tensor structure.

The above functor $\Gamma$ naturally gives rise to a functor 
\begin{equation} \label{glob Hecke sections}
\Gamma^{\Hecke}: {\mc D}_{\ka_c}^{\Hecke}\mod_{\reg} \to \gmod_{\reg}.
\end{equation}

This is in fact a general property. Suppose for simplicity that we
have an abelian category ${\mc C}$ which is acted upon by the tensor
category ${\mc Rep} \, H$, where $H$ is an algebraic group; we denote
this functor by $\CM \mapsto \CM \star V, V \in {\mc Rep} \, H$. Let
${\mc C}^\Hecke$ be the category whose objects are collections
$(\CM,\{\al_V\}_{V \in {\mc Rep} \, H})$, where $\CM \in {\mc C}$ and
$\{ \al_V \}$ is a compatible system of isomorphisms
$$
\al_V: \CM \star V \overset{\sim}\longrightarrow \underline{V}
\underset{\C}\otimes \CM, \qquad V \in {\mc Rep} \, H,
$$
where $\underline{V}$ is the vector space underlying $V$. One may
think of ${\mc C}^\Hecke$ as the ``de-equivarian\-tized'' category
${\mc C}$ with respect to the action of $H$. It carries a natural
action of the group $H$: for $h \in H$, we have
$$
h \cdot (\CM,\{\al_V\}_{V \in {\mc Rep} \, H}) = (\CM,\{ (h\otimes
\on{id}_\CM) \circ \al_V\}_{V \in {\mc Rep} \, H}).
$$
In other words, $\CM$ remains unchanged, but the isomorphisms $\al_V$
get composed with $h$.

The category ${\mc C}$ may in turn be reconstructed as the category of
$H$-equivariant objects of ${\mc C}^\Hecke$ with respect to this
action, see \cite{Ga}.

Suppose that we have a functor $\sG: {\mc C} \to {\mc C}'$, such that
we have functorial isomorphisms
\begin{equation}    \label{syst}
\sG(\CM \star V) \simeq \sG(\CM) \underset{\C}\otimes \underline{V},
\qquad V \in {\mc Rep} \, H,
\end{equation}
compatible with the tensor structure. Then, according to \cite{AG},
there exists a functor ${\sG}^\Hecke: \CC^\Hecke\to \CC'$ such that
$\sG \simeq \sG^\Hecke \circ \on{Ind}$, where the functor $\on{Ind}:
{\mc C} \to {\mc C}^\Hecke$ sends $\CM$ to $\CM \star {\mc O}_{H}$,
where ${\mc O}_{H}$ is the regular representation of $H$. The functor
$\sG^\Hecke$ may be explicitly described as follows: the isomorphisms
$\al_V$ and \eqref{syst} give rise to an action of the algebra ${\mc
O}_{H}$ on $\sG(\CM)$, and ${\sG}^\Hecke(\CM)$ is obtained by taking
the fiber of $\sG(\CM)$ at $1 \in H$.

We take ${\mc C}={\mc D}_{\ka_c}\mod$, ${\mc C}'=\gmod_{\reg}$, and
$\sG=\Gamma$. The only difference is that now we are working over the
base $\on{Op}_{^L G}(D)$, which we have to take into account. Thus, we
obtain a functor \eqref{glob Hecke sections} (see
\cite{FG:local,FG:equiv} for more details). Moreover, the left action of
the group $G\ppart$ on $\Gr$ gives rise to its action on the category
${\mc D}_{\ka_c}^{\Hecke}\mod_{\reg}$, and the functor
$\Gamma^{\Hecke}$ intertwines this action with the action of $G\ppart$
on $\gmod_{\reg}$.

The following was conjectured in \cite{FG:local}:

\begin{conj}    \label{equiv with Dmod}
The functor $\Gamma^{\Hecke}$ in formula \eqref{glob Hecke sections}
defines an equivalence of the categories ${\mc
D}_{\ka_c}^{\Hecke}\mod_{\reg}$ and $\gmod_{\reg}$.
\end{conj}

It was proved in \cite{FG:local} that the functor $\Gamma^{\Hecke}$,
when extended to the derived categories, is fully
faithful. Furthermore, it was proved in \cite{FG:equiv} that it sets
up an equivalence of the corresponding $I^0$-equivariant categories,
where $I^0 = [I,I]$ is the radical of the Iwahori subgroup.

Let us specialize \conjref{equiv with Dmod} to a point $\chi =
(\F,\nabla,\F_{^L B}) \in \on{Op}_{^L G}(D)$. Then on the right hand
side we consider the category $\gmod_\chi$, and on the left hand
side we consider the category ${\mc D}_{\ka_c}^{\Hecke}\mod_\chi$. Its
objects consist of a ${\mc D}_{\ka_c}$-module $\CM$ and a
collection of isomorphisms
\begin{equation}    \label{Hecke for Dmod}
\al_V: \CH_V \star \CM \overset\sim\longrightarrow
V_{\F_0} \otimes \CM, \qquad V \in {\mc Rep} \, {}^L G.
\end{equation}
Here $V_{\F_0}$ is the twist of the representation $V$ by
the $^L G$-torsor $\F_0$. These isomorphisms have to be compatible
with the tensor structure on the category ${\mc H}(G\ppart,G[[t]])$.

\conjref{equiv with Dmod} implies that there is a canonical
equivalence of categories
\begin{equation}    \label{Dmod fixed chi}
{\mc D}_{\ka_c}^{\Hecke}\mod_\chi \simeq \gmod_\chi.
\end{equation}
It is this conjectural equivalence that should be viewed as an
analogue of the Beilinson-Bernstein equivalence.

{}From this point of view, one can think of each of the categories
${\mc D}_{\ka_c}^{\Hecke}\mod_\chi$ as the second incarnation of the
sought-after Langlands category ${\mc C}_{\sigma_0}$ corresponding to
the trivial $^L G$-local system.

Now we give another explanation why it is natural to view the
category ${\mc D}_{\ka_c}^{\Hecke}\mod_\chi$ as a categorification of
an unramified representation of the group $G(F)$. First of all,
observe that these categories are all equivalent to each other and to
the category ${\mc D}_{\ka_c}^{\Hecke}\mod$, whose objects are ${\mc
D}_{\ka_c}$-modules $\CM$ together with a collection of isomorphisms
\begin{equation}    \label{Hecke for Dmod1}
\al_V: \CH_V \star \CM \overset\sim\longrightarrow
\underline{V} \otimes \CM, \qquad V \in {\mc Rep} \, {}^L G.
\end{equation}
Comparing formulas \eqref{Hecke for Dmod} and \eqref{Hecke for Dmod1},
we see that there is an equivalence
$$
{\mc D}_{\ka_c}^{\Hecke}\mod_\chi \simeq {\mc D}_{\ka_c}^{\Hecke}\mod,
$$
for each choice of trivialization of the $^L G$-torsor $\F_0$
(the fiber at $0 \in D$ of the principal $^L G$-bundle $\F$ on $D$
underlying the oper $\chi$).

Now recall from \secref{unr rep} that to each semi-simple conjugacy
class $\gamma$ in $^L G$ corresponds an irreducible unramified
representation $(R_\gamma,\pi_\gamma)$ of $G(F)$ via the Satake
correspondence \eqref{unramified langlands1}. It is known that there
is a non-degenerate pairing
$$
\langle,\rangle: R_\gamma \times R_{\gamma^{-1}} \to \C,
$$
in other words, $R_{\ga^{-1}}$ is the representation of $G(F)$ which
is contragredient to $R_\ga$ (it may be realized in the space of
smooth vectors in the dual space to $R_\ga$).

Let $v \in R_{\ga^{-1}}$ be a non-zero vector such that $K_0 v = v$
(this vector is unique up to a scalar). It then satisfies the Hecke
eigenvector property \eqref{Hecke eigenfunction} (in which we need to
replace $\ga$ by $\ga^{-1}$). This allows us to embed $R_\ga$ into the
space of locally constant smooth right $K_0$-invariant functions on
$G(F)$ (equivalently, on $G(F)/K_0$), by using matrix coefficients, as
follows:
$$
u \in R_\gamma \mapsto f_u, \qquad f_u(g) = \langle u,g v
\rangle.
$$
The Hecke eigenvector property \eqref{Hecke eigenfunction} implies
that the functions $f_u$ are right $K_0$-invariant and satisfy the
condition
\begin{equation}    \label{hecke for functions}
f \star H_V = \on{Tr}(\gamma^{-1},V) f,
\end{equation}
where $\star$ denotes the convolution product \eqref{conv}. Let
$C(G(F)/K_0)_\gamma$ be the space of locally constant smooth functions
on $G(F)/K_0$ satisfying \eqref{hecke for functions}. It carries a
representation of $G(F)$ induced by its left action on $G(F)/K_0$. We
have constructed an injective map $R_\gamma \to C(G(R)/G(R))_\gamma$,
and one can show that for generic $\gamma$ it is an isomorphism.

Thus, we obtain a realization of an irreducible unramified
representation of $G(F)$ in the space of functions on the quotient
$G(F)/K_0$ satisfying the Hecke eigenfunction condition \eqref{hecke
for functions}. The Hecke eigenmodule condition \eqref{Hecke for
Dmod1} may be viewed as a categorical analogue of \eqref{hecke for
functions}. Therefore the category ${\mc D}_{\ka_c}^{\Hecke}\mod$ of
twisted ${\mc D}$-modules on $\Gr = G\ppart/K_0$ satisfying the Hecke
eigenmodule condition \eqref{Hecke for Dmod1}, equipped with a
$G\ppart$-action appears to be a natural categorification of the
irreducible unramified representations of $G(F)$.

\subsection{Equivalences between categories of modules}    \label{eq
of cat}

All opers in $\on{Op}_{^L G}(D)$ correspond to one and the same $^L
G$-local system, namely, the trivial local system. Therefore,
according to \conjref{equiv loc syst}, we expect that the categories
$\gmod_\chi$ are equivalent to each other. More precisely, for each
isomorphism between the underlying local systems of any two opers in
$\on{Op}_{^L G}(D)$ we wish to have an equivalence of the
corresponding categories, and these equivalences should be compatible
with respect to the operation of composition of these isomorphisms.

Let us spell this out in detail. Let $\chi=(\CF,\nabla,\CF_{^L B})$
and $\chi'=(\CF',\nabla',\CF'_{^L B})$ be two opers in $\on{Op}_{^L
G}(D)$. Then an isomorphism between the underlying local systems
$(\CF,\nabla)
\overset\sim\longrightarrow (\CF',\nabla')$ is the same as an
isomorphism $\CF_0 \overset\sim\longrightarrow \CF'_0$ between the $^L
G$-torsors $\CF_0$ and $\CF'_0$, which are the fibers of the $^L
G$-bundles $\CF$ and $\CF'$, respectively, at $0 \in D$. Let us denote
this set of isomorphisms by $\on{Isom}_{\chi,\chi'}$. Then we have
$$
\on{Isom}_{\chi,\chi'} = \CF_0 \underset{^L G}\times {}^L G
\underset{^L G}\times \CF'_0,
$$
where we twist $^L G$ by $\CF_0$ with respect to the left action and
by $\CF'_0$ with respect to the right action. In particular,
$$
\on{Isom}_{\chi,\chi} = {}^L G_{\CF_0} = \CF_0 \underset{^L G}\times
\on{Ad} {}^L G
$$
is just the group of automorphisms of $\CF_0$.

It is instructive to combine the sets $\on{Isom}_{\chi,\chi'}$
into a groupoid $\on{Isom}$ over $\on{Op}_{^L G}(D)$. Thus, by
definition $\on{Isom}$ consists of triples $(\chi,\chi',\phi)$,
where $\chi,\chi' \in \on{Op}_{^L G}(D)$ and $\phi \in
\on{Isom}_{\chi,\chi}$ is an isomorphism of the underlying local
systems. The two morphisms $\on{Isom} \to \on{Op}_{^L G}(D)$
correspond to sending such a triple to $\chi$ and $\chi'$. The
identity morphism $\on{Op}_{^L G}(D) \to \on{Isom}$ sends $\chi$ to
$(\chi,\chi,\on{Id})$, and the composition morphism
$$
\on{Isom} \underset{\on{Op}_{^L G}(D)}\times \on{Isom} \to \on{Isom}
$$
corresponds to composing two isomorphisms.

\conjref{equiv loc syst} has the following more precise formulation
for regular opers:

\begin{conj}    \label{equiv reg opers}
For each $\phi \in \on{Isom}_{\chi,\chi'}$ there exists an
equivalence
$$
E_\phi: \gmod_\chi \to \gmod_{\chi'},
$$
which intertwines the actions of $G\ppart$ on the two categories, such
that $E_{\on{Id}} = \on{Id}$ and there exist isomorphisms
$\beta_{\phi,\phi'}: E_{\phi \circ \phi'} \simeq E_\phi \circ
E_{\phi'}$ satisfying $$\beta_{\phi \circ \phi',\phi''}
\beta_{\phi,\phi'} = \beta_{\phi,\phi' \circ \phi''}
\beta_{\phi',\phi''}$$ for all isomorphisms $\phi,\phi',\phi''$,
whenever they may be composed in the appropriate order.

In other words, the groupoid $\Isom$ over $\on{Op}_{^L G}(D)$ acts on
the category $\gmod_{\reg}$ fibered over $\on{Op}_{^L G}(D)$,
preserving the action of $G\ppart$ along the fibers.
\end{conj}

In particular, this conjecture implies that the group $^L G_{\CF_0}$
acts on the category $\gmod_\chi$ for any $\chi \in \on{Op}_{^L
G}(D)$.

Now we observe that \conjref{equiv with Dmod} implies \conjref{equiv
reg opers}. Indeed, by \conjref{equiv with Dmod}, there is a
canonical equivalence of categories \eqref{Dmod fixed chi},
$$
{\mc D}_{\ka_c}^{\Hecke}\mod_\chi \simeq \gmod_\chi.
$$
It follows immediately from the definition of the category ${\mc
D}_{\ka_c}^{\Hecke}\mod_\chi$ (namely, formula \eqref{Hecke for Dmod})
that for each isomorphism $\phi \in
\on{Isom}_{\chi,\chi'}$, i.e., an isomorphism of the $^L G$-torsors
$\F_0$ and $\F'_0$ underlying the opers $\chi$ and $\chi'$, there is a
canonical equivalence
$$
{\mc D}_{\ka_c}^{\Hecke}\mod_\chi \simeq {\mc
D}_{\ka_c}^{\Hecke}\mod_{\chi'}.
$$
Therefore we obtain the sought-after equivalence $E_\phi: \gmod_\chi
\to \gmod_{\chi'}$. Furthermore, it is clear that these equivalences
satisfy the conditions of \conjref{equiv reg opers}. In particular,
they intertwine the actions of $G\ppart$, which affects the ${\mc
D}$-module $\CM$ underlying an object of ${\mc
D}_{\ka_c}^{\Hecke}\mod_\chi$, but does not affect the isomorphisms
$\al_V$.

Equivalently, we can express this by saying that the groupoid $\Isom$
naturally acts on the category ${\mc
D}_{\ka_c}^{\Hecke}\mod_{\reg}$. By \conjref{equiv with Dmod}, this
gives rise to an action of $\Isom$ on $\gmod_{\reg}$.

In particular, we construct an action of the group $(^L G)_{\F_0}$,
the twist of $^L G$ by the $^L G$-torsor $\F_0$ underlying a
particular oper $\chi$, on the category ${\mc
D}_{\ka_c}^{\Hecke}\mod_{\chi}$. Indeed, each element $g \in (^L
G)_{\F_0}$ acts on the $\F_0$-twist $V_{\F_0}$ of any
finite-dimensional representation $V$ of $^L G$. Given an object
$(\CM,(\al_V))$ of ${\mc D}_{\ka_c}^{\Hecke}\mod_{\chi'}$, we
construct a new object, namely, $(\CM,((g \otimes \on{Id}_\CM) \circ
\al_V))$. Thus, we do not change the ${\mc D}$-module $\CM$, but we
change the isomorphisms $\al_V$ appearing in the Hecke eigenmodule
condition
\eqref{Hecke for Dmod} by composing them with the action of $g$ on
$V_{\F_0}$. According to \conjref{equiv with Dmod}, the category ${\mc
D}_{\ka_c}^{\Hecke}\mod_\chi$ is equivalent to $\gmod_\chi$. Therefore
this gives rise to an action of the group $(^L G)_{\F_0}$ on
$\gmod_\chi$. But this action is much more difficult to describe in
terms of $\ghat_{\ka_c}$-modules.

\subsection{Generalization to other dominant integral weights}

We have extensively studied above the categories $\gmod_\chi$ and
$\gmod_\chi^{G[[t]]}$ associated to regular opers $\chi \in
\on{Op}_{^L G}(D)$. However, according to \lemref{fiber over trivial},
the (set-theoretic) fiber of the map $\al:
\on{Op}_{^L G}(D^\times) \to \Loc_{^L G}(D^\times)$ over the trivial
local system $\sigma_0$ is the disjoint union of the subsets
$\on{Op}_{^L G}^{\la}, \la \in P^+$. Here we discuss briefly
the categories $\gmod_\chi$ and $\gmod_\chi^{G[[t]]}$ for $\chi \in
\on{Op}_{^L G}^\la$, where $\la \neq 0$.

Consider the Weyl module $\BV_\la$ with highest weight $\la$,
$$
\BV_{\la} = U(\ghat_{\ka_c}) \underset{U(\g[[t]] \oplus \C{\mb
1})}\otimes V_\la.
$$
According to \cite{FG:weyl}, we have
\begin{equation}    \label{End V and reg opers1}
\on{End}_{\ghat_{\ka_c}} \BV_\la \simeq \on{Fun} \on{Op}_{^L G}^\la.
\end{equation}
Let $\chi \in \on{Op}_{^L G}^\la \subset \on{Op}_{^L
G}(D^\times)$. Then $\chi$ defines a character of the algebra
$\on{End}_{\ghat_{\ka_c}} \BV_\la$. Let $\BV_\la(\chi)$ be the
quotient of $\BV_\la$ by the kernel of this character. The following
conjecture of \cite{FG:weyl} is an analogue of \thmref{single nabla}:

\begin{conj}    \label{single nabla1}
Let $\chi \in \on{Op}_{^L G}^\la \subset \on{Op}_{^L
G}(D^\times)$. Then the category $\gmod^{G[[t]]}_\chi$ is equivalent
to the category of vector spaces: its unique, up to isomorphism,
irreducible object is $\BV_\la(\chi)$ and any other object is
isomorphic to the direct sum of copies of $\BV_\la(\chi)$.
\end{conj}

Note that this is consistent with \conjref{equiv loc syst}, which
tells us that the categories $\gmod^{G[[t]]}_\chi$ should be
equivalent to each other for all opers which are gauge equivalent to
the trivial local system on $D$.

\section{Local Langlands correspondence: tamely ramified case}
\label{tr case}

In the previous section we have considered categorical analogues of
the irreducible unramified representations of a reductive group $G(F)$
over a local non-archi\-me\-dian field $F$. We recall that these are
the representations containing non-zero vectors fixed by the maximal
compact subgroup $K_0 \subset G(F)$. The corresponding Langlands
parameters are unramified admissible homomorphisms from the
Weil-Deligne group $W'_F$ to $^L G$, i.e., those which factor through
the quotient
$$
W'_F \to W_F \to \Z,
$$
and whose image in $^L G$ is semi-simple. Such homomorphisms are
parameterized by semi-simple conjugacy classes in $^L G$.

We have seen that the categorical analogues of unramified
representations of $G(F)$ are the categories $\gmod_\chi$ (equipped
with an action of the loop group $G\ppart$), where $\chi$ is a $^L
G$-oper on $D^\times$ whose underlying $^L G$-local system is
trivial. These categories can be called unramified in the sense that
they contain non-zero $G[[t]]$-equivariant objects. The corresponding
Langlands parameter is the trivial $^L G$-local system $\sigma_0$ on
$D^\times$, which should be viewed as an analogue of an unramified
homomorphism $W'_F \to {} ^L G$. However, the local system $\sigma_0$
is realized by many different opers, and this introduces an additional
complication into our picture: at the end of the day we need to show
that the categories $\gmod_\chi$, where $\chi$ is of the above type,
are equivalent to each other. In particular, \conjref{equiv reg opers},
which describes what we expect to happen when $\chi \in \on{Op}_{^L
G}(D)$.

The next natural step is to consider categorical analogues of
representations of $G(F)$ that contain vectors invariant under the
Iwahori subgroup $I
\subset G[[t]]$, the preimage of a fixed Borel subgroup $B \subset G$
under the evaluation homomorphism $G[[t]] \to G$. We begin this
section by recalling a classification of these representations, due
to D. Kazhdan and G. Lusztig \cite{KL} and V. Ginzburg
\cite{Ginzburg}. We then discuss the categorical analogues of these
representations following \cite{FG:local}--\cite{FG:wak} and the
intricate interplay between the classical and the geometric
pictures.

\subsection{Tamely ramified representations}    \label{tame reps}

The Langlands parameters corresponding to irreducible representations
of $G(F)$ with $I$-invariant vectors are {\bf tamely ramified}
homomorphisms $W'_F \to {} ^L G$. Recall from \secref{langlands param}
that $W'_F = W_F \ltimes \C$. A homomorphism $W'_F \to {} ^L G$ is
called tamely ramified if it factors through the quotient
$$
W'_F \to \Z \ltimes \C.
$$
According to the relation \eqref{rel in weil}, the group $\Z \ltimes
\C$ is generated by two elements $F = 1 \in \Z$ (Frobenius) and
$M = 1 \in \C$ (monodromy) satisfying the relation 
\begin{equation}    \label{rel in weil1}
F M F^{-1} = q M.
\end{equation}

Under an admissible tamely ramified homomorphism the generator $F$
goes to a semi-simple element $\ga \in {} ^L G$ and the generator $M$
goes to a unipotent element $N \in {}^L G$. According to formula
\eqref{rel in weil1}, they have to satisfy the relation
\begin{equation}
\gamma N \gamma^{-1} = N^q.
\end{equation}
Alternatively, we may write $N = \exp(u)$, where $u$ is a
nilpotent element of $^L \g$. Then this relation becomes
$$
\gamma u \gamma^{-1} = qu.
$$

Thus, we have the following bijection between the sets of equivalence
classes
\begin{equation}    \label{tamely ramified langlands}
\boxed{\begin{matrix} \text{tamely ramified admissible} \\
\text{homomorphisms } W'_F \to {}^L G \end{matrix}} \quad
\Longleftrightarrow  \quad \boxed{\begin{matrix} \text{pairs }
\ga \in {}^L G, \text{ semi-simple}, \\
u \in {}^L \g, \text{ nilpotent}, \gamma u \gamma^{-1} = qu
\end{matrix}}
\end{equation}
In both cases equivalence relation amounts to conjugation by an
element of $^L G$.

Now to each Langlands parameter of this type we wish to attach an
irreducible representation of $G(F)$ which contains non-zero
$I$-invariant vectors. It turns out that if $G=GL_n$ there is indeed a
bijection, proved in \cite{BZ}, between the sets of equivalence
classes of the following objects:
\begin{equation}    \label{tamely ramified langlands1}
\boxed{\begin{matrix} \text{tamely ramified admissible} \\
\text{homomorphisms } W'_F \to GL_n \end{matrix}} \quad
\Longleftrightarrow  \quad \boxed{\begin{matrix} \text{irreducible
representations} \\ (R,\pi) \text{ of } GL_n(F), R^{\pi(I)} \neq 0 
\end{matrix}}
\end{equation}

However, such a bijection is no longer true for other reductive groups:
two new phenomena appear, which we discuss presently.

The first one is the appearance of $L$-{\bf packets}. One no longer
expects to be able to assign to a particular admissible homomorphism
$W'_F \to {} ^L G$ a single irreducible smooth representations of
$G(F)$. Instead, a finite collection of such representations (more
precisely, a collection of equivalence classes of representations) is
assigned, called an $L$-packet. In order to distinguish
representations in a given $L$-packet, one needs to introduce an
additional parameter. We will see how this is done in the case at
hand shortly. However, and this is the second subtlety alluded to
above, it turns out that not all irreducible representations of $G(F)$
within the $L$-packet associated to a given tamely ramified
homomorphism $W'_F
\to {} ^L G$ contain non-zero $I$-invariant vectors. Fortunately,
there is a certain property of the extra parameter used to distinguish
representations inside the $L$-packet that tells us whether the
corresponding representation of $G(F)$ has $I$-invariant vectors.

In the case of tamely ramified homomorphisms $W'_F \to {}^L G$ this
extra parameter is an irreducible representation $\rho$ of the finite
group $C(\ga,u)$ of components of the simultaneous centralizer of
$\ga$ and $u$ in $^L G$, on which the center of $^L G$ acts trivially
(see \cite{Lus}). In the case of $G=GL_n$ these centralizers are
always connected, and so this parameter never appears. But for other
reductive groups $G$ this group of components is often non-trivial. The
simplest example is when $^L G = G_2$ and $u$ is a
subprincipal nilpotent element of the Lie algebra $^L
\g$.\footnote{The term ``subprincipal'' means that the adjoint orbit
of this element has codimension $2$ in the nilpotent cone.} In this
case for some $\ga$ satisfying $\ga u \ga^{-1} = qu$ the group of
components $C(\ga,u)$ is the symmetric group $S_3$, which has three
irreducible representations (up to equivalence). Each of them
corresponds to a particular member of the $L$-packet associated with
the tamely ramified homomorphism $W'_F
\to {} ^L G$ defined by $(\ga,u)$. Thus, the $L$-packet consists
of three (equivalence classes of) irreducible smooth representations
of $G(F)$. However, not all of them contain non-zero $I$-invariant
vectors.

The representations $\rho$ of the finite group $C(\ga,u)$ which
correspond to representations of $G(F)$ with $I$-invariant vectors are
distinguished by the following property. Consider the {\bf Springer
fiber} $\on{Sp}_u$. We recall that
\begin{equation}    \label{Springer1}
\on{Sp}_{u} = \{ \bb' \in {}^L G/{}^L B \, | \, u \in \bb' \}.
\end{equation}
The group $C(\ga,u)$ acts on the homology of the variety
$\on{Sp}^\ga_u$ of $\ga$-fixed points of $\on{Sp}_u$. A representation
$\rho$ of $C(\ga,u)$ corresponds to a representation of $G(F)$ with
non-zero $I$-invariant vectors if and only if $\rho$ occurs in
the homology of $\on{Sp}_u$, $H_\bullet(\on{Sp}^\ga_u)$.

In the case of $G_2$ the Springer fiber $\on{Sp}_{u}$ of the
subprincipal element $u$ is a union of four projective lines connected
with each other as in the Dynkin diagram of $D_4$. For some $\ga$ the
set $\on{Sp}^\ga_u$ is the union of a projective line (corresponding
to the central vertex in the Dynkin diagram of $D_4$) and three points
(each in one of the remaining three projective lines). The
corresponding group $C(\ga,u) = S_3$ on $\on{Sp}^\ga_{u}$ acts
trivially on the projective line and by permutation of the three
points. Therefore the trivial and the two-dimensional representations
of $S_3$ occur in $H_\bullet(\on{Sp}^\ga_u)$, but the sign
representation does not. The irreducible representations of $G(F)$
corresponding to the first two contain non-zero $I$-invariant vectors,
whereas the one corresponding to the sign representation of $S_3$ does
not.

The ultimate form of the local Langlands correspondence for
representations of $G(F)$ with $I$-invariant vectors is then as
follows (here we assume, as in \cite{KL,Ginzburg}), that the group $G$
is split and has connected center):
\begin{equation}    \label{tamely ramified langlands2}
\boxed{\begin{matrix} \text{triples } (\ga,u,\rho), \ga u \ga^{-1} =
qu, \\ \rho \in {\mc Rep} \; C(\ga,u) \text{ occurs in }
H_\bullet(\on{Sp}^\ga_u,\C) \end{matrix}} \quad
\Longleftrightarrow  \quad \boxed{\begin{matrix} \text{irreducible
representations} \\ (R,\pi) \text{ of } G(F), R^{\pi(I)} \neq 0 
\end{matrix}}
\end{equation}
Again, this should be understood as a bijection between two sets of
equivalence classes of the objects listed. This bijection is due to
\cite{KL} (see also \cite{Ginzburg}). It was conjectured by Deligne
and Langlands, with a subsequent modification (addition of $\rho$)
made by Lusztig.

How to set up this bijection? The idea is to replace irreducible
representations of $G(F)$ appearing on the right hand side of
\eqref{tamely ramified langlands2} with irreducible modules
over the corresponding Hecke algebra $H(G(F),I)$. Recall from
\secref{Kinv} that this is the algebra of compactly supported $I$
bi-invariant functions on $G(F)$, with respect to convolution. It
naturally acts on the space of $I$-invariant vectors of any smooth
representation of $G(F)$ (see formula
\eqref{conv action}). Thus, we obtain a functor from the category of
smooth representations of $G(F)$ to the category of $H(G(F),I)$.
According to a theorem of A. Borel \cite{Borel}, it induces a
bijection between the set of equivalence classes of irreducible
representations of $G(F)$ with non-zero $I$-invariant vectors and the
set of equivalence classes of irreducible $H(G(F),I)$-modules.

The algebra $H(G(F),I)$ is known as the {\bf affine Hecke algebra} and
has the standard description in terms of generators and
relations. However, for our purposes we need another description,
due to \cite{KL,Ginzburg}, which identifies it with the
equivariant $K$-theory of the {\bf Steinberg variety}
$$
\on{St} = \wt{\mc N} \underset{{\mc N}}\times \wt{\mc N},
$$
where ${\mc N} \subset {} ^L \g$ is the nilpotent cone and $\wt{\mc
N}$ is the {\bf Springer resolution}
$$
\wt{\mc N} = \{ x \in {\mc N}, \bb' \in {}^L G/{}^L B \; | \;
x \in \bb' \}.
$$
Thus, a point of $\on{St}$ is a triple consisting of a nilpotent
element of $^L \g$ and two Borel subalgebras containing it. The group
$^L G \times \C^\times$ naturally acts on $\on{St}$, with $^L G$
conjugating members of the triple and $\C^\times$ acting by
multiplication on the nilpotent elements,
\begin{equation}    \label{C star action}
a \cdot (x,\bb',\bb'') = (a^{-1} x,\bb',\bb'').
\end{equation}

According to a theorem of \cite{KL,Ginzburg}, there is an isomorphism
\begin{equation}    \label{KL}
H(G(F),I) \simeq K^{^L G \times \C^\times}(\on{St}).
\end{equation}
The right hand side is the $^L G \times \C^\times$-equivariant
$K$-theory of $\on{St}$. It is an algebra with respect to a natural
operation of convolution (see \cite{Ginzburg} for details). It is also
a free module over its center, isomorphic to
$$
K^{^L G \times \C^\times}(\on{pt}) = \on{Rep} {} ^L G \otimes
\C[{\mb q},{\mb q}^{-1}].
$$
Under the isomorphism \eqref{KL} the element ${\mb q}$ goes to the
standard parameter ${\mb q}$ of the affine Hecke algebra $H(G(F),I)$
(here we consider $H(G(F),I)$ as a $\C[{\mb q},{\mb q}^{-1}]$-module).

Now, the algebra $K^{^L G \times \C^\times}(\on{St})$, and hence the
algebra $H(G(F),I)$, has a natural family of modules which are
parameterized precisely by the conjugacy classes of pairs $(\ga,u)$ as
above. On these modules $H(G(F),I)$ acts via a central character
corresponding to a point in $\on{Spec} \on{Rep} {} ^L G
\underset{\C}\otimes \C[{\mb q},{\mb q}^{-1}]$, which is just a pair
$(\ga,q)$, where $\ga$ is a semi-simple conjugacy class in $^L G$ and
$q \in \C^\times$. In our situation $q$ is the cardinality of the
residue field of $F$ (hence a power of a prime), but in what follows
we will allow a larger range of possible values of $q$: all non-zero
complex numbers except for the roots of unity. Consider the quotient
of $H(G(F),I)$ by the central character defined by $(\ga,u)$. This is
just the algebra $K^{^L G \times \C^\times}(\on{St})$, specialized at
$(\ga,q)$.  We denote it by $K^{^L G \times
\C^\times}(\on{St})_{(\ga,q)}$.

Now for a nilpotent element $u \in {\mc N}$ consider the Springer
fiber $\on{Sp}_u$. The condition that $\ga u \ga^{-1} = q u$ means
that $u$, and hence $\on{Sp}_u$, is stabilized by the action of
$(\ga,q) \in {}^L G \times \C^\times$ (see formula \eqref{C star
action}). Let $A$ be the smallest algebraic subgroup of $^L G \times
\C^\times$ containing $(\ga,q)$. The algebra $K^{^L G \times
\C^\times}(\on{St})_{(\ga,q)}$ naturally acts on the equivariant
$K$-theory $K^A(\on{Sp}_u)$ specialized at $(\ga,q)$,
$$
K^A(\on{Sp}_u)_{(\ga,q)} = K^A(\on{Sp}_u) \underset{\on{Rep} A}
\otimes \C_{(\ga,q)}.
$$
It is known that $K^A(\on{Sp}_u)_{(\ga,q)}$ is isomorphic to the
homology $H_\bullet(\on{Sp}^\ga_u)$ of the $\ga$-fixed subset of
$\on{Sp}_u$ (see \cite{KL,Ginzburg}). Thus, we obtain that
$K^A(\on{Sp}_u)_{(\ga,q)}$ is a module over $H(G(F),I)$.

Unfortunately, these $H(G(F),I)$-modules are not irreducible in
general, and one needs to work harder to describe the irreducible
modules over $H(G(F),I)$. For $G=GL_n$ one can show that each of these
modules has a unique irreducible quotient, and this way one recovers
the bijection
\eqref{tamely ramified langlands1}. But for a general groups $G$ the
finite groups $C(\ga,u)$ come into play. Namely, the group $C(\ga,u)$
acts on $K^A(\on{Sp}_u)_{(\ga,q)}$, and this action commutes with the
action of $K^{^L G \times
\C^\times}(\on{St})_{(\ga,q)}$. Therefore we have a decomposition
$$
K^A(\on{Sp}_u)_{(\ga,q)} = \bigoplus_{\rho \in \on{Irrep} C(\ga,u)}
\rho \otimes K^A(\on{Sp}_u)_{(\ga,q,\rho)},
$$
of $K^A(\on{Sp}_u)_{(\ga,q)}$ as a representation of $C(\ga,u) \times
H(G(F),I)$. One shows (see
\cite{KL,Ginzburg} for details) that each $H(G(F),I)$-module
$K^A(\on{Sp}_u)_{(\ga,q,\rho)}$ has a unique irreducible quotient, and
this way one obtains a parameterization of irreducible modules by the
triples appearing in the left hand side of
\eqref{tamely ramified langlands2}. Therefore we obtain that the same
set is in bijection with the right hand side of \eqref{tamely ramified
langlands2}. This is how the tame local Langlands correspondence
\eqref{tamely ramified langlands2}, also known as the
Deligne--Langlands conjecture, is proved.

\subsection{Categories admitting $(\ghat_{\ka_c},I)$ Harish-Chandra
modules}    \label{admitting I}

We now wish to find categorical analogues of the above results in the
framework of the categorical Langlands correspondence for loop groups.

As we explained in \secref{eq mod}, in the categorical setting a
representation of $G(F)$ is replaced by a category $\gmod_\chi$
equipped with an action of $G\ppart$, and the space of $I$-invariant
vectors is replaced by the subcategory of $(\ghat_{\ka_c},I)$
Harish-Chandra modules in $\gmod_\chi$. Hence the analogue of the
question which representations of $G(F)$ admit non-zero $I$-invariant
vectors becomes the following question: for what $\chi$ does the
category $\gmod_\chi$ contain non-zero $(\ghat_{\ka_c},I)$
Harish-Chandra modules?

To answer this question, we introduce the space $\on{Op}_{^L
G}^{\on{RS}}(D)$ of {\bf opers with regular singularity}. By
definition (see \cite{BD}, Sect. 3.8.8), an element of this space is
an $^L N[[t]]$-conjugacy class of operators of the form
\begin{equation}    \label{oper with RS}
\nabla = \pa_t + t^{-1} \left( p_{-1} + {\mb v}(t) \right),
\end{equation}
where ${\mb v}(t) \in {}^L \bb[[t]]$. One can show that a natural map
$\on{Op}_{^L G}^{\on{RS}}(D) \to \on{Op}_{^L G}(D^\times)$ is an
embedding.

Following \cite{BD}, we associate to an oper with regular singularity
its {\bf residue}. For an operator \eqref{oper with RS} the residue is
by definition equal to $p_{-1} + {\mb v}(0)$. Clearly, under gauge
transformations by an element $x(t)$ of $^L N[[t]]$ the residue gets
conjugated by $x(0) \in N$. Therefore its projection onto
$$
^L \g/{}^L G = \on{Spec} (\on{Fun} {}^L \g)^{^L G} = \on{Spec}
(\on{Fun} {}^L \h)^W = \h^*/W
$$
is well-defined.

Given $\mu \in \h^*$, we write $\varpi(\mu)$ for the projection of
$\mu$ onto $\h^*/W$. Finally, let $P$ be the set of integral (not
necessarily dominant) weights of $\g$, viewed as a subset of $\h^*$.
The next result follows from \cite{F:wak,FG:local}.

\begin{lem}    \label{inv vect1}
The category $\gmod_\chi$ contains a non-zero $(\ghat_{\ka_c},I)$
Harish-Chandra module if and only if
\begin{equation}    \label{tamely ram subscheme}
\chi \in \bigsqcup_{\nu \in P/W} \on{Op}_{^L
G}^{\on{RS}}(D)_{\varpi(\nu)}.
\end{equation}
\end{lem}

Thus, the opers $\chi$ for which the corresponding category
$\gmod_\chi$ contain non-trivial $I$-equivariant objects are precisely
the points of the subscheme \eqref{tamely ram subscheme} of
$\on{Op}_{^L G}(D^\times)$. The next question is what are the
corresponding $^L G$-local systems.

Let $\Loc_{^L G}^{\on{RS,uni}} \subset \Loc_{^L G}(D^\times)$ be the
locus of $^L G$-local systems on $D^\times$ with regular singularity
and unipotent monodromy. Such a local system is determined, up to an
isomorphism, by the conjugacy class of its monodromy (see, e.g.,
\cite{BV}, Sect. 8). Therefore
$\Loc_{^L G}^{\on{RS,uni}}$ is an algebraic stack isomorphic to ${\mc
N}/{}^L G$.  The following result is proved in a way similar to the
proof of \lemref{fiber over trivial}.

\begin{lem}
If the local system underlying an oper $\chi \in \on{Op}_{^L
G}(D^\times)$ belongs to $\Loc_{^L G}^{\on{RS,uni}}$, then $\chi$
belongs to the subset \eqref{tamely ram subscheme} of $\on{Op}_{^L
G}(D^\times)$.
\end{lem}

Indeed, the subscheme \eqref{tamely ram subscheme} is precisely the
(set-theoretic) preimage of $\Loc_{^L G}^{\on{RS,uni}} \subset
\Loc_{^L G}(D^\times)$ under the map $\al: \on{Op}_{^L G}(D^\times)
\to \Loc_{^L G}(D^\times)$.

This hardly comes as a surprise. Indeed, by analogy with the classical
Langlands correspondence we expect that the categories $\gmod_\chi$
containing non-trivial $I$-equivariant objects correspond to the
Langlands parameters which are the geometric counterparts of tamely
ramified homomorphisms $W'_F \to {}^L G$. The most obvious candidates
for those are precisely the $^L G$-local systems on $D^\times$ with
regular singularity and unipotent monodromy. For this reason we will
call such local systems {\bf tamely ramified}.

Let us summarize: suppose that $\sigma$ is a tamely ramified $^L
G$-local system on $D^\times$, and let $\chi$ be a $^L G$-oper that is
in the gauge equivalence class of $\sigma$. Then $\chi$ belongs to the
subscheme \eqref{tamely ram subscheme}, and the corresponding category
$\gmod_\chi$ contains non-zero $I$-equivariant objects, by \lemref{inv
vect1}. Let $\gmod_\chi^I$ be the corresponding category of
$I$-equivariant (or, equivalently, $(\ghat_{\ka_c},I)$ Harish-Chandra)
modules. Note that according to
\conjref{equiv loc syst}, the categories $\gmod_\chi$ (resp.,
$\gmod_\chi^I$) should be equivalent to each other for all $\chi$
which are gauge equivalent to each other as $^L G$-local systems.

In the next section, following \cite{FG:local}, we will give a
conjectural description of the categories $\gmod_\chi^I$ for $\chi \in
\on{Op}_{^L G}^{\on{RS}}(D)_{\varpi(-\rho)}$ in terms of the category
of coherent sheaves on the Springer fiber corresponding to the residue
of $\chi$. This description in particular implies that at least the
derived categories of these categories are equivalent to each other
for the opers corresponding to the same local system. We have a
similar conjecture for $\chi \in \on{Op}_{^L
G}^{\on{RS}}(D)_{\varpi(\nu)}$ for other $\nu \in P$, which the reader
may easily reconstruct from our discussion of the case $\nu=-\rho$.

\subsection{Conjectural description of the categories of
$(\ghat_{\ka_c},I)$ Harish-Chandra modules}    \label{conj descr}

Let us consider one of the connected components of the subscheme
\eqref{tamely ram subscheme}, namely, $\on{Op}_{^L
G}^{\on{RS}}(D)_{\varpi(-\rho)}$. Here it will be convenient to use a
different realization of this space, as the space $\nOp_{^L G}$ of
{\bf nilpotent opers} introduced in \cite{FG:local}. By definition, an
element of this space is an $^L N[[t]]$-gauge equivalence class of
operators of the form
\begin{equation}    \label{nilp oper mu1}
\nabla = \pa_t + p_{-1} + {\mb v}(t) + \frac{v}{t},
\end{equation}
where ${\mb v}(t) \in {}^L \bb[[t]]$ and $v \in {}^L \n$. It is shown
in \cite{FG:local} that $\nOp_{^L G} \simeq \on{Op}_{^L
G}^{\on{RS}}(D)_{\varpi(-\rho)}$. In particular, $\nOp_{^L G}$ is a
subspace of $\on{Op}_{^L G}(D^\times)$.

We have the (secondary) residue map
$$
\on{Res}: \nOp_{^L G} \to {}^L \n_{\F_{^L B,0}} = \F_{^L B,0}
\underset{^L B}\times {}^L \n,
$$
sending a gauge equivalence class of operators \eqref{nilp oper mu1}
to $v$. By abuse of notation, we will denote the corresponding map
$$
\nOp_{^L G} \to {}^L \n/{}^L B = \wt{\mc N}/{}^L G
$$
also by $\on{Res}$.

For any $\chi \in \nOp_{^L G}$ the $^L G$-gauge equivalence class of
the corresponding connection is a tamely ramified $^L G$-local system
on $D^\times$. Moreover, its monodromy conjugacy class is equal to
$\exp(2 \pi i \on{Res}(\chi))$.

We wish to describe the category $\gmod_\chi^I$ of $(\ghat_{\ka_c},I)$
Harish-Chandra modules with the central character $\chi \in \nOp_{^L
G}$. However, here we face the first major complication as compared to
the unramified case. While in the ramified case we worked with the
abelian category $\gmod_\chi^{G[[t]]}$, this does not seem to be
possible in the tamely ramified case. So from now on we will work with
the appropriate derived category $D^b(\gmod_\chi)^I$. By definition,
this is the full subcategory of the bounded derived category
$D^b(\gmod_\chi)$ whose objects are complexes with cohomologies in
$\gmod_\chi^I$.

Roughly speaking, the conjecture of \cite{FG:local} is that the
category $D^b(\gmod_\chi)^I$ is equivalent to the derived category
$D^b(\on{QCoh}(\on{Sp}_{\on{Res}(\chi)}))$ of the category
$\on{QCoh}(\on{Sp}_{\on{Res}(\chi)})$ of quasicoherent sheaves on
the Springer fiber of $\on{Res}(\chi)$. However, we need to make
some adjustments to this statement. These adjustments are needed to
arrive at a ``nice'' statement, \conjref{main conj FG} below. We
now explain what these adjustments are the reasons behind them.

The first adjustment is that we need to consider a slightly larger
category of representations than $D^b(\gmod_\chi)^I$. Namely, we wish
to include extensions of $I$-equivariant $\ghat_{\ka_c}$-modules which
are not necessarily $I$-equivariant, but only $I^0$-equivariant, where
$I^0 = [I,I]$. To explain this more precisely, let us choose a Cartan
subgroup $H \subset B \subset I$ and the corresponding Lie subalgebra
$\hh \subset \bb \subset \on{Lie} I$. We then have an isomorphism $I =
H \ltimes I^0$. An $I$-equivariant $\ghat_{\ka_c}$-module is the same
as a module on which $\hh$ acts diagonally with eigenvalues given by
integral weights and the Lie algebra $\on{Lie} I^0$ acts locally
nilpotently. However, there may exist extensions between such modules
on which the action of $\hh$ is no longer semi-simple. Such modules
are called $I$-{\bf monodromic}. More precisely, an $I$-monodromic
$\ghat_{\ka_c}$-module is a module that admits an increasing
filtration whose consecutive quotients are $I$-equivariant. It is
natural to include such modules in our category. However, it is easy
to show that an $I$-monodromic object of $\gmod_\chi$ is the same as
an $I^0$-equivariant object of $\gmod_\chi$ for any $\chi \in \nOp_{^L
G}$ (see \cite{FG:local}). Therefore instead of $I$-monodromic modules
we will use $I^0$-equivariant modules. Denote by
$D^b(\gmod_\chi)^{I^0}$ the the full subcategory of $D^b(\gmod_\chi)$
whose objects are complexes with the cohomologies in
$\gmod_\chi^{I^0}$.

The second adjustment has to do with the non-flatness of the Springer
resolution $\wt{\mc N} \to {\mc N}$. By definition, the Springer
fiber $\on{Sp}_u$ is the fiber product $\wt{\mc N}
\underset{{\mc N}}\times \on{pt}$, where $\on{pt}$ is the point $u
\in {\mc N}$. This means that the structure sheaf of $\on{Sp}_u$
is given by
\begin{equation}    \label{tensor product}
{\mc O}_{\on{Sp}_u} = {\mc O}_{\wt{\mc N}} \underset{{\mc O}_{{\mc
N}}}\otimes \C.
\end{equation}
However, because the morphism $\wt{\mc N} \to {\mc N}$ is not flat,
this tensor product functor is not left exact, and there are
non-trivial higher derived tensor products (the $Tor$'s). Our
(conjectural) equivalence is not going to be an exact functor: it
sends a general object of the category $\gmod_\chi^{I^0}$ not to an
object of the category of quasicoherent sheaves, but to a complex of
sheaves, or, more precisely, an object of the corresponding derived
category. Hence we are forced to work with derived categories, and so
the higher derived tensor products need to be taken into account.

To understand better the consequences of this non-exactness, let us
consider the following model example. Suppose that we have established
an equivalence between the derived category $D^b(\QCoh(\wt{\mc N}))$
and another derived category $D^b({\mc C})$. In particular, this means
that both categories carry an action of the algebra $\on{Fun} {\mc N}$
(recall that ${\mc N}$ is an affine algebraic variety). Let us suppose
that the action of $\on{Fun} {\mc N}$ on $D^b({\mc C})$ comes from its
action on the abelian category ${\mc C}$. Thus, ${\mc C}$ fibers over
${\mc N}$, and let ${\mc C}_u$ the fiber category corresponding to $u
\in {\mc N}$. This is the full subcategory of ${\mc C}$ whose objects
are objects of ${\mc C}$ on which the ideal of $u$ in $\on{Fun} {\mc
N}$ acts by $0$.\footnote{The relationship between ${\mc C}$ and ${\mc
C}_u$ is similar to the relationship between $\gmod$ and and
$\gmod_{\chi}$, where $\chi \in \on{Op}_{^L G}(D^\times)$.} What is
the category $D^b({\mc C}_u)$ equivalent to?

It is tempting to say that it is equivalent to
$D^b(\QCoh(\on{Sp}_u)$. However, this does not follow from the
equivalence of $D^b(\QCoh({\mc N}))$ and $D^b({\mc C})$ because of the
tensor product \eqref{tensor product} having non-trivial higher
derived functors. The correct answer is that $D^b({\mc C}_u)$ is
equivalent to the category $D^b(\QCoh(\on{Sp}^{\on{DG}}_u)$, where
$\on{Sp}^{\on{DG}}_u$ is the ``DG fiber'' of $\wt{\mc N} \to {\mc N}$
at $u$. By definition, a quasicoherent sheaf on $\on{Sp}^{\on{DG}}_u$
is a DG module over the DG algebra
\begin{equation}    \label{O Sp}
{\mc O}_{\on{Sp}^{\on{DG}}_u} = {\mc O}_{\wt{\mc N}}
\overset{L}{\underset{{\mc O}_{{\mc N}}}\otimes} \C_u,
\end{equation}
where we now take the full derived functor of tensor product. Thus,
$D^b(\QCoh(\on{Sp}^{\on{DG}}_u))$ may be thought of as the
derived category of quasicoherent sheaves on the ``DG scheme''
$\on{Sp}^{\on{DG}}_u$ (see \cite{CK} for a precise definition of DG
scheme).

Finally, the last adjustment is that we should consider the
non-reduced Springer fibers. This means that instead of the Springer
resolution $\wt{\mc N}$ we should consider the ``thickened''
Springer resolution
$$
\wt{\wt{\mc N}} = {} ^L \wt\g \underset{^L \g}\times {\mc N},
$$
where $^L \wt\g$ is the so-called {\bf Grothendieck alteration},
$$
^L \wt\g = \{ x \in {}^L \g, \bb' \in {}^L G/{}^L B \; | \; x \in \bb'
\}.
$$
The variety $\wt{\wt{\mc N}}$ is non-reduced, and the underlying
reduced variety is the Springer resolution $\wt{\mc N}$. For instance,
the fiber of $\wt{\mc N}$ over a regular element in ${\mc N}$ consists
of a single point, but the corresponding fiber of $\wt{\wt{\mc N}}$ is
the spectrum of the Artinian ring $h_0 = \on{Fun} {}^L \h/(\on{Fun}
{}^L \h)^W_+$. Here $(\on{Fun} {}^L \h)^W_+$ is the ideal in $\on{Fun}
{}^L \h$ generated by the augmentation ideal of the subalgebra of
$W$-invariants. Thus, $\on{Spec} h_0$ is the scheme-theoretic fiber of
$\varpi: {}^L \h \to {}^L \h/W$ at $0$. It turns out that in order to
describe the category $D^b(\gmod_\chi)^{I^0}$ we need to use the
``thickened'' Springer resolution.

Let us summarize: in order to construct the sought-after equivalence
of categories we take, instead of individual Springer fibers, the
whole Springer resolution, and we further replace it by the
``thickened'' Springer resolution $\wt{\wt{\mc N}}$ defined above.  In
this version we will be able to formulate our equivalence in such a
way that we avoid DG schemes.

This means that instead of considering the categories $\gmod_\chi$ for
individual nilpotent opers $\chi$, we should consider the
``universal'' category $\gmod_{\nilp}$ which is the ``family version''
of all of these categories. By definition, the category
$\gmod_{\nilp}$ is the full subcategory of $\gmod$ whose objects have
the property that the action of $Z(\ghat) = \on{Fun} \on{Op}_{^L
G}(D)$ on them factors through the quotient $\on{Fun} \on{Op}_{^L
G}(D) \to \on{Fun} \nOp_{^L G}$. Thus, the category $\gmod_{\nilp}$ is
similar to the category $\gmod_{\reg}$ that we have considered
above. While the former fibers over $\nOp_{^L G}$, the latter fibers
over $\on{Op}_{^L G}(D)$. The individual categories $\gmod_\chi$ are
now realized as fibers of these categories over particular opers
$\chi$.

Our naive idea was that for each $\chi \in \nOp_{^L G}$ the category
$D^b(\gmod_\chi)^{I^0}$ is equivalent to
$\on{QCoh}(\on{Sp}_{\on{Res}(\chi)})$. We would like to formulate
now a ``family version'' of such an equivalence. To this end we form
the fiber product
$$
^L \wt{\wt\n} = {}^L \wt\g \underset{^L\g}\times {}^L \n.
$$
It turns out that this fiber product does not suffer from the
problem of the individual Springer fibers, as the following lemma
shows:

\begin{lem}[\cite{FG:local},Lemma 6.4]    \label{L tensor}
The derived tensor product
$$
\on{Fun} {} ^L\wt\g \overset{L}{\underset{\on{Fun} {} ^L\g}\otimes}
\on{Fun} {}^L \n
$$
is concentrated in cohomological dimension $0$.
\end{lem}

The variety $^L \wt{\wt\n}$ may be thought of as the family of
(non-reduced) Springer fibers parameterized by $^L \n \subset {} ^L
\g$. It is important to note that it is singular, reducible and
non-reduced. For example, if $\g=\sw_2$, it has two components, one of
which is $\pone$ (the Springer fiber at $0$) and the other is the
doubled affine line (i.e., $\on{Spec} \C[x,y]/(y^2)$).

We note that the corresponding reduced scheme is
\begin{equation}    \label{wt L n}
^L \wt\n = \wt{\mc N} \underset{{\mc N}}\times {}^L \n.
\end{equation}
However, the derived tensor product corresponding to \eqref{wt L n} is
not concentrated in cohomological dimension $0$, and this is the
reason why we prefer to use $^L \wt{\wt\n}$ rather than $^L \wt\n$.

Now we set
$$
\nMOp = \nOp_{^L G} \underset{^L \n/{}^L
B}\times {}^L \wt{\wt{\n}}/{}^L B,
$$
where we use the residue morphism $\on{Res}: \nOp_{^L G} \to {}^L
\n/{}^L B$. Thus, informally $\nMOp$ may be thought as the
family over $\nOp_{^L G}$ whose fiber over $\chi \in \nOp_{^L
G}$ is the (non-reduced) Springer fiber of $\on{Res}(\chi)$.

The space $\nMOp$ is the space of {\bf Miura opers} whose underlying
opers are nilpotent, introduced in \cite{FG:local}.

We also introduce the category $\gmod_{\nilp}^{I^0}$ which is a full
subcategory of $\gmod_{\nilp}$ whose objects are
$I^0$-equivariant. Let $D^b(\gmod_{\nilp})^{I^0}$ be the corresponding
derived category.

Now we can formulate the Main Conjecture of \cite{FG:local}:

\begin{conj}    \label{main conj FG}
There is an equivalence of categories
\begin{equation}    \label{equiv FG}
D^b(\gmod_{\nilp})^{I^0} \simeq D^b(\QCoh(\nMOp))
\end{equation}
which is compatible with the action of the algebra $\on{Fun}
\nOp_{^L G}$ on both categories.
\end{conj}

Note that the action of $\on{Fun} \nOp_{^L G}$ on the first
category comes from the action of the center $Z(\ghat)$, and on the
second category it comes from the fact that $\nMOp$ is a
scheme over $\nOp_{^L G}$.

Another important remark is that the equivalence \eqref{equiv FG} does
not preserve the $t$-structures on the two categories. In other words,
\eqref{equiv FG} is expected in general to map objects of the abelian
category $\gmod_{\nilp}^{I^0}$ to complexes in $D^b(\QCoh(\nMOp))$,
and vice versa.

There are similar conjectures for the categories corresponding to the
spaces $\nOpla_{^L G}$ of nilpotent opers with dominant integral
weights $\la \in P^+$.

In the next section we will discuss the connection between
\conjref{main conj FG} and the classical tamely ramified Langlands
correspondence. We then present some evidence for this conjecture.

\subsection{Connection between the classical and the geometric
settings}

Let us discuss the connection between the equivalence
\eqref{equiv FG} and the realization of representations of affine
Hecke algebras in terms of $K$-theory of the Springer fibers. As we
have explained, we would like to view the category
$D^b(\gmod_{\chi})^{I^0}$ for $\chi \in
\nOp_{^L G}$ as, roughly, a categorification of the
space $R^{\pi(I)}$ of $I$-invariant vectors in an irreducible
representation $(R,\pi)$ of $G(F)$. Therefore, we expect that the
Grothendieck group of the category $D^b(\gmod_{\chi})^{I^0}$ is
somehow related to the space $R^{\pi(I)}$.

Let us try to specialize the statement of \conjref{main conj FG} to a
particular oper
$$
\chi = (\F,\nabla,\F_{^L B}) \in \nOp_{^L G}.
$$
Let $\wt{\on{Sp}}^{\on{DG}}_{\on{Res}(\chi)}$ be the DG fiber of
$\nMOp$ over $\chi$. By definition (see \secref{conj descr}), the
residue $\on{Res}(\chi)$ of $\chi$ is a vector in the twist of $^L \n$
by the $^L B$-torsor $\F_{^L B,0}$. It follows that
$\wt{\on{Sp}}^{\on{DG}}_{\on{Res}(\chi)}$ is the DG fiber over
$\on{Res}(\chi)$ of the $\F_{^L B,0}$-twist of the Grothendieck
alteration.

If we trivialize $\F_{^L B,0}$, then $u = \on{Res}(\chi)$ becomes
an element of $^L \n$. By definition, the (non-reduced) DG Springer
fiber $\wt{\on{Sp}}_u^{\on{DG}}$ is the DG fiber of the Grothendieck
alteration $^L \wt\g \to {}^L \g$ at $u$. In other words, the
corresponding structure sheaf is the DG algebra
$$
{\mc O}_{\wt{\on{Sp}}^{\on{DG}}_u} = {\mc O}_{^L \wt\g}
\overset{L}{\underset{{\mc O}_{^L\g}}\otimes} \C_u
$$
(compare with formula \eqref{O Sp}).

To see what these DG fibers look like, let $u=0$. Then the naive
Springer fiber is just the flag variety $^L G/{}^L B$ (it is reduced
in this case), and ${\mc O}_{\wt{\on{Sp}}_0}$ is the structure sheaf
of $^L G/{}^L B$. But the sheaf ${\mc O}_{\wt{\on{Sp}}^{\on{DG}}_0}$
is a sheaf of DG algebras, which is quasi-isomorphic to the complex of
differential forms on $^L G/{}^L B$, with the zero differential. In
other words, $\wt{\on{Sp}}^{\on{DG}}_0$ may be viewed as a
``$\Z$-graded manifold'' such that the corresponding supermanifold,
obtained by replacing the $\Z$-grading by the corresponding
$\Z/2\Z$-grading, is $\Pi T({}^L G/{}^L B)$, the tangent bundle to $^L
G/{}^L B$ with the parity of the fibers changed from even to odd.

We expect that the category $\gmod_{\nilp}^{I^0}$ is flat over
$\nOp_{^L G}$. Therefore, specializing \conjref{main conj FG} to a
particular oper $\chi \in \nOp_{^L G}$, we obtain as a corollary an
equivalence of categories
\begin{equation}    \label{nilp fiber}
D^b(\gmod_{\chi})^{I^0} \simeq
D^b(\QCoh(\wt{\on{Sp}}^{\on{DG}}_{\on{Res}(\chi)})).
\end{equation}
This bodes well with \conjref{equiv loc syst} saying that the
categories $\gmod_{\chi_1}$ and $\gmod_{\chi_2}$ (and hence
$D^b(\gmod_{\chi_1})^{I^0}$ and $D^b(\gmod_{\chi_2})^{I^0}$) should be
equivalent if the underlying local systems of the opers $\chi_1$ and
$\chi_2$ are isomorphic. For nilpotent opers $\chi_1$ and $\chi_2$
this is so if and only if their monodromies are conjugate to each
other. Since their monodromies are obtained by exponentiating their
residues, this is equivalent to saying that the residues,
$\on{Res}(\chi_1)$ and $\on{Res}(\chi_2)$, are conjugate
with respect to the $\F_{^L B,0}$-twist of $^L G$. But in this case
the DG Springer fibers corresponding to $\chi_1$ and $\chi_2$ are also
isomorphic, and so $D^b(\gmod_{\chi_1})^{I^0}$ and
$D^b(\gmod_{\chi_2})^{I^0}$ are equivalent to each other, by
\eqref{nilp fiber}.

The Grothendieck group of the category
$D^b(\QCoh(\wt{\on{Sp}}^{\on{DG}}_{u}))$, where $u$ is any nilpotent
element, is the same as the Grothendieck group of
$\QCoh(\on{Sp}_{u})$. In other words, the Grothendieck group does
not ``know'' about the DG or the non-reduced structure of
$\wt{\on{Sp}}^{\on{DG}}_{u}$. Hence it is nothing but the algebraic
$K$-theory $K(\on{Sp}_{u})$. As we explained at the end of
\secref{tame reps}, equivariant variants of this algebraic $K$-theory
realize the ``standard modules'' over the affine Hecke algebra
$H(G(F),I)$. Moreover, the spaces of $I$-invariant vectors
$R^{\pi(I)}$ as above, which are naturally modules over the affine
Hecke algebra, may be realized as subquotients of
$K(\on{Sp}_{u})$. This indicates that the equivalences \eqref{nilp
fiber} and \eqref{equiv FG} are compatible with the classical results.

However, at first glance there are some important differences between
the classical and the categorical pictures, which we now discuss in
more detail.

In the construction of $H(G(F),I)$-modules outlined in \secref{tame
reps} we had to pick a semi-simple element $\ga$ of $^L G$ such that
$\ga u \ga^{-1} = qu$, where $q$ is the number of elements in the
residue field of $F$. Then we consider the specialized
$A$-equivariant $K$-theory $K^A(\on{Sp}_u)_{(\ga,q)}$
where $A$ is the the smallest algebraic subgroup of $^L G \times
\C^\times$ containing $(\ga,q)$. This gives $K(\on{Sp}_{u})$ the
structure of an $H(G(F),I)$-module. But this module carries a residual
symmetry with respect to the group $C(\ga,u)$ of components
of the centralizer of $\ga$ and $u$ in $^L G$, which commutes with the
action of $H(G(F),I)$. Hence we consider the $H(G(F),I)$-module
$$
K^A(\on{Sp}_u)_{(\ga,q,\rho)} =
\on{Hom}_{C(\ga,u)}(\rho,K(\on{Sp}_{u})),
$$
corresponding to an irreducible representation $\rho$ of
$C(\ga,u)$. Finally, each of these components has a unique irreducible
quotient, and this is an irreducible representation of $H(G(F),I)$
which is realized on the space $R^{\pi(I)}$, where $(R,\pi)$ is an
irreducible representation of $G(F)$ corresponding to $(\ga,u,\rho)$
under the bijection
\eqref{tamely ramified langlands2}. How is this intricate structure
reflected in the categorical setting?

Our category $D^b(\QCoh(\wt{\on{Sp}}^{\on{DG}}_{u}))$, where
$u=\on{Res}(\chi)$, is a particular categorification of the
(non-equivariant) $K$-theory $K(\on{Sp}_{u})$. Note that in the
classical local Langlands correspondence \eqref{tamely ramified
langlands2} the element $u$ of the triple $(\ga,u,\rho)$ is
interpreted as the logarithm of the monodromy of the corresponding
representation of the Weil-Deligne group $W'_F$. This is in agreement
with the interpretation of $\on{Res}(\chi)$ as the logarithm of
the monodromy of the $^L G$-local system on $D^\times$ corresponding
to $\chi$, which plays the role of the local Langlands parameter for
the category $\gmod_\chi$ (up to the inessential factor $2\pi i$).

But what about the other parameters, $\ga$ and $\rho$? And why does
our category correspond to the non-equivariant $K$-theory of the
Springer fiber, and not the equivariant $K$-theory, as in the
classical setting?

The element $\ga$ corresponding to the Frobenius in $W'_F$ does not
seem to have an analogue in the geometric setting. We have already
seen this above in the unramified case: while in the classical setting
unramified local Langlands parameters are the semi-simple conjugacy
classes $\ga$ in $^L G$, in the geometric setting we have only one
unramified local Langlands parameter, namely, the trivial local
system.

To understand better what's going on here, we revisit the unramified
case. Recall that the spherical Hecke algebra $H(G(F),K_0)$ is
isomorphic to the representation ring $\on{Rep} {}^L G$. The
one-dimensional space of $K_0$-invariants in an irreducible
unramified representation $(R,\pi)$ of $G(F)$ realizes a one-dimensional
representation of $H(G(F),K_0)$, i.e., a homomorphism $\on{Rep} {}^L G
\to \C$. The unramified Langlands parameter $\ga$ of $(R,\pi)$, which
is a semi-simple conjugacy class in $^L G$, is the point in
$\on{Spec}(\on{Rep} {}^L G)$ corresponding to this homomorphism. What
is a categorical analogue of this homomorphism? The categorification
of $\on{Rep} {}^L G$ is the category ${\mc Rep} \, {}^L G$. The
product structure on $\on{Rep} {}^L G$ is reflected in the structure
of tensor category on ${\mc Rep}
\, {}^L G$. On the other hand, the categorification of the algebra
$\C$ is the category ${\mc Vect}$ of vector spaces. Therefore a
categorical analogue of a homomorphism $\on{Rep} {}^L G
\to \C$ is a functor ${\mc Rep} \, {}^L G \to
{\mc Vect}$ respecting the tensor structures on both categories. Such
functors are called the fiber functors. The fiber functors form a
category of their own, which is equivalent to the category of $^L
G$-torsors. Thus, any two fiber functors are isomorphic, but not
canonically. In particular, the group of automorphisms of each fiber
functor is isomorphic to $^L G$. (Incidentally, this is how $^L G$ is
reconstructed from a fiber functor in the Tannakian formalism.) Thus,
we see that while in the categorical world we do not have analogues of
semi-simple conjugacy classes $\ga$ (the points of $\on{Spec}(\on{Rep}
{}^L G)$), their role is in some sense played by the group of
automorphisms of a fiber functor.

This is reflected in the fact that while in the categorical setting we
have a unique unramified Langlands parameter, namely, the trivial $^L
G$-local system $\sigma_0$ on $D^\times$, this local system has a
non-trivial group of automorphisms, namely, $^L G$. We therefore
expect that the group $^L G$ should act by automorphisms of the
Langlands category ${\mc C}_{\sigma_0}$ corresponding to $\sigma_0$,
and this action should commute with the action of the loop group
$G\ppart$ on ${\mc C}_{\sigma_0}$. It is this action of $^L G$ that is
meant to compensate for the lack of unramified Langlands parameters,
as compared to the classical setting.

We have argued in \secref{unram case} that the category $\gmod_\chi$,
where $\chi = (\F,\nabla,\F_{^L B}) \in \on{Op}_{^L G}(D)$, is a
candidate for the Langlands category ${\mc C}_{\sigma_0}$. Therefore
we expect that the group $^L G$ (more precisely, its twist $^L
G_{\F}$) acts on the category $\gmod_\chi$. In \secref{eq of cat} we
showed how to obtain this action using the conjectural equivalence
between $\gmod_\chi$ and the category ${\mc
D}_{\ka_c}^{\Hecke}\mod_\chi$ of Hecke eigenmodules on the affine
Grassmannian $\on{Gr}$ (see \conjref{equiv with Dmod}). The category
${\mc D}_{\ka_c}^{\Hecke}\mod_\chi$ was defined in \secref{reps and
Dmod} as a ``de-equivariantization'' of the category ${\mc
D}_{\ka_c}\mod$ of twisted ${\mc D}$-modules on $\on{Gr}$ with respect
to the monoidal action of the category ${\mc Rep} \, {}^L G$.

Now comes a crucial observation which will be useful for understanding
the way things work in the tamely ramified case: the category ${\mc
Rep} \, {}^L G$ may be interpreted as the category of $^L
G$-equivariant quasicoherent sheaves on the variety $\on{pt} =
\on{Spec} \C$. In other words, ${\mc Rep} \, {}^L G$ may be
interpreted as the category of quasicoherent sheaves on the stack
$\on{pt}/{}^L G$. The existence of monoidal action of the category
${\mc Rep} \, {}^L G$ on ${\mc D}_{\ka_c}\mod$ should be viewed as the
statement that the category ${\mc D}_{\ka_c}\mod$ ``lives'' over the
stack $\on{pt}/{}^L G$. The statement of \conjref{equiv with Dmod} may
then be interpreted as saying that
$$
\gmod_\chi \simeq {\mc D}_{\ka_c}\mod \underset{\on{pt}/{}^L G}\times
\on{pt}.
$$
In other words, if ${\mc C}$ is the conjectural Langlands category
fibering over the stack $\Loc_{^L G}(D^\times)$ of all $^L G$-local
systems on $D^\times$, then
$$
{\mc D}_{\ka_c}\mod \simeq {\mc C} \underset{\Loc_{^L
G}(D^\times)}\times \on{pt}/{}^L G,
$$
whereas
$$
\gmod_\chi \simeq {\mc C} \underset{\Loc_{^L G}(D^\times)}\times
\on{pt},
$$
where the morphism $\on{pt} \to \Loc_{^L G}(D^\times)$ corresponds to
the oper $\chi$.

Thus, in the categorical setting there are two different ways to think
about the trivial local system $\sigma_0$: as a point (defined by a
particular $^L G$-bundle on $D$ with connection, such as a regular
oper $\chi$), or as a stack $\on{pt}/{}^L G$. The base change of the
Langlands category in the first case gives us a category with an
action of $^L G$, such as the categories $\gmod_\chi$ or ${\mc
D}_{\ka_c}^{\Hecke}\mod$. The base change in the second case gives us
a category with a monoidal action of ${\mc Rep} \, {}^L G$, such as
the category ${\mc D}_{\ka_c}\mod$. We can go back and forth between
the two by applying the procedures of equivariantization and
de-equivariantization with respect to $^L G$ and ${\mc Rep} \, {}^L
G$, respectively.

Now we return to the tamely ramified case. The semi-simple element
$\ga$ appearing in the triple $(\ga,u,\rho)$ plays the same role as
the unramified Langlands parameter $\ga$. However, now it must satisfy
the identity $\ga u \ga^{-1} = qu$. Recall that the center $Z$ of
$H(G(F),I)$ is isomorphic to $\on{Rep} {}^L G$, and so $\on{Spec} Z$
is the set of all semi-simple elements in $^L G$. For a fixed
nilpotent element $u$ the equation $\ga u \ga^{-1} = qu$ cuts out a
locus $C_u$ in $\on{Spec} Z$ corresponding to those central characters
which may occur on irreducible $H(G(F),I)$-modules corresponding to
$u$. In the categorical setting (where we set $q=1$) the analogue of
$C_u$ is the centralizer $Z(u)$ of $u$ in $^L G$, which is precisely
the group $\on{Aut}(\sigma)$ of automorphisms of a tame local system
$\sigma$ on $D^\times$ with monodromy $\exp(2\pi iu)$. On general
grounds we expect that the group $\on{Aut}(\sigma)$ acts on the
Langlands category ${\mc C}_\sigma$, just as we expect the group $^L
G$ of automorphisms of the trivial local system $\sigma_0$ to act on
the category ${\mc C}_{\sigma_0}$. It is this action that replaces the
parameter $\ga$ in the geometric setting.

In the classical setting we also have one more parameter, $\rho$. Let
us recall that $\rho$ is a representation of the group $C(\ga,u)$ of
connected components of the centralizer $Z(\ga,u)$ of $\ga$ and
$u$. But the group $Z(\ga,u)$ is a subgroup of $Z(u)$, which becomes
the group $\on{Aut}(\sigma)$ in the geometric setting. Therefore one
can argue that the parameter $\rho$ is also absorbed into the action
of $\on{Aut}(\sigma)$ on the category ${\mc C}_\sigma$.

If we have an action of $\on{Aut}(\sigma)$ on the category ${\mc
C}_\sigma$, or on one of its many incarnations $\gmod_\chi, \chi \in
\nOp_{^L G}$, it means that these categories must be
``de-equivariantized'', just like the categories $\gmod_\chi, \chi \in
\on{Op}_{^L G}(D)$, in the unramified case. This is the reason why in
the equivalence \eqref{nilp fiber} (and in \conjref{main conj FG}) we
have the non-equivariant categories of quasicoherent sheaves (whose
Grothendieck groups correspond to the non-equivariant $K$-theory of
the Springer fibers).

However, there is also an equivariant version of these
categories. Consider the substack of tamely ramified local systems in
$\Loc_{^L G}(D^\times)$ introduced in \secref{admitting I}. Since a
tamely ramified local system is completely determined by the logarithm
of its (unipotent) monodromy, this substack is isomorphic to ${\mc
N}/{}^L G$. This substack plays the role of the substack $\on{pt}/{}^L
G$ corresponding to the trivial local system. Let us set
$$
{\mc C}_{\on{tame}} = {\mc C} \underset{\Loc_{^L
G}(D^\times)}\times {\mc N}/{}^L G.
$$
Then, according to our general conjecture expressed by the Cartesian
diagram \eqref{base change}, we expect to have
\begin{equation}    \label{base change1}
\gmod_{\nilp} \simeq {\mc C}_{\on{tame}} \underset{{\mc N}/{}^L
G}\times \nOp_{^L G}.
\end{equation}
Let $D^b({\mc C}_{\on{tame}})^{I^0}$ be the $I^0$-equivariant derived
category corresponding to ${\mc C}_{\on{tame}}$. Combining \eqref{base
change1} with
\conjref{main conj FG}, and noting that
$$
\nMOp \simeq \nOp_{^L G} \underset{{\mc N}/{}^L
G}\times \wt{\wt{{\mc N}}}/{}^L G,
$$
we obtain the following conjecture (see \cite{FG:local}):
\begin{equation}    \label{tame eq}
D^b({\mc C}_{\on{tame}})^{I^0} \simeq D^b(\QCoh(\wt{\wt{{\mc N}}}/{}^L
G)).
\end{equation}
The category on the right hand side may be interpreted as the derived
category of $^L G$-equivariant quasicoherent sheaves on the
``thickened'' Springer resolution $\wt{\wt{{\mc N}}}$.

Together, the conjectural equivalences \eqref{nilp fiber} and
\eqref{tame eq} should be thought of as the categorical versions of
the realizations of modules over the affine Hecke algebra in the
$K$-theory of the Springer fibers.

One corollary of the equivalence \eqref{nilp fiber} is the following:
the classes of irreducible objects of the category $\gmod_\chi^{I_0}$
in the Grothendieck group of $\gmod_\chi^{I_0}$ give rise to a basis
in the algebraic $K$-theory $K(\on{Sp}_u)$, where
$u=\on{Res}(\chi)$. Presumably, this basis is closely related to
the bases in (equivariant version of) this $K$-theory constructed by
G. Lusztig in \cite{Lus1} (from the perspective of unrestricted
$\g$-modules in positive characteristic).

\subsection{Evidence for the conjecture}    \label{evidence}

We now describe some evidence for \conjref{main conj FG}. It consists
of the following four groups of results:

\begin{itemize}

\item Interpretation of the Wakimoto modules as
$\ghat_{\ka_c}$-modules corresponding to the skyscraper sheaves on
$\nMOp$;

\item Connection to R. Bezrukavnikov's theory;

\item Proof of the equivalence of certain quotient categories of
$D^b(\gmod_{\nilp})^{I^0}$ and $D^b(\QCoh(\nMOp))$, \cite{FG:local}.

\item Proof of the restriction of the equivalence \eqref{equiv FG} to
regular opers, \cite{FG:equiv}.

\end{itemize}

We start with the discussion of Wakimoto modules.

Suppose that we have proved the equivalence of categories \eqref{equiv
FG}. Then each quasicoherent sheaf on $\nMOp$ should correspond to an
object of the derived category $D^b(\gmod_{\nilp})^{I^0}$. The
simplest quasicoherent sheaves on $\nMOp$ are the {\bf skyscraper
sheaves} supported at the $\C$-points of $\nMOp$. It follows from the
definition that a $\C$-point of $\nMOp$, which is the same as a
$\C$-point of the reduced scheme $\on{MOp}_{G}^{0}$, is a pair
$(\chi,\bb')$, where $\chi = (\F,\nabla,\F_{^L B})$ is a nilpotent $^L
G$-oper in $\nOp_{^L G}$ and $\bb'$ is a point of the Springer fiber
corresponding to $\on{Res}(\chi)$, which is the variety of Borel
subalgebras in $^L \g_{\F_0}$ that contain
$\on{Res}(\chi)$. Thus, if
\conjref{main conj FG} is true, we should have a family of objects of
the category $D^b(\gmod_{\nilp})^{I^0}$ parameterized by these
data. What are these objects?

The answer is that these are the {\bf Wakimoto modules}. These modules
were originally introduced by M. Wakimoto \cite{Wak} for $\g=\sw_2$
and by B. Feigin and myself in general in \cite{FF:usp,FF:si} (see
also \cite{F:wak}). We recall from \cite{F:wak} that Wakimoto modules
of critical level are parameterized by the space $\Con_{D^\times}$ of
connections on the $^L H$-bundle $\Omega^{-\rho}$ over
$D^\times$. This is the push-forward of the $\C^\times$-bundle
corresponding to the canonical line bundle $\Omega$ with respect to
the homomorphism $\rho: \C^\times \to {}^L H$. Let us denote the
Wakimoto module corresponding to $\ol\nabla \in \Con_{D^\times}$ by
$W_{\ol\nabla}$. According to \cite{F:wak}, Theorem 12.6, the center
$Z(\ghat)$ acts on $W_{\ol\nabla}$ via the central character
$\mu(\ol\nabla)$, where
$$
\mu: \Con_{D^\times} \to \on{Op}_{^L G}(D^\times)
$$
is the {\bf Miura transformation}.

It is not difficult to show that if $\chi \in \nOp_{^L G}$, then
$W_{\ol\nabla}$ is an object of the category $\gmod^I_\chi$ for any
$\ol\nabla \in \mu^{-1}(\chi)$. Now, according to the results
presented in \cite{FG:local}, the points of the fiber $\mu^{-1}(\chi)$
of the Miura transformation over $\chi$ are in bijection with the
points of the Springer fiber $\on{Sp}_{\on{Res}(\chi)}$ corresponding
to the nilpotent element $\on{Res}(\chi)$. Therefore to each point of
$\on{Sp}_{\on{Res}(\chi)}$ we may assign a Wakimoto module, which is
an object of the category $\gmod^{I^0}_\chi$ (and hence of the
corresponding derived category). In other words, Wakimoto modules are
objects of the category $\gmod^I_{\nilp}$ parameterized by the
$\C$-points of $\nMOp$. It is natural to assume that they correspond
to the skyscraper sheaves on $\nMOp$ under the equivalence
\eqref{equiv FG}. This was in fact one of our motivations for this
conjecture.

Incidentally, this gives us a glimpse into how the group of
automorphisms of the $^L G$-local system underlying the oper $\chi$
acts on the category $\gmod_\chi$. This group is $Z(\on{Res}(\chi))$,
the centralizer of the residue $\on{Res}(\chi)$, and it acts on the
Springer fiber $\on{Sp}_{\on{Res}(\chi)}$. Therefore $g \in
Z(\on{Res}(\chi))$ sends the skyscraper sheaf supported at a point $p
\in \on{Sp}_{\on{Res}(\chi)}$ to the skyscraper sheaf supported at $g
\cdot p$. Thus, we expect that $g$ sends the Wakimoto module
corresponding to $p$ to the Wakimoto module corresponding to $g \cdot
p$.

If the Wakimoto modules indeed correspond to the skyscraper sheaves,
then the equivalence \eqref{equiv FG} may be thought of as a kind of
``spectral decomposition'' of the category $D^b(\gmod_{\nilp})^{I^0}$,
with the basic objects being the Wakimoto modules $W_{\ol\nabla}$,
where $\ol\nabla$ runs over the locus in $\Con_{D^\times}$ which is
isomorphic, pointwise, to $\nMOp$ (see \cite{FG:wak} for more
details).

Now we discuss the second piece of evidence, connection with
Bezruka\-ni\-kov's theory.

To motivate it, let us recall that in \secref{action of sph hecke} we
discussed the action of the categorical spherical algebra ${\mc
H}(G\ppart,G[[t]])$ on the category $\gmod_\chi$, where $\chi$ is a
regular oper. The affine Hecke algebra $H(G(F),I)$ also has a
categorical analogue. Consider the {\bf affine flag variety} $\on{Fl}
= G\ppart/I$. The categorical Hecke algebra is the category ${\mc
H}(G\ppart,I)$ which is the full subcategory of the derived category
of ${\mc D}$-modules on $\on{Fl} = G\ppart/I$ whose objects are
complexes with $I$-equivariant cohomologies. This category naturally
acts on the derived category $D^b(\gmod_\chi)^I$. What does this
action correspond to on the other side of the equivalence \eqref{equiv
FG}?

The answer is given by a theorem of R. Bezrukavnikov \cite{Bez}, which
may be viewed as a categorification of the isomorphism \eqref{KL}:
\begin{equation}    \label{bez iso}
D^b({\mc D}^{\on{Fl}}_{\ka_c}\mod)^{I^0} \simeq
D^b(\QCoh(\wt{\on{St}})),
\end{equation}
where ${\mc D}^{\on{Fl}}_{\ka_c}\mod$ is the category of twisted ${\mc
D}$-modules on $\on{Fl}$ and $\wt{\on{St}}$ is the ``thickened''
Steinberg variety
$$
\wt{\on{St}} = \wt{\mc N} \underset{{\mc N}}\times \wt{\wt{\mc N}} =
\wt{\mc N} \underset{^L\g}\times {}^L \wt{\g}.
$$

Morally, we expect that the two categories in \eqref{bez iso} act on
the two categories in \eqref{equiv FG} in a compatible way. However,
strictly speaking, the left hand side of
\eqref{bez iso} acts like this:
$$
D^b(\gmod_{\nilp})^I \to D^b(\gmod_{\nilp})^{I^0},
$$
and the right hand side of \eqref{bez iso} acts like this:
$$
D^b(\QCoh(\on{MOp}^0_{^L G})) \to D^b(\QCoh(\nMOp)).
$$
So one needs a more precise statement, which may be found in
\cite{Bez}, Sect. 4.2. Alternatively, one can consider the
corresponding actions of the affine braid group of $^L G$, as in
\cite{Bez}.

A special case of this compatibility concerns some special objects of
the category $D^b({\mc D}^{\on{Fl}}_{\ka_c}\mod)^{I}$, the central
sheaves introduced in \cite{Ga:central}. They correspond to the
central elements of the affine Hecke algebra $H(G(F),I)$. These
central elements act as scalars on irreducible $H(G(F),I)$-modules, as
well as on the standard modules $K^A(\on{Sp}_u)_{(\ga,q,\rho)}$
discussed above. We have argued that the categories $\gmod_\chi^{I^0},
\chi \in \nOp_{^L G}$, are categorical versions of these
representations. Therefore it is natural to expect that its objects
are ``eigenmodules'' with respect to the action of the central sheaves
from $D^b({\mc D}^{\on{Fl}}_{\ka_c}\mod)^{I}$ (in the sense of
\secref{action of sph hecke}). This has indeed been proved in
\cite{FG:fusion}.

This discussion indicates an intimate connection between the
category $D^b(\gmod_{\nilp})$ and the category of twisted ${\mc
D}$-modules on the affine flag variety, which is similar to the
connection between $\gmod_{\reg}$ and the category of twisted ${\mc
D}$-modules on the affine Grassmannian which we discussed in
\secref{reps and Dmod}. A more precise conjecture relating
$D^b(\gmod_{\nilp})$ and $D^b({\mc D}^{\on{Fl}}_{\ka_c}\mod)$ was
formulated in \cite{FG:local} (see the Introduction and Sect. 6),
where we refer the reader for more details. This conjecture may be
viewed as an analogue of \conjref{equiv with Dmod} for nilpotent
opers. As explained in \cite{FG:local}, this conjecture is supported
by the results of \cite{AB,ABG} (see also \cite{Bez1,Bez}). Together,
these results and conjectures provide additional evidence for the
equivalence \eqref{equiv FG}.

\section{Ramified global Langlands correspondence}
\label{global}

We now discuss the implications of the local Langlands correspondence
for the global geometric Langlands correspondence.

We begin by briefly discussing the setting of the classical global
Langlands correspondence (see \cite{F:rev} for more details).

\subsection{The classical setting}

Let $X$ be a smooth projective curve over $\Fq$. Denote by $F$ the
field $\Fq(X)$ of rational functions on $X$. For any closed point $x$
of $X$ we denote by $F_x$ the completion of $F$ at $x$ and by $\OO_x$
its ring of integers. If we choose a local coordinate $t_x$ at $x$
(i.e., a rational function on $X$ which vanishes at $x$ to order one),
then we obtain isomorphisms $F_x \simeq {\mathbb F}_{q_x}((t_x))$ and
$\OO_x \simeq {\mathbb F}_{q_x}[[t_x]]$, where ${\mathbb F}_{q_x}$ is
the residue field of $x$; in general, it is a finite extension of
$\Fq$ containing $q_x = q^{\on{ord}_x}$ elements.

Thus, we now have a local field attached to each point of $X$.  The
ring $\AD=\AD_F$ of {\bf ad\`eles} of $F$ is by definition the {\bf
restricted} product of the fields $F_x$, where $x$ runs over the set
$|X|$ of all closed points of $X$. The word ``restricted'' means that
we consider only the collections $(f_x)_{x \in |X|}$ of elements of
$F_x$ in which $f_x \in \OO_x$ for all but finitely many $x$. The ring
$\AD$ contains the field $F$, which is embedded into $\AD$ diagonally,
by taking the expansions of rational functions on $X$ at all points.

While in the local Langlands correspondence we considered irreducible
smooth representations of the group $GL_n$ over a local field, in the
global Langlands correspondence we consider irreducible {\bf
automorphic representations} of the group $GL_n(\AD)$. The word
``automorphic'' means, roughly, that the representation may be
realized in a reasonable space of functions on the quotient
$GL_n(F)\bs GL_n(\AD)$ (on which the group  $GL_n(\AD)$ acts from the
right).

On the other side of the correspondence we consider $n$-dimensional
representations of the Galois group $\on{Gal}(\ol{F}/F)$, or, more
precisely, the Weil group $W_F$, which is a subgroup of
$\on{Gal}(\ol{F}/F)$ defined in the same way as in the local case.

Roughly speaking, the global Langlands correspondence is a
bijection between the set of equivalence classes of $n$-dimensional
representations of $W_F$ and the set of equivalence classes of
irreducible automorphic representations of $GL_n(\AD)$:

$$
\boxed{\begin{matrix} n\text{-dimensional representations} \\
    \text{of } W_F \end{matrix}} \quad \Longleftrightarrow
    \quad \boxed{\begin{matrix} \text{irreducible automorphic} \\
	\text{representations of } GL_n(\AD)
\end{matrix}}
$$
\bigskip

The precise statement is more subtle. For example, we should consider
the so-called $\ell$-adic representations of the Weil group (while in
the local case we considered the admissible complex representations of
the Weil-Deligne group; the reason is that in the local case those are
equivalent to the $\ell$-adic representations). Moreover, under this
correspondence important invariants attached to the objects appearing
on both sides (Frobenius eigenvalues on the Galois side and the Hecke
eigenvalues on the other side) are supposed to match. We refer the
reader to Part I of the review \cite{F:rev} for more details.

The global Langlands correspondence has been proved for $GL_2$ in the
80's by V. Drinfeld \cite{Dr1}--\cite{Dr4} and more recently by
L. Lafforgue \cite{Laf} for $GL_n$ with an arbitrary $n$.

Like in the local story, we may also wish to replace the group $GL_n$
by an arbitrary reductive algebraic group defined over $F$. Then on
one side of the global Langlands correspondence we have homomorphisms
$\sigma: W_F \to {}^L G$ satisfying some properties (or perhaps, some
more refined data, as in \cite{Arthur}). We expect to be able to
attach to each $\sigma$ an {\bf automorphic representation} $\pi$ of
$GL_n(\AD_F)$.\footnote{In this section, by abuse of notation, we will
use the same symbol to denote a representation of the group and the
vector space underlying this representation.} The word ``automorphic''
again means, roughly, that the representation may be realized in a
reasonable space of functions on the quotient $GL_n(F)\bs GL_n(\AD)$
(on which the group $GL_n(\AD)$ acts from the right). We will not try
to make this precise. In general, we expect not one but several
automorphic representations assigned to $\sigma$ which are the global
analogues of the $L$-packets discussed above (see
\cite{Arthur}). Another complication is that the multiplicity of a
given irreducible automorphic representation in the space of functions
on $GL_n(F)\bs GL_n(\AD)$ may be greater than one. We will mostly
ignore all of these issues here, as our main interest is in the
geometric theory (note also that these issues do not arise if
$G=GL_n$).

An irreducible automorphic representation may always be decomposed as
the restricted tensor product $\bigotimes'_{x \in X} \pi_x$, where
each $\pi_x$ is an irreducible representation of $G(F_x)$. Moreover,
for all by finitely many $x \in X$ the factor $\pi_x$ is an {\bf
unramified} representation of $G(F_x)$: it contains a non-zero vector
invariant under the maximal compact subgroup $K_{0,x}=G(\OO_x)$ (see
\secref{unr rep}). Let us choose such a vector $v_x \in \pi_x$ (it is
unique up to a scalar). The word ``restricted'' means that we consider
the span of vectors of the form $\otimes_{x \in X} u_x$, where $u_x
\in \pi_x$ and $u_x = v_x$ for all but finitely many $x \in X$.

An important property of the global Langlands correspondence is its
compatibility with the local one. We can embed the Weil group
$W_{F_x}$ of each of the local fields $F_x$ into the global Weil group
$W_F$. Such an embedding is not unique, but it is well-defined up to
conjugation in $W_F$. Therefore an equivalence class of $\sigma: W_F
\to {}^L G$ gives rise to a well-defined equivalence class of
$\sigma_x: W_{F_x} \to {}^L G$. We will impose the condition on
$\sigma$ that for all but finitely many $x \in X$ the homomorphism
$\sigma_x$ is unramified (see \secref{unr rep}).

By the local Langlands correspondence, to $\sigma_x$ one can attach an
equivalence class of irreducible smooth representations $\pi_x$ of
$G(F_x)$.\footnote{Here we are considering $\ell$-adic homomorphisms
from the Weil group $W_{F_x}$ to $^L G$, and therefore we do not need
to pass from the Weil group to the Weil-Deligne group.} Moreover, an
unramified $\sigma_x$ will correspond to an unramified irreducible
representation $\pi_x$. The compatibility between local and global
correspondences is the statement that the automorphic representation
of $G(\AD)$ corresponding to $\sigma$ should be isomorphic to the
restricted tensor product $\bigotimes'_{x \in X} \pi_x$.
Schematically, this is represented as follows:
\begin{align*}
\sigma &\overset{\on{global}}\longleftrightarrow \pi = \bigotimes_{x
\in X}{}' \pi_x \\ \sigma_x &\overset{\on{local}}\longleftrightarrow
\pi_x.
\end{align*}

In this section we discuss an analogue of this local-to-global
principle in the geometric setting and the implications of our
local results and conjectures for the global geometric Langlands
correspondence. We focus in particular on the unramified and tamely
ramified Langlands parameters. At the end of the section we also
discuss connections with irregular singularities.

\subsection{The unramified case, revisited}    \label{outside the
  locus}

An important special case is when $\sigma: W_F \to {}^L G$ is
everywhere unramified. Then for each $x \in X$ the corresponding
homomorphism $\sigma_x: W_{F_x} \to {}^L G$ is unramified, and hence
corresponds, as explained in \secref{unr rep}, to a semi-simple
conjugacy class $\ga_x$ in $^L G$, which is the image of the Frobenius
element under $\sigma_x$. This conjugacy class in turn gives rise to
an unramified irreducible representation $\pi_x$ of $G(F_x)$ with a
unique, up to a scalar, vector $v_x$ such that $G(\OO_x) v_x =
v_x$. The spherical Hecke algebra $H(G(F_x),G(\OO_x)) \simeq \on{Rep}
{}^L G$ acts on this vector according to formula \eqref{Hecke
eigenfunction}:
\begin{equation}    \label{Hecke eigenfunction1}
H_{V,x} \star v_x = \on{Tr}(\ga_x,V) v_x, \qquad [V] \in
\on{Rep} {}^L G.
\end{equation}
The tensor product $v = \otimes_{x \in X}$ of this vectors is a
$G(\OO)$-invariant vector in $\pi = \bigotimes'_{x \in X} \pi_x$,
which, according to the global Langlands conjecture is
automorphic. This means that $\pi$ is realized in the space of
functions on $G(F)\bs G(\AD_F)$. In this realization vector $v$
corresponds to a right $G(\OO)$-invariant function on $G(F)\bs
G(\AD_F)$, or equivalently, a function on the double quotient
\begin{equation}    \label{double quotient}
G(F)\bs G(\AD_F)/G(\OO).
\end{equation}

Thus, an unramified global Langlands parameter $\sigma$ gives rise to
a function on \eqref{double quotient}. This function is the {\bf
automorphic function} corresponding to $\sigma$. We denote it by
$f_\pi$. Since it corresponds to a vector in an irreducible
representation $\pi$ of $G(\AD_F)$, the entire representation $\pi$
may be reconstructed from this function. Thus, we do not lose any
information by passing from $\pi$ to $f_\pi$.

Since $v \in \pi$ is an eigenvector of the Hecke operators, according
to formula \eqref{Hecke eigenfunction1}, we obtain that the function
$f_\pi$ is a {\bf Hecke eigenfunction} on the double quotient
\eqref{double quotient}. In fact, the local Hecke algebras
$H(G(F_x),G(\OO_x))$ act naturally (from the right) on the space of
functions on \eqref{double quotient}, and $f_\pi$ is an
eigenfunction of this action. It satisfies the same property
\eqref{Hecke eigenfunction1}.

To summarize, the unramified global Langlands correspondence in the
classical setting may be viewed as a correspondence between unramified
homomorphisms $\sigma: W_F \to {}^L G$ and Hecke eigenfunctions on
\eqref{double quotient} (some irreducibility condition on $\sigma$
needs to be added to make this more precise, but we will ignore
this).

What should be the geometric analogue of this correspondence, when $X$
is a complex algebraic curve?

As explained in \secref{fund grp}, the geometric analogue of an
unramified homomorphism $W_F \to {}^L G$ is a homomorphism $\pi_1(X)
\to {}^L G$, or equivalently, since $X$ is assumed to be compact, a
holomorphic $^L G$-bundle on $X$ with a holomorphic connection (it
automatically gives rise to a flat connection). The global geometric
Langlands correspondence should therefore associate to a flat
holomorphic $^L G$-bundle on $X$ a geometric object on a geometric
version of the double quotient \eqref{double quotient}. As we argued
in \secref{from functions to sheaves}, this should be a ${\mc
D}$-module on an algebraic variety whose set of points is
\eqref{double quotient}.

Now, it is known that \eqref{double quotient} is in bijection with the
set of isomorphism classes of $G$-bundles on $X$. This key result is
due to A. Weil (see, e.g., \cite{F:rev}, Sect. 3.2). This suggests
that \eqref{double quotient} is the set of points of the moduli space
of $G$-bundles on $X$. Unfortunately, in general this is not an
algebraic variety, but an algebraic stack, which locally looks like
the quotient of an algebraic variety by an action of an algebraic
group. We denote it by $\Bun_G$. The theory of ${\mc D}$-modules has
been developed in the setting of algebraic stacks like $\Bun_G$ in
\cite{BD}, and so we can use it for our purposes. Thus, we would like
to attach to a flat holomorphic $^L G$-bundle $E$ on $X$ a ${\mc
D}$-module $\Aut_E$ on $\Bun_G$. This ${\mc D}$-module should satisfy
an analogue of the Hecke eigenfunction condition, which makes it into
a {\bf Hecke eigensheaf} with eigenvalue $E$. This notion is spelled
out in \cite{F:rev}, Sect. 6.1 (following \cite{BD}), where we refer
the reader for details.

This brings us to the following question:

\medskip

\noindent {\em How to relate this global correspondence to the local
geometric Langlands correspondence discussed above?}

\medskip

As we have already seen in \secref{unram first}, the key element in
answering this question is a {\bf localization functor}
$\Delta_{\ka_c,x}$ from $(\ghat_{\ka_c,x},G(\OO_x))$-modules to
(twisted) ${\mc D}$-modules on $\Bun_G$. In \secref{unram first} we
have applied this functor to the object $\BV_0(\chi_x)$ of the
Harish-Chandra category $\gmodx_{\chi_x}^{G(\OO_x)}$, where $\chi_x
\in \on{Op}_{^L G}(D_x)$. For an oper $\chi_x$ which extends from
$D_x$ to the entire curve $X$ we have obtained this way the Hecke
eigensheaf associated to the underlying $^L G$-local system (see
\thmref{bd}).

For a $^L G$-local system $E=(\F,\nabla)$ on $X$ which does not admit
the structure of a regular oper on $X$, the above construction may be
modified as follows (see the discussion in \cite{F:rev}, Sect. 9.6,
based on an unpublished work of Beilinson and Drinfeld). In this case
one can choose an $^L B$-reduction $\F_{^L B}$ satisfying the oper
condition away from a finite set of points $y_1,\ldots,y_n$ and such
that the restriction $\chi_{y_i}$ of the corresponding oper $\chi$ on
$X \bs \{ y_1,\ldots,y_n \}$ to $D_{y_i}^\times$ belongs to
$\on{Op}_{^L G}^{\la_i}(D_{y_i}) \subset \on{Op}_{^L
G}(D_{y_i}^\times)$ for some $\la_i \in P^+$. Then one can construct a
Hecke eigensheaf corresponding to $E$ by applying a multi-point
version of the localization functor to the tensor product of the
quotients $\BV_{\la_i}(\chi_{y_i})$ of the Weyl modules
$\BV_{\la_i,y_i}$ (see \cite{F:rev}, Sect. 9.6).

The main lesson of this construction is that in the geometric setting
the localization functor gives us a powerful tool for converting local
Langlands categories, such as $\gmodx^{G(\OO_x)}_{\chi_x}$, into
global categories of Hecke eigensheaves. The category
$\gmodx^{G(\OO_x)}_{\chi_x}$ turns out to be very simple: it has a
unique irreducible object, $\BV_0(\chi_x)$. That is why it is
sufficient to consider its image under the localization functor, which
turns out to be the desired Hecke eigensheaf $\Aut_{E_{\chi}}$. For
general opers, with ramification, the corresponding local categories
are more complicated, as we have seen above, and so are the
corresponding categories of Hecke eigensheaves. We will consider
examples of these categories in the next section.

\subsection{Classical Langlands correspondence with ramification}

Let us first consider ramified global Langlands correspondence in the
classical setting. Suppose that we are given a homomorphism $\sigma:
W_F \to {}^L G$ that is ramified at finitely many points
$y_1,\ldots,y_n$ of $X$. Then we expect that to such $\sigma$
corresponds an automorphic representation $\bigotimes'_{x \in X}
\pi_x$ (more precisely, an $L$-packet of representations). Here $\pi_x$
is still unramified for all $x \in X \bs \{ y_1,\ldots,y_n \}$, but is
{\em ramified} at $y_1,\ldots,y_n$, i.e., the space of
$G(\OO_{y_i})$-invariant vectors in $\pi_{y_i}$ is zero. In
particular, consider the special case when each $\sigma_{y_i}:
W_{F_{y_i}} \to {}^L G$ is tamely ramified (see \secref{tame reps} for
the definition). Then, according to the results presented in
\secref{tame reps}, the corresponding $L$-packet of representations
of $G(F_{y_i})$ contains an irreducible representation $\pi_{y_i}$
with non-zero invariant vectors with respect to the Iwahori subgroup
$I_{y_i}$. Let us choose such a representation for each point $y_i$.

Consider the subspace
\begin{equation}    \label{subsp}
\bigotimes_{i=1}^n \pi_{y_i}^{I_{y_i}} \otimes \bigotimes_{x \neq y_i}
v_x \subset \bigotimes_{x \in X}{}' \pi_x,
\end{equation}
where $v_x$ is a $G(\OO_x)$-vector in $\pi_x, x \neq y_i,
i=1,\ldots,n$. Then, because $\bigotimes'_{x \in X} \pi_x$ is realized
in the space of functions on $G(F)\bs G(\AD_F)$, we obtain that
the subspace \eqref{subsp} is realized in the space of functions on
the double quotient
\begin{equation}    \label{double quotient1}
G(F)\bs G(\AD_F)/\prod_{i=1}^n I_{y_i} \times \prod_{x \neq y_i}
G(\OO_x).
\end{equation}

The spherical Hecke algebras $H(G(F_x),G(\OO_x)), x \neq y_i$, act on
the subspace \eqref{subsp}, and all elements of \eqref{subsp} are
eigenfunctions of these algebras (they satisfy formula \eqref{Hecke
eigenfunction1}). At the points $y_i$ we have, instead of the action
of the commutative spherical Hecke algebra
$H(G(F_{y_i}),G(\OO_{y_i})$, the action of the non-commutative affine
Hecke algebra $H(G(F_{y_i}),I_{y_i})$. Thus, we obtain a subspace of
the space of functions on \eqref{double quotient1}, which consists of
Hecke eigenfunctions with respect to the spherical Hecke algebras
$H(G(F_x),G(\OO_x)), x \neq y_i$, and which realize a module over
$\bigotimes_{i=1}^n H(G(F_{y_i}),I_{y_i})$ (which is irreducible,
since each $\pi_{y_i}$ is irreducible).

This subspace encapsulates the automorphic representation
$\bigotimes'_{x \in X} \pi_x$ the way the automorphic function $f_\pi$
encapsulates an unramified automorphic representation. The difference
is that in the unramified case the function $f_\pi$ spans the
one-dimensional space of invariants of the maximal compact subgroup
$G(\OO)$ in $\bigotimes'_{x \in X} \pi_x$, whereas in the tamely
ramified case the subspace \eqref{subsp} is in general a
multi-dimensional vector space.

\subsection{Geometric Langlands correspondence in the tamely ramified
  case}    \label{tram}

Now let us see how this plays out in the geometric setting. As we
discussed before, the analogue of a homomorphism $\sigma: W_F \to {}^L
G$ tamely ramified at the points $y_1,\ldots,y_n \in X$ is now a local
system $E = (\F,\nabla)$, where $\F$ a $^L G$-bundle $\F$ on $X$ with
a connection $\nabla$ that has regular singularities at
$y_1,\ldots,y_n$ and unipotent monodromies around these points. We
will call such a local system {\bf tamely ramified} at
$y_1,\ldots,y_n$. What should the global geometric Langlands
correspondence attach to such a local system? It is clear that we need
to find a geometric object replacing the finite-dimensional vector
space \eqref{subsp} realized in the space of functions on
\eqref{double quotient1}.

Just as \eqref{double quotient} is the set of points of the moduli
stack $\Bun_G$ of $G$-bundles, the double quotient \eqref{double
quotient1} is the set of points of the moduli stack $\Bun_{G,(y_i)}$
of $G$-bundles on $X$ with the {\bf parabolic structures} at $y_i,
i=1,\ldots,n$. By definition, a parabolic structure of a $G$-bundle
${\mc P}$ at $y \in X$ is a reduction of the fiber ${\mc P}_y$ of
${\mc P}$ at $y$ to a Borel subgroup $B \subset G$. Therefore, as
before, we obtain that a proper replacement for \eqref{subsp} is a
category of ${\mc D}$-modules on $\Bun_{G,(y_i)}$. As in the
unramified case, we have the notion of a Hecke eigensheaf on
$\Bun_{G,(y_i)}$. But because the Hecke functors are now defined using
the Hecke correspondences over $X \bs \{ y_1,\ldots,y_n \}$ (and not
over $X$ as before), an ``eigenvalue'' of the Hecke operators is now
an $^L G$-local system on $X \bs \{ y_1,\ldots,y_n \}$ (rather than on
$X$). Thus, we obtain that the global geometric Langlands
correspondence now should assign to a $^L G$-local system $E$ on $X$,
tamely ramified at the points $y_1,\ldots,y_n$, a {\bf category} ${\mc
Aut}_E$ of ${\mc D}$-modules on $\Bun_{G,(y_i)}$ with the eigenvalue
$E|_{X \bs \{ y_1,\ldots,y_n \}}$,
$$
E \mapsto {\mc Aut}_E.
$$

We now construct these categories using a generalization of the
localization functor we used in the unramified case (see
\cite{FG:local}). For the sake of notational simplicity, let us assume
that our $^L G$-local system $E = (\F,\nabla)$ is tamely ramified at a
single point $y \in X$. Suppose that this local system on $X \bs y$
admits the structure of a $^L G$-oper $\chi = (\F,\nabla,\F_{^L B})$
whose restriction $\chi_y$ to the punctured disc $D_y^\times$ belongs
to the subspace $\nOp_{^L G}(D_y)$ of nilpotent $^L G$-opers.

For a simple Lie group $G$, the moduli stack $\Bun_{G,y}$ has a
realization analogous to \eqref{global uniform}:
$$
\Bun_{G,y} \simeq G_{\out} \bs G(\K_y)/I_y.
$$
Let ${\mc D}_{\ka_c,I_y}$ be the sheaf of twisted differential
operators on $\Bun_{G,y}$ acting on the line bundle corresponding to
the critical level (it is the pull-back of the square root of the
canonical line bundle $K^{1/2}$ on $\Bun_G$ under the natural
projection $\Bun_{G,y} \to \Bun_G$). Applying the formalism of the
previous section, we obtain a localization functor
$$
\Delta_{\ka_c,I_y}: \gmody^{I_y} \to {\mc D}_{\ka_c,I_y}\mod.
$$
However, in order to make contact with the results obtained above we
also consider the larger category $\gmody^{I^0_y}$ of
$I^0_y$-equiva\-riant modules, where $I^0_y = [I_y,I_y]$.

Set
$$
\Bun'_{G,y} = G_{\out} \bs G(\K_y)/I^0_y,
$$
and let ${\mc D}_{\ka_c,I^0_y}$ be the sheaf of twisted differential
operators on $\Bun'_{G,y}$ acting on the pull-back of the line bundle
$K^{1/2}$ on $\Bun_G$. Applying the general formalism, we
obtain a localization functor
\begin{equation}    \label{I0}
\Delta_{\ka_c,I^0_y}: \gmody^{I_y^0} \to {\mc
  D}_{\ka_c,I^0_y}\mod.
\end{equation}
We note that a version of the categorical affine Hecke algebra ${\mc
H}(G(\K_y),I_y)$ discussed in \secref{evidence} naturally acts on the
derived categories of the above categories, and the functors
$\Delta_{\ka_c,I_y}$ and $\Delta_{\ka_c,I^0_y}$ intertwine these
actions. Equivalently, one can say that this functor intertwines the
corresponding actions of the affine braid group associated to $^L G$
on the two categories (as in \cite{Bez}).

We now restrict the functors $\Delta_{\ka_c,I_y}$ and
$\Delta_{\ka_c,I^0_y}$ to the subcategories $\gmody^{I_y}_{\chi_y}$
and $\gmody^{I_y^0}_{\chi_y}$, respectively. By using the same
argument as in \cite{BD}, we obtain the following analogue of
\thmref{bd}.

\begin{thm}    \label{bd1}
Fix $\chi_y \in \nOp_{^L G}(D_y)$ and let $M$ be an object of the
category $\gmody^{I_y}_{\chi_y}$
(resp., $\gmody^{I^0_y}_{\chi_y}$). Then

{\em (1)} $\Delta_{\ka_c,I_y}(M) = 0$ (resp.,
$\Delta_{\ka_c,I^0_y}(M) = 0$) unless $\chi_y$ is the restriction of
a regular oper $\chi = (\F,\nabla,\F_{^L B})$ on $X \bs y$ to
$D_y^\times$.

{\em (2)} In that case $\Delta_{\ka_c,y}(M)$ (resp.,
$\Delta_{\ka_c,I^0_y}(M)$) is a Hecke eigensheaf with the eigenvalue
$E_\chi = (\F,\nabla)$.
\end{thm}

Thus, we obtain that if $\chi_y = \chi|_{D_y^\times}$, then the image
of any object of $\gmody^{I_y}_{\chi_y}$ under the functor
$\Delta_{\ka_c,I_y}$ belongs to the category ${\mc
Aut}^{I_y}_{E_\chi}$ of Hecke eigensheaves on $\Bun_{G,y}$. Now
consider the restriction of the functor $\Delta_{\ka_c,I^0_y}$ to
$\gmody^{I^0_y}_{\chi_y}$. As discussed in \secref{conj descr}, the
category $\gmody^{I^0_y}_{\chi_y}$ coincides with the corresponding
category $\gmody^{I_y,m}_{\chi_y}$ of $I_y$-monodromic
modules. Therefore the image of any object of
$\gmody^{I^0_y}_{\chi_y}$ under the functor $\Delta_{\ka_c,I^0_y}$
belongs to the subcategory ${\mc D}^{m}_{\ka_c,I^0_y}\mod$ of ${\mc
D}_{\ka_c,I^0_y}\mod$ whose objects admit an increasing filtration
such that the consecutive quotients are pull-backs of ${\mc
D}_{\ka_c,I_y}$-modules from $\Bun_{G,y}$. Such ${\mc
D}_{\ka_c,I^0_y}$--modules are called {\bf monodromic}.

Let ${\mc Aut}^{I_y,m}_{E_\chi}$ be the subcategory of ${\mc
D}^{m}_{\ka_c,I^0_y}\mod$ whose objects are Hecke eigensheaves with
eigenvalue $E_\chi$.

Thus, we obtain functors
\begin{equation}    \label{Delta y1}
\Delta_{\ka_c,I_y}: \gmody^{I_y}_{\chi_y} \to {\mc Aut}^{I_y}_{E_\chi},
\qquad \Delta_{\ka_c,I^0_y}: \gmody^{I^0_y}_{\chi_y} \to {\mc
Aut}^{I_y,m}_{E_\chi}.
\end{equation}
It is tempting to conjecture (see \cite{FG:local}) that these functors
are equivalences of categories, at least for generic $\chi$. Suppose
that this is true. Then we may identify the {\em global} categories
${\mc Aut}^{I_y}_{E_\chi}$ and ${\mc Aut}^{I_y,m}_{E_\chi}$ of Hecke
eigensheaves on $\Bun_{G,I_y}$ and $\Bun'_{G,I^0_y}$ with the {\em
local} categories $\gmody^{I_y}_{\chi_y}$ and
$\gmody^{I^0_y}_{\chi_y}$, respectively.  Therefore we can use our
results and conjectures on the local Langlands categories, such as
$\gmody^{I^0_y}_{\chi_y}$, to describe the global categories of Hecke
eigensheaves on the moduli stacks of $G$-bundles on $X$ with parabolic
structures.

We have the following conjectural description of the derived category
of $I^0_y$-equivariant modules, $D^b(\gmody_{\chi_y})^{I^0_y}$ (see
formula \eqref{nilp fiber}):
\begin{equation}    \label{nilp fiber2}
D^b(\gmody_{\chi_y})^{I^0_y} \simeq
D^b(\QCoh(\wt{\on{Sp}}^{\on{DG}}_{\on{Res}(\chi_y)})).
\end{equation}

The corresponding $I_y$-equivariant version is
\begin{equation}    \label{nilp fiber1}
D^b(\gmody_{\chi_y})^{I_y} \simeq
D^b(\QCoh(\on{Sp}^{\on{DG}}_{\on{Res}(\chi_y)})),
\end{equation}
where we replace the non-reduced DG Springer fiber by the reduced one:
it is defined as the DG fiber of the Springer resolution $\wt{\mc N}
\to {\mc N}$ (see formula \eqref{O Sp}).

If the functors \eqref{Delta y1} are equivalences, then by combining
them with \eqref{nilp fiber2} and \eqref{nilp fiber1}, we obtain the
following conjectural equivalences of categories:
\begin{equation}    \label{grand equiv}
D^b({\mc Aut}^{I_y}_{E_\chi}) \simeq
D^b(\QCoh({\on{Sp}}^{\on{DG}}_{\on{Res}(\chi_y)})), \qquad D^b({\mc
Aut}^{I_y,m}_{E_\chi}) \simeq
D^b(\QCoh(\wt{\on{Sp}}^{\on{DG}}_{\on{Res}(\chi_y)})).
\end{equation}
In other words, the derived category of a global Langlands category
(monodromic or not) corresponding to a local system tamely ramified at
$y \in X$ is equivalent to the derived category of quasicoherent
sheaves on the DG Springer fiber of its residue at $y$ (non-reduced or
reduced).

Again, these equivalences are supposed to intertwine the natural
actions on the above categories of the categorical affine Hecke
algebra ${\mc H}(G(\K_y),I_y)$ (or, equivalently, the affine braid
group associated to $^L G$).

The categories appearing in \eqref{grand equiv} actually make sense
for an arbitrary $^L G$-local system $E$ on $X$ tamely ramified at
$y$. It is therefore tempting to conjecture that these equivalences
still hold in general:
\begin{equation}    \label{grand equiv1}
D^b({\mc Aut}^{I_y}_{E}) \simeq
D^b(\QCoh({\on{Sp}}^{\on{DG}}_{\on{Res}(E)})), \qquad D^b({\mc
Aut}^{I_y,m}_{E}) \simeq
D^b(\QCoh(\wt{\on{Sp}}^{\on{DG}}_{\on{Res}(E)})).
\end{equation}
The corresponding localization functors are constructed as follows: we
represent a general local system $E$ on $X$ with tame ramification at
$y$ by an oper $\chi$ on the complement of finitely many points
$y_1,\ldots,y_n$, whose restriction to $D_{y_i}^\times$ belongs to
$\on{Op}_{^L G}^{\la_i}(D_{y_i}) \subset \on{Op}_{^L
G}(D_{y_i}^\times)$ for some $\la_i \in P^+$. Then, in the same way as
in the unramified case, we construct localization functors from
$\gmody_{\chi_y}^{I_y}$ to ${\mc Aut}^{I_y}_{E}$ and from
$\gmody_{\chi_y}^{I^0_y}$ to ${\mc Aut}^{I_y,m}_{E}$ (here, as
before, $\chi_y = \chi|_{D_y^\times}$), and this leads us to the
conjectural equivalences \eqref{grand equiv1}.

The equivalences \eqref{grand equiv1} also have family versions in
which we allow $E$ to vary. It is analogous to the family version
\eqref{equiv FG} of the local equivalences. As in the local case, in
a family version we can avoid using DG schemes.

The above construction may be generalized to allow local systems
tamely ramified at finitely many points $y_1,\ldots,y_n$. The
corresponding Hecke eigensheaves are then $\D$-modules on the moduli
stack of $G$-bundles on $X$ with parabolic structures at
$y_1,\ldots,y_n$. Non-trivial examples of these Hecke eigensheaves
arise already in genus zero. These sheaves were constructed explicitly
in \cite{F:icmp} (see also \cite{F:faro,F:flag}), and they are closely
related to the Gaudin integrable system.

\subsection{Connections with regular singularities}

So far we have only considered the categories of
$\ghat_{\ka_c}$-modules corresponding to $^L G$-opers on $X$ which are
regular everywhere except at a point $y \in X$ (or perhaps, at several
points) and whose restriction to $D_y^\times$ is a nilpotent oper
$\chi_y$ in $\nOp_{^L G}(D_y)$. In other words, $\chi_y$ is an oper
with regular singularity at $y$ with residue $\varpi(-\rho)$ (where
$\varpi: \h^* \to \h^*/W$), see the definition in \secref{admitting
I}. However, we can easily generalize the localization functor to the
categories of $\ghat_{\ka_c}$-modules corresponding to $^L G$-opers
which have regular singularity at $y$ with {\em arbitrary}
residue. The monodromy of such an oper could lie in an arbitrary
conjugacy class in $^L G$, not necessarily a unipotent one.

So suppose we are given an oper $\chi \in \on{Op}^{\on{RS}}_{^L
G}(D)_{\varpi(-\la-\rho)}$ with regular singularity and residue
$\varpi(-\la-\rho)$, where $\la \in \h^*$. In this case the monodromy
of this oper around $y$ is is such that its semi-simple part is
conjugate to
$$
M_{\on{ss}} = \exp(2 \pi i (\la+\rho)) = \exp(2 \pi i \la).
$$
We then have the category $\gmod^{I^0}_\chi$ of $I^0$-equivariant
$\ghat_{\ka_c}$-modules with central character $\chi$. The case of
$\la=0$ is an ``extremal'' case when the category $\gmod^{I^0}_\chi$
is most complicated. On the other ``extreme'' is the case of generic
opers $\chi$, corresponding to a generic $\la$. In this case one can
show that the category $\gmod^{I^0}_\chi$ is quite simple: it contains
irreducible objects $\BM_{w(\la+\rho)-\rho}(\chi)$ labeled by elements
$w$ of the Weyl group $W$ of $\g$, and each object of
$\gmod^{I^0}_\chi$ is a direct sum of these irreducible modules. Here
$\BM_{w(\la+\rho)-\rho}(\chi)$ is the quotient of the Verma module
$$
\BM_{w(\la+\rho)-\rho} = \on{Ind}_{\wh\bb_+ \oplus \C {\mb
    1}}^{\ghat_{\ka_c}} \C_{w(\la+\rho)-\rho}, \qquad w \in W,
$$
where $\wh\bb_+ = (\bb_+ \otimes 1) \oplus (\g \otimes t\C[[t]])$, by
the central character corresponding to $\chi$.

For other values of $\la$ the structure of $\gmod^{I^0}_\chi$
is somewhere in-between these two extreme cases.

Recall that we have a localization functor \eqref{I0},
$$
\Delta^\la_{\ka_c,I^0_y}: \gmody^{I^0_y} \to {\mc
D}_{\ka_c,I^0_y}\mod,
$$
from $\gmody^{I^0_y}_{\chi_y}$ to a category of ${\mc D}$-modules on
$\Bun'_{G,I_y}$ twisted by the pull-back of the line bundle $K^{1/2}$
on $\Bun_G$. We now restrict this functor to the subcategory
$\gmody^{I^0_y}_{\chi_y}$ where $\chi_y$ is a $^L G$-oper on $D_y$
with regular singularity at $y$ and residue $\varpi(-\la-\rho)$.

Consider first the case when $\la \in \h^*$ is generic. Suppose that
$\chi_y$ extends to a regular oper $\chi$ on $X \bs y$. One then shows
in the same way as in \thmref{bd1} that for any object $M$ of
$\gmody^{I^0_y}_{\chi_y}$ the corresponding ${\mc
D}_{\ka_c,I^0_y}$-module $\Delta_{\ka_c,I^0_y}(M)$ is a Hecke
eigensheaf with eigenvalue $E_\chi$, which is the $^L G$-local system
on $X$ with regular singularity at $y$ underlying $\chi$ (if $\chi_y$
cannot be extended to $X \bs y$, then $\Delta^\la_{\ka_c,I_y}(M)=0$,
as before). Therefore we obtain a functor
$$
\Delta_{\ka_c,I^0_y}: \gmody^{I^0_y}_{\chi_y} \to {\mc
  Aut}^{I^0_y}_{E_\chi},
$$
where ${\mc Aut}^{I^0_y}_{E_\chi}$ is the category of Hecke eigensheaves
on $\Bun'_{G,I_y}$ with eigenvalue $E_\chi$.

Since we have assumed that the residue of the oper $\chi_y$ is
generic, the monodromy of $E_\chi$ around $y$ belongs to a regular
semi-simple conjugacy class of $^L G$ containing $\exp(2 \pi i
\la)$. In this case the category $\gmody^{I^0_y}_{\chi_y}$ is
particularly simple, as we have discussed above: there are $|W|$
irreducible objects, and there are no extensions between them. We
expect that the functor $\Delta_{\ka_c,I^0_y}$ sets up an equivalence
between $\gmody^{I^0_y}_{\chi_y}$ and ${\mc Aut}^{I^0_y}_{E_\chi}$.

We can formulate this more neatly as follows. For $M \in {}^L G$ let
${\mc B}_M$ be the variety of Borel subgroups containing $M$. Observe
that if $M$ is regular semi-simple, then ${\mc B}_M$ is a set of
points which is in bijection with $W$. Therefore our conjecture is
that ${\mc Aut}^{I^0_y}_{E_\chi}$ is equivalent to the category
$\QCoh({\mc B}_M)$ of quasicoherent sheaves on ${\mc B}_M$, where $M$
is a representative of the conjugacy class of the monodromy of
$E_\chi$.

Consider now an arbitrary $^L G$-local system $E$ on $X$ with regular
singularity at $y \in X$ whose monodromy around $y$ is regular
semi-simple. It is then tempting to conjecture that, at least if $E$
is generic, this category has the same structure as in the case when
$E$ has the structure of an oper, i.e., it is equivalent to the
category $\QCoh({\mc B}_M)$, where $M$ is a representative of the
conjugacy class of the monodromy of $E$ around $y$.

On the other hand, if the monodromy around $y$ is unipotent, then
${\mc B}_M$ is nothing but the Springer fiber $\on{Sp}_{u}$, where $M
= \exp(2\pi i u)$. The corresponding category ${\mc Aut}^{I^0_y}_E$
was discussed in \secref{tram} (we expect that it coincides with ${\mc
Aut}^{I_y,m}_E$). Thus, we see that in both ``extreme'' cases:
unipotent monodromy and regular semi-simple monodromy, our conjectures
identify the derived category of ${\mc Aut}^{I^0_y}_E$ with the
derived category of the category $\QCoh({\mc B}_M)$ (where ${\mc B}_M$
should be viewed as a DG scheme $\wt{\on{Sp}}^{\on{DG}}_u$ in the
unipotent case). One is then led to conjecture, most ambitiously, that
for {\em any} $^L G$-local system $E$ on $X$ with regular singularity
at $y \in X$ the derived category of ${\mc Aut}^{I^0_y}_E$ is
equivalent to the derived category of quasicoherent sheaves on a
suitable DG version of the scheme ${\mc B}_M$, where $M$ is a
representative of the conjugacy class of the monodromy of $E$ around
$y$:
$$
D^b({\mc Aut}^{I^0_y}_E) \simeq D^b(\QCoh({\mc B}^{\on{DG}}_M)).
$$
This has an obvious generalization to the case of multiple
ramification points, where on the right hand side we take the
Cartesian product of the varieties ${\mc B}^{\on{DG}}_{M_i}$
corresponding to the monodromies. Thus, we obtain a conjectural
realization of the categories of Hecke eigensheaves, whose eigenvalues
are local systems with regular singularities, in terms of categories of
quasicoherent sheaves.

It is useful to note that the Hecke eigensheaves on $\Bun'_{G,I_y}$
obtained above via the localization functors may be viewed as
pull-backs of twisted ${\mc D}$-modules on $\Bun_{G,I_y}$ (or, more
generally, extensions of such pull-backs).

More precisely, for each $\la \in \h^*$ we have the sheaf of twisted
differential operators acting on a ``line bundle'' $\wt{\Ll}_\la$ on
$\Bun_{G,y}$. If $\la$ were an integral weight, this would be an
actual line bundle, which is constructed as follows: note that the map
$p: \Bun_{G,I_y} \to \Bun_G$, corresponding to forgetting the
parabolic structure, is a fibration with the fibers isomorphic to the
flag manifold $G/B$. For each integral weight $\la$ we have the
$G$-equivariant line bundle $\ell_\la = G \underset{B}\times \C_\la$
on $G/B$. The line bundle $\Ll_\la$ on $\Bun_{G,I_y}$ is defined in
such a way that its restriction to each fiber of the projection $p$ is
isomorphic to $\ell_\la$. We then set $\wt{\Ll}_\la = {\mc L}_\la
\otimes p^*(K^{1/2})$, where $K^{1/2}$ is the square root of the
canonical line bundle on $\Bun_G$ corresponding to the critical level.
Now, it is well-known (see, e.g., \cite{BB}) that even though the line
bundle $\wt\Ll_\la$ does not exist if $\la$ is not an integral weight,
the corresponding sheaf of $\wt\Ll_\la$-twisted differential operators
on $\Bun_{G,I_y}$ is still well-defined. We denote it by ${\mc
D}^\la_{\ka_c,I_y}$.

Observe that we have an action of the group $H = I/I^0$ on
$\Bun'_{G,y}$. There is an equivalence between the category ${\mc
D}^\la_{\ka_c,I_y}\mod$ and the category of weakly $H$-equivariant
${\mc D}_{\ka_c,I^0_y}$-modules on $\Bun'_{G,y}$ on which $\h$ acts
via the character $\la: \h \to \C$ (see \cite{BL1}). If ${\mc F}$ is
an object of ${\mc D}^\la_{\ka_c,I_y}\mod$, then the corresponding
weakly $H$-equivariant ${\mc D}_{\ka_c,I^0_y}$-module on $\Bun'_{G,y}$
is $\pi^*({\mc F})$, where $\pi: \Bun'_{G,y} \to \Bun_{G,I_y}$.

Now, it is easy to see that the ${\mc D}_{\ka_c,I^0_y}$-modules on
$\Bun'_{G,y}$ obtained by applying the localization functor
$\Delta_{\ka_c,I^0_y}$ to objects of $\gmody^{I^0_y}_{\chi_y}$ are
always weakly $H$-equivariant. Consider, for example, the case when
$\chi_y$ is an oper with regular singularity at $y$ with residue
$\varpi(-\la-\rho)$, where $\la$ is a regular element of $\h^*$. Then
its monodromy is $M = \exp(2 \pi i \la)$. The corresponding category
$\gmody^{I^0_y}_{\chi_y}$ has objects $\BM_{w(\la+\rho)-\rho}(\chi_y),
w \in W$, that we introduced above. The Cartan subalgebra $\h$ of
$\ghat_{\ka_c,y}$ acts on $\BM_{w(\la+\rho)-\rho}(\chi_y)$ semi-simply
with the eigenvalues given by the weights of the form
$w(\la+\rho)-\rho+\mu$, where $\mu$ is an integral weight. In other
words,
$$
\BM_{w(\la+\rho)-\rho}(\chi_y) \otimes \C_{-w(\la+\rho)+\rho}
$$
is $I_y$-equivariant. Therefore we find that the ${\mc
D}_{\ka_c,I^0_y}$-module
$\Delta_{\ka_c,I^0_y}(\BM_{w(\la+\rho)-\rho}(\chi_y))$ is weakly
$H$-equivariant, and the corresponding action of $\h$ is given by
$w(\la+\rho)-\rho: \h \to \C$. Thus,
$\Delta_{\ka_c,I^0_y}(\BM_{w(\la+\rho)-\rho}(\chi_y))$ is the pull-back
of a ${\mc D}^{w(\la+\rho)-\rho}_{\ka_c,I_y}$-module on
$\Bun_{G,y}$. This ${\mc D}^{w(\la+\rho)-\rho}_{\ka_c,I_y}$-module is
a Hecke eigensheaf with eigenvalue $E_\chi$ provided that $\chi_y =
\chi|_{D^\times_y}$, where $\chi$ is a regular oper on $X \bs y$.

Thus, for a given generic oper $\chi_y$ we have $|W|$ different Hecke
eigensheaves
$$
\Delta_{\ka_c,I^0_y}(\BM_{w(\la+\rho)-\rho}(\chi_y)), \qquad w
\in W,
$$
on $\Bun'_{G,y}$. However, each of them is the pull-back of a twisted
${\mc D}$-module on $\Bun_{G,y}$ corresponding to a particular twist:
namely, by a ``line bundle'' $\wt{\Ll}_{w(\la+\rho)-\rho}$. Since we
have assumed that $\la$ is generic, all of these twists are different;
note also that if $\mu = w(\la+\rho)-\rho$, then $\exp(2\pi i \mu)$ is
in the conjugacy class of the monodromy $\exp(2 \pi i \la)$. It is
therefore natural to conjecture that there is a unique irreducible
Hecke eigensheaf on $\Bun_{G,y}$ with eigenvalue $E_\chi$, which is a
twisted ${\mc D}$-module with the twisting given by
$\wt{\Ll}_{w(\la+\rho)-\rho}$.

More generally, suppose that $E$ is a local system on $X$ with regular
singularity at $y$ and regular semi-simple monodromy.  Let us choose a
representative $M$ of the monodromy which belongs to the Cartan
subgroup $^L H \subset {}^L G$. Choose $\mu \in \h^* = {}^L \h$ to be
such that $M = \exp(2\pi i \mu)$. Note that there are exactly $|W|$
such choices up to a shift by an integral weight $\nu$. Let ${\mc
Aut}^{I_y,\mu}_E$ be the category of Hecke eigensheaves with
eigenvalue $E$ in the category of twisted ${\mc D}$-modules on
$\Bun_{G,I_y}$ with the twisting given by $\wt{\Ll}_{\mu}$. Then we
expect that for generic $E$ of this type the category ${\mc
Aut}^{I_y,\mu}_E$ has a unique irreducible object. Its pull-back to
$\Bun'_{G,y}$ is one of the $|W|$ irreducible objects of ${\mc
Aut}^{I^0_y}_E$.

Note that tensoring with the line bundle ${\mc L}_\nu$, where $\nu$ is
an integral weight, we identify the categories ${\mc Aut}^{I_y,\mu}_E$
and ${\mc Aut}^{I_y,\mu'}_E$ if $\mu'=\mu+\nu$, which accounts for the
ambiguity in the choice of $\mu$ with respect to shifts by integral
weights.

Similarly, one can describe the Hecke eigensheaves on $\Bun'_{G,y}$
obtained by applying the functor $\Delta_{\ka_c,I^0_y}$ to the
categories $\gmody^{I^0_y}_{\chi_y}$ for other opers $\chi_y$ in terms
of twisted ${\mc D}$-modules on $\Bun_{G,y}$. In the opposite extreme
case, when $\chi_y$ is a nilpotent oper (and so its monodromy is
unipotent), this is explained in \secref{tram}. In this case we may
choose to consider monodromic ${\mc D}$-modules; note that this is not
necessary if the monodromy is regular semi-simple, because in this
case there are no non-trivial extensions between the corresponding
${\mc D}$-modules.

For a general local system $E$ on $X$ with regular singularity and
monodromy conjugate to $M \in {}^L G$, let $M_{\on{ss}}$ be the
semi-simplification of $M$. Let us choose $\mu \in \h^* = {}^L \h$
such that $M_{\on{ss}} = \exp(2\pi i \mu)$. For each such $\mu$ we
have a non-trivial category ${\mc Aut}^{I_y,\mu}_E$ of Hecke
eigensheaves with eigenvalue $E$, which are twisted ${\mc D}$-modules
on $\Bun_{G,I_y}$ with the twisting given by $\wt{\Ll}_{\mu}$. We
expect that this category may be described as a category of
quasicoherent sheaves on some variety of Borel subalgebras (at least,
for generic $E$ with the monodromy $M$).

Finally, it is natural to ask whether these equivalences for
individual local systems may be combined into a family version
encompassing all of them. It is instructive to view the global
geometric Langlands correspondence in the unramified case as a kind of
non-abelian Fourier-Mukai transform relating the (derived) category of
${\mc D}$-modules on $\Bun_G$ and the (derived) category of
quasicoherent sheaves on $\on{Loc}_{^L G}(X)$, the stack of $^L
G$-local systems on the curve $X$. Under this correspondence, the
skyscraper sheaf supported at a $^L G$-local system $E$ is supposed to
go to the Hecke eigensheaf $\Aut_E$ on $\Bun_G$. Thus, one may think
of $\on{Loc}_{^L G}(X)$ as a parameter space of a ``spectral
decomposition'' of the derived category of ${\mc D}$-modules on
$\Bun_G$ (see, e.g., \cite{F:rev}, Sect. 6.2, for more details).

The above results and conjectures suggest that one may also view the
geometric Langlands correspondence in the tamely ramified case in a
similar way. Now the role of $\on{Loc}_{^L G}(X)$ should be played by
the stack $\on{Loc}_{^L G,y}(X)$ of parabolic $^L G$-local systems
with regular singularity at $y \in X$ (or, more generally, multiple
points) and unipotent monodromy. This stack classifies triples
$(\F,\nabla,\F_{^L B,y})$, where $\F$ is a $^L G$-bundle on $X$,
$\nabla$ is a connection on $\F$ with regular singularity at $y$ and
unipotent monodromy, and $\F_{^L B,y}$ is a $^L B$-reduction of the
fiber $\F_y$ of $\F$ at $y$, which is preserved by $\nabla$. This
stack is now a candidate for a parameter space of a ``spectral
decomposition'' of the derived category of ${\mc D}$-modules on the
moduli stack $\Bun_{G,y}$ of $G$-bundles with parabolic structure at
$y$.\footnote{One may also try to extend this ``spectral
decomposition'' to the case of all connections with regular
singularities, but here the situation is more subtle, as can already
be seen in the abelian case.}

\subsection{Irregular connections}

We now generalize the above results to the case of connections with
irregular singularities. Let $\F$ be a $^L G$-bundle on $X$ with
connection $\nabla$ that is regular everywhere except for a point $y
\in X$, where it has a pole of order greater than $1$. As before, we
assume first that $(\F,\nabla)$ admits the structure of a $^L G$-oper
on $X \bs y$, which we denote by $\chi$. Let $\chi_y$ be the the
restriction of $\chi$ to $D_y^\times$. A typical example of such an
oper is a $^L G$-oper with pole of order $\leq n$ on the disc $D_y$,
which is, by definition (see \cite{BD}, Sect. 3.8.8), an $^L
N[[t]]$-conjugacy class of operators of the form
\begin{equation} \label{oper with RS1}
\nabla = \pa_t + \frac{1}{t^n} \left( p_{-1} +
{\mb v}(t) \right), \qquad {\mb v}(t) \in {}^L \bb[[t]].
\end{equation}
We denote the space of such opers by $\on{Op}_{^L G}^{\leq n}(D_y)$.

One can show that for $\chi_y \in \on{Op}_{^L G}^{\leq n}(D_y)$ the
category $\gmody_{\chi_y}^K$ is non-trivial if $K$ is the congruence
subgroup $K_{m,y} \subset G(\OO_y)$ with $m \geq n$. (We recall that
for $m>0$ we have $K_{m,y} = \exp(\g \otimes (\mm_y)^m)$, where
$\mm_y$ is the maximal ideal of $\OO_y$.) Let us take the category
$\gmody_{\chi_y}^{K_{n,y}}$ (in principle, we could take another
``compact'' subgroup $K$ instead of $K_{n,y}$). Then our general
formalism gives us a localization functor
$$
\Delta_{\ka_c,K_{n,y}}: \gmody_{\chi_y}^{K_{n,y}} \to {\mc
D}_{\ka_c,K_{n,y}}\mod,
$$
where ${\mc D}_{\ka_c,K_{n,y}}\mod$ is the category of critically
twisted\footnote{this refers to the twisting by the line
bundle on $\Bun_{G,y,n}$ obtained by pull-back of the line bundle
$K^{1/2}$ on $\Bun_G$, as before} ${\mc D}$-modules on
$$
\Bun_{G,y,n} \simeq G_{\out} \bs G(\K_y)/K_{n,y}.
$$
This is the moduli stack of $G$-bundles on $X$ with a level $n$
structure at $y \in X$ (which is a trivialization of the restriction
of the $G$-bundle to the $n$th infinitesimal neighborhood of
$y$).

In the same way as above, one shows that the ${\mc D}$-modules
obtained by applying $\Delta_{\ka_c,K_{n,y}}$ to objects of
$\gmody_{\chi_y}^{K_{n,y}}$ are Hecke eigensheaves with the eigenvalue
$E_\chi|_{X \bs y}$, where $E_\chi$ is the $^L G$-local system
underlying the oper $\chi$. Let ${\mc Aut}^{K_{n,y}}_{E_\chi}$ be the
category of these eigensheaves. Thus, we really obtain a functor
$$
\gmody_{\chi_y}^{K_{n,y}} \to {\mc Aut}^{K_{n,y}}_{E_\chi}.
$$
By analogy with the case of regular connections, we expect that
this functor is an equivalence of categories.

As before, this functor may be generalized to an arbitrary flat bundle
$E = (\F,\nabla)$, where $\nabla$ has singularity at $y$, by
representing it as an oper with mild ramification at additional points
$y_1,\ldots,y_m$ on $X$. Let $\chi_y$ be the restriction of this oper
to $D_y^\times$. Then it belongs to $\on{Op}_{^L G}^{\leq n}(D_y)$ for
some $n$, and we obtain a functor
$$
\gmody_{\chi_y}^{K_{n,y}} \to {\mc Aut}^{K_{n,y}}_{E},
$$
which we expect to be an equivalence of categories for generic
$E$. This also has an obvious multi-point generalization.

This way we obtain a conjectural description of the categories of
Hecke eigensheaves corresponding to (generic) connections on $X$ with
arbitrary singularities at finitely many points in terms of categories
of Harish-Chandra modules of critical level over $\ghat$. However, in
the case of regular singularities, we also have an alternative
description of these categories: in terms of (derived) categories of
quasicoherent sheaves on the varieties ${\mc B}^{\on{DG}}_M$. It would
be desirable to obtain such a description for irregular connections as
well.

Finally, we remark that the above construction has a kind of limiting
version where we take the infinite level structure at $y$, i.e., a
trivialization of the restriction of a $G$-bundle to the disc
$D_y$. Let $\Bun_{G,y,\infty}$ be the moduli stack of $G$-bundles on
$X$ with an infinite level structure at $y$. Then
$$
\Bun_{G,y,\infty} \simeq G_{\out} \bs G(\K_y).
$$
We now have a localization functor
$$
\gmody_{\chi_y} \to {\mc Aut}^\infty_{E},
$$
where $E$ and $\chi_y$ are as above, and ${\mc Aut}^\infty_{E}$ is the
category of Hecke eigensheaves on $\Bun_{G,y,\infty}$ with eigenvalue
$E|_{X\bs y}$. Thus, instead of the category
$\gmody_{\chi_y}^{K_{n,y}}$ of Harish-Chandra modules we now have the
category $\gmody_{\chi_y}$ of all (smooth) $\ghat_{\ka_c,y}$-modules
with fixed central character (corresponding to $\chi$).

According to our general local conjecture, this is precisely the local
Langlands category associated to the restriction of the local system
$E$ to $D_y^\times$ (equipped with an action of the loop group
$G(\K_y)$). It is natural to assume that for generic $E$ this functor
establishes an equivalence between this category and the category
${\mc Aut}^\infty_{E}$ of Hecke eigensheaves on $\Bun_{G,y,\infty}$
(also equipped with an action of the loop group $G(\K_y)$). This may
be thought of as the ultimate form of the local--to--global
compatibility in the geometric Langlands Program:
$$
\begin{CD}
E @>>> {\mc Aut}^\infty_E \\
@VVV @AAA \\
E|_{D_y^\times} @>>> \gmod_{\chi_y}.
\end{CD}
$$

\bigskip

Let us summarize: by using representation theory of affine Kac-Moody
algebras at the critical level we have constructed the local Langlands
categories corresponding to the local Langlands parameters: $^L
G$-local systems on the punctured disc. We then applied the technique
of localization functors to produce from these local categories, the
global categories of Hecke eigensheaves on the moduli stacks of
$G$-bundles on a curve $X$ with parabolic (or level) structures. These
global categories correspond to the global Langlands parameters: $^L
G$-local systems on $X$ with ramification. We have used our results
and conjectures on the structure of the local categories to
investigate these global categories.  We hope that in this way
representation theory of affine Kac-Moody algebras may one day fulfill
the dream of uncovering the mysteries of the geometric Langlands
correspondence.


\begin{thebibliography}{BeLu}

\bibitem[AB]{AB} A.~Arkhipov and R.~Bezrukavnikov,
{\em Perverse sheaves on affine flags and Langlands dual group},
Preprint math.RT/0201073.

\bibitem[ABG]{ABG} S. Arkhipov, R. Bezrukavnikov and V. Ginzburg, {\em
Quantum groups, the loop Grassmannian, and the Springer resolution},
Journal of AMS {\bf 17} (2004) 595--678.

\bibitem[AG]{AG} S.~Arkhipov and D.~Gaitsgory, {\em Another
realization of the category of modules over the small quantum group},
Adv. Math. {\bf 173} (2003) 114--143.

\bi[A]{Arthur} J. Arthur, {\em Unipotent automorphic representations:
conjectures}, Asterisque {\bf 171-172} (1989) 13--71.

\bi[BV]{BV} D.G. Babbitt and V.S. Varadarajan, {\em Formal reduction
  theory of meromorphic differential equations: a group theoretic
  view}, Pacific J. Math. {\bf 109} (1983) 1--80.

\bi[BeLa]{BL} A.~Beauville and Y.~Laszlo, {\em Un lemme de descente},
C.R. Acad.  Sci. Paris, S\'{e}r. I Math. {\bf 320} (1995) 335--340.

\bi[Bei]{Be} A. Beilinson, {\em Langlands parameters for Heisenberg
  modules}, Preprint math.QA/0204020.

\bi[BB]{BB} A. Beilinson and J. Bernstein, {\em A proof of Jantzen
conjectures}, Advances in Soviet Mathematics {\bf 16}, Part 1,
pp. 1--50, AMS, 1993.

\bi[BD1]{BD} A. Beilinson and V. Drinfeld, {\em Quantization of
Hitchin's integrable system and Hecke eigensheaves}, Preprint,
available at www.math.uchicago.edu/$\sim$arinkin

\bibitem[BD2]{BD:CHA} A. Beilinson and V. Drinfeld, {\em Chiral
algebras}, Colloq. Publ. {\bf 51}, AMS, 2004.

\bibitem[BD3]{BD:opers} A. Beilinson and V. Drinfeld, {\em Opers},
Preprint math.AG/0501398.

\bibitem[BeLu]{BL1} J. Bernstein and V. Lunts, {\em Localization for
  derived categories of $(\g,K)$-modules}, Journal of AMS {\bf 8}
  (1995) 819--856.

\bibitem[BZ]{BZ} J. Bernstein and A. Zelevinsky, {\em Induced
  representations of reductive $p$-adic groups}, I, Ann. Sci. ENS {\bf
  10} (1977) 441--472.

\bibitem[Bez1]{Bez1} R. Bezrukavnikov, {\em Perverse sheaves on affine
  flags and nilpotent cone of the Langlands dual group}, Preprint
  math.RT/0201256.

\bibitem[Bez2]{Bez} R. Bezrukavnikov, {\em Noncommutative counterparts
of the Springer resolution}, Preprint math.RT/0604445.

\bi[B1]{Borel} A. Borel, {\em Admissible representations of a
  semi-simple group over a local field with vectors fixed under an
  Iwahori subgroup}, Inv. Math. {\bf 35} (1976) 233--259.

\bi[B2]{Dmodules} A.~Borel, e.a., {\em Algebraic $D$--modules},
Academic Press, 1987.

\bibitem[CG]{Ginzburg} N. Chriss and V. Ginzburg, {\em Representation
  theory and complex geometry}, Birkh\"auser 1997.

\bibitem[CK]{CK} I. Ciocan-Fountanine and M. Kapranov, {\em Derived
  Quot schemes}, Ann. Sci. ENS {\bf 34} (2001) 403--440.

\bibitem[De1]{D} P. Deligne, {\em Equations diff\'erentielles \'a
  points singuliers r\'eguliers}, Lect. Notes in Math. {\bf 163},
  Springer, 1970.

\bibitem[De2]{Del} P. Deligne, {\em Les constantes des \'equations
fonctionnelles des fonctions L}, in {\em Modular Functions one
Variable II}, Proc. Internat. Summer School, Univ. Antwerp 1972,
Lect. Notes Math. {\bf 349}, pp. 501--597, Springer 1973.

\bi[Dr1]{Dr1} V.G. Drinfeld, {\em Langlands conjecture for $GL(2)$ over
  function field}, Proc. of Int. Congress of Math. (Helsinki, 1978),
  pp. 565--574.

\bi[Dr2]{Dr2} V.G. Drinfeld, {\em Two-dimensional $\ell$--adic
representations of the fundamental group of a curve over a finite
field and automorphic forms on $GL(2)$}, Amer. J. Math. {\bf 105}
(1983) 85--114.

\bi[Dr3]{Dr3} V.G. Drinfeld, {\em Moduli varieties of $F$--sheaves},
  Funct. Anal. Appl. {\bf 21} (1987) 107--122.

\bi[Dr4]{Dr4} V.G. Drinfeld, {\em The proof of Petersson's conjecture
  for $GL(2)$ over a global field of characteristic $p$},
  Funct. Anal. Appl. {\bf 22} (1988) 28--43.

\bibitem[DS]{DS} V. Drinfeld and V. Sokolov, {\em Lie algebras and KdV
type equations}, J. Sov. Math. {\bf 30} (1985) 1975--2036.

\bi[DrSi]{DSimp} V. Drinfeld and C. Simpson, {\em $B$--structures on
$G$--bundles and local triviality}, Math. Res. Lett. {\bf 2} (1995)
823--829.

\bi[EF]{EF} D. Eisenbud and E. Frenkel, Appendix to M. Mustata, {\em
Jet schemes of locally complete intersection canonical singularities},
Invent. Math. {\bf 145} (2001) 397--424.

\bi[FF1]{FF:usp} B. Feigin and E. Frenkel, {\em A family of
  representations of affine Lie algebras}, Russ. Math. Surv.  {\bf
  43}, N 5 (1988) 221--222.

\bi[FF2]{FF:si} B. Feigin and E. Frenkel, {\em Affine Kac-Moody
  Algebras and semi-infinite flag manifolds}, Comm.  Math. Phys. {\bf
  128}, 161--189 (1990).

\bibitem[FF3]{FF:gd} B. Feigin and E. Frenkel, {\em Affine Kac--Moody
algebras at the critical level and Gelfand--Dikii algebras},
Int. Jour. Mod. Phys. {\bf A7}, Supplement 1A (1992) 197--215.

\bi[FK]{Weil} E. Freitag and R. Kiehl, {\em Etale Cohomology and the
  Weil conjecture}, Springer, 1988.

\bibitem[F1]{F:icmp}
E. Frenkel, {\em Affine algebras, Langlands duality
and Bethe Ansatz}, in Proceedings of the International Congress of
Mathematical Physics, Paris, 1994, ed. D. Iagolnitzer, pp. 606--642,
International Press, 1995; arXiv: q-alg/9506003.

\bi[F2]{F:bull} E. Frenkel, {\em Recent advances in the Langlands
   Program}, Bull. Amer. Math. Soc. {\bf 41} (2004) 151--184
   (math.AG/0303074).

\bi[F3]{F:wak} E. Frenkel, {\em Wakimoto modules, opers and the center
    at the critical level}, Advances in Math. {\bf 195} (2005)
    297--404.

\bi[F4]{F:faro} E. Frenkel, {\em Gaudin model and opers}, in Infinite
    Dimensional Algebras and Quantum Integrable Systems,
    eds. P. Kulish, e.a., Progress in Math. {\bf 237}, pp. 1--60,
    Birkh\"auser, 2005.

\bi[F5]{F:flag} E. Frenkel, {\em Opers on the projective line, flag
manifolds and Bethe Ansatz}, Mosc. Math. J. {\bf 4} (2004) 655--705.

\bi[F6]{F:rev} E. Frenkel, {\em Lectures on the Langlands Program and
conformal field theory}, Preprint hep-th/0512172.

\bi[F7]{newbook} E. Frenkel, {\em Langlands Correspondence for Loop
  Groups. An Introduction}, to be published by Cambridge University
  Press; draft available at http://math.berkeley.edu/$\sim$frenkel/

\bi[FB]{vertex} E. Frenkel and D. Ben-Zvi, {\em Vertex algebras and
algebraic curves}, Second Edition, Mathematical Surveys and
Monographs, vol. 88. AMS 2004.

\bibitem[FG1]{FG:exact} E.~Frenkel and D.~Gaitsgory, {\em D-modules on
the affine Grassmannian and representations of affine Kac-Moody
algebras}, Duke Math. J. {\bf 125} (2004) 279--327.

\bi[FG2]{FG:local} E. Frenkel and D. Gaitsgory, {\em Local geometric
  Langlands correspondence and affine Kac-Moody algebras}, Preprint
  math.RT/0508382.

\bi[FG3]{FG:fusion} E. Frenkel and D. Gaitsgory, {\em Fusion and
convolution: applications to affine Kac-Moody algebras at the critical
level}, Preprint math.RT/0511284.

\bi[FG4]{FG:equiv} E. Frenkel and D. Gaitsgory, {\em Localization of
$\ghat$-modules on the affine Grassmannian}, Preprint 
math.RT/0512562.

\bi[FG5]{FG:wak} E. Frenkel and D. Gaitsgory, {\em Geometric realizations
of Wakimoto modules at the critical level}, Preprint math.RT/0603524.

\bi[FG6]{FG:weyl} E. Frenkel and D. Gaitsgory, {\em Weyl modules and
  opers without monodromy}, to appear.

\bibitem[FGV]{FGV1} E. Frenkel, D. Gaitsgory and K. Vilonen, {\em On
    the geometric Langlands conjecture}, Journal of AMS {\bf 15}
    (2001) 367--417.

\bibitem[FT]{FT} E. Frenkel and C. Teleman, {\em Self-extensions of
Verma modules and differential forms on opers}, Compositio Math. {\bf
  142} (2006) 477--500.

\bibitem[Ga1]{Ga:central} D.~Gaitsgory, {\em Construction of central
elements in the Iwahori Hecke algebra via nearby cycles},
Inv. Math. {\bf 144} (2001) 253--280.

\bi[Ga2]{Ga:vanish} D. Gaitsgory, {\em On a vanishing conjecture
appearing in the geometric Langlands correspondence}, Ann. Math. {\bf
160} (2004) 617--682.

\bibitem[Ga3]{Ga} D.~Gaitsgory, {\em The notion of category over an
algebraic stack}, Preprint math.AG/0507192.

\bi[GM]{GM} S.I. Gelfand, Yu.I. Manin, {\em Homological algebra},
Encyclopedia of Mathematical Sciences {\bf 38}, Springer, 1994.

\bi[GW]{GW} S. Gukov and E. Witten, {\em Gauge theory, ramification,
  and the geometric Langlands Program}, to appear.

\bi[HT]{HT} M. Harris and R. Taylor, {\em The geometry and cohomology
of some simple Shimura varieties}, Annals of Mathematics Studies {\bf
151}, Princeton University Press, 2001.

\bi[He]{Henniart} G. Henniart, {\em Une preuve simple des conjectures
de Langlands pour ${\rm GL}(n)$ sur un corps $p$-adique},
Invent. Math. {\bf 139} (2000) 439--455.

\bi[K1]{Kac:laplace} V. Kac, {\em Laplace operators of
infinite-dimensional Lie algebras and theta functions},
Proc. Nat. Acad. Sci. U.S.A. {\bf 81} (1984) no. 2, Phys. Sci.,
645--647.

\bi[K2]{Kac} V.G. Kac, {\em Infinite-dimensional Lie Algebras}, 3rd
Edition, Cambridge University Press, 1990.

\bi[KW]{KW} A. Kapustin and E. Witten, {\em Electric-magnetic duality
  and the geometric Langlands Program}, Preprint hep-th/0604151.

\bi[KL]{KL} D. Kazhdan and G. Lusztig, {\em Proof of the
  Deligne--Langlands conjecture for Hecke algebras}, Inv. Math. {\bf
  87} (1987) 153--215.

\bi[Ku]{Kudla} S. Kudla, {\em The local Langlands correspondence: the
non-Archimedean case}, in {\em Motives} (Seattle, 1991), pp. 365--391,
Proc. Sympos. Pure Math. {\bf 55}, Part 2, AMS, 1994.

\bi[Laf]{Laf} L. Lafforgue, {\em Chtoucas de Drinfeld et
   correspondance de Langlands}, Invent. Math. {\bf 147} (2002)
   1--241.

\bi[L]{Langlands} R.P.~Langlands, {\em Problems in the theory of
automorphic forms}, in Lect. Notes in Math. {\bf 170}, pp. 18--61,
Springer Verlag, 1970.

\bi[La]{Laumon} G. Laumon, {\em Transformation de Fourier,
  constantes d'\'equations fonctionelles et conjecture de Weil},
  Publ. IHES {\bf 65} (1987) 131--210.

\bi[LRS]{LRS} G. Laumon, M. Rapoport and U. Stuhler, {\em ${\mc
D}$-elliptic sheaves and the Langlands
correspondence}, Invent. Math. {\bf 113} (1993) 217--338.

\bi[Lu1]{Lus} G. Lusztig, {\em Classification of unipotent
  representations of simple $p$-adic groups}, Int. Math. Res. Notices
  (1995) no. 11, 517--589.

\bi[Lu2]{Lus1} G. Lusztig, {\em Bases in K-theory}, Represent. Theory
{\bf 2} (1998) 298--369; {\bf 3} (1999) 281--353.

\bi[Mi]{Milne} J.S. Milne, {\em \'Etale cohomology}, Princeton
University Press, 1980.

\bi[MV]{MV} I. Mirkovi\'c and K. Vilonen, {\em Geometric Langlands
duality and representations of algebraic groups over commutative
rings}, Preprint math.RT/0401222.

\bi[Sat]{Satake} I. Satake, {\em Theory of spherical functions on
reductive algebraic groups over $p$--adic fields}, IHES
Publ. Math. {\bf 18} (1963) 5--69.

\bibitem[Vog]{Vogan} D.A. Vogan, {\em The local Langlands conjecture},
  Contemporary Math. {\bf 145}, pp. 305--379, AMS, 1993.

\bibitem[Wak]{Wak} M. Wakimoto, {\em Fock representations of affine Lie
algebra $A_1^{(1)}$}, Comm. Math. Phys. {\bf 104} (1986) 605--609.

\end{thebibliography}
\end{document}